\documentclass[
	a4paper,
	11pt,
]{article}
\usepackage{amsmath}
\usepackage{amssymb}
\usepackage{amsthm}
\usepackage{MnSymbol}
\usepackage{dsfont}
\usepackage{paralist}
\usepackage{enumitem}
\usepackage{hyperref,url}
\usepackage{bm} 
\usepackage{placeins}
\usepackage{color}
\usepackage{graphicx,subfig} 
\usepackage[margin=30mm]{geometry}
\usepackage[font=small]{caption}

\usepackage{array}
\newcolumntype{C}[1]{>{\centering\arraybackslash}p{#1}}

\definecolor{FrameColor}{rgb}{0.85,0.85,0.85}

\newcommand{\EOC}{\mathrm{EOC}}
\newcommand{\V}{\mathcal{V}}

\newcommand{\NN}{\mathcal{N}}
\newcommand{\NNG}{\mathcal{N}_\Gamma}
\newcommand{\N}{\mathbb{N}}
\newcommand{\R}{\mathbb{R}}

\newcommand{\CL}{\mathcal{L}}
\renewcommand{\div}{ \mathrm{div}  }

\newcommand{\dG}{\, d\Gamma}

\newcommand{\ds}{\, ds}
\newcommand{\dx}{\, dx}
\newcommand{\dt}{\, dt}

\newcommand{\pd}{\partial}
\newcommand{\bn}{\bm{n}}
\newcommand{\pdnu}{\pd_{\bm{n}}}
\newcommand{\abs}[1]{\left| #1 \right|}
\newcommand{\norm}[1]{\| #1 \|}
\newcommand{\inn}[2]{ \langle #1 , #2  \rangle}

\newcommand{\mean}[1]{\langle #1 \rangle}
\newcommand{\eps}{\varepsilon}
\newcommand{\Lx}{\Delta}

\newcommand{\LB}{\Delta_{\Gamma}}
\newcommand{\muG}{\mu_\Gamma}
\newcommand{\K}[1]{#1^{K}}

\newcommand{\Th}{\mathcal{T}_h}

\newcommand{\grad}{\nabla}
\newcommand{\gradg}{\nabla_\Gamma}
\newcommand{\mo}{m_\Omega}
\newcommand{\mg}{m_\Gamma}

\newcommand{\muGN}{\mu_{\Gamma,N}}

\newcommand{\bu}{\bm{u}}

\newcommand{\feta}{\bm{\eta}}
\newcommand{\buN}{\bm{u}_N}

\newcommand{\lu}{\bar u}
\newcommand{\lv}{\bar v}
\newcommand{\lbu}{\bar \bu}
\newcommand{\lmu}{\bar \mu}
\newcommand{\lmuG}{{\bar\mu}_{\Gamma}}
\newcommand{\lmuGN}{{\bar\mu}_{\Gamma,N}}

\newcommand{\Vk}{{\mathcal{V}^\kappa}}
\newcommand{\Vm}{{\mathcal{V}^\kappa_m}}
\newcommand{\Vo}{{\mathcal{V}^\kappa_0}}
\newcommand{\Vod}{{(\Vo)'}}

\newcommand{\intO}{\int_\Omega}
\newcommand{\intG}{\int_\Gamma}

\newcommand{\intS}{\int_{\Sigma_T}}

\renewcommand{\H}{\mathcal{H}}
\newcommand{\Hm}{{\mathcal{H}_m}}
\newcommand{\Ho}{{\mathcal{H}_0}}
\newcommand{\Hod}{{(\Ho)'}}

\newcommand{\Xk}{\mathcal{X}^\kappa}
\newcommand{\Yk}{\mathcal{Y}^\kappa}
\newcommand{\Zmk}{\mathcal{Z}_m^\kappa}

\newcommand{\wto}{\rightharpoonup}
\newcommand{\wsto}{\overset{\ast}{\rightharpoonup}}
\newcommand{\emb}{\hookrightarrow}

\newcommand{\PM}{\mathcal{P}_M}

\newcommand{\bl}{\color{blue}}

\newcommand{\CG}{\mathcal{G}}

\theoremstyle{plain}
\newtheorem{thm}{Theorem}[section]

\theoremstyle{plain}
\newtheorem{remark}{Remark}[section]

\numberwithin{equation}{section}


\makeatletter
\renewcommand\paragraph{\@startsection{paragraph}{4}{\z@}%
	{1ex \@plus1ex \@minus.2ex}%
	{-1em}%
	{\normalfont\normalsize\bfseries}}
\renewcommand\subparagraph{\@startsection{paragraph}{4}{\z@}%
	{1ex \@plus1ex \@minus.2ex}%
	{-1em}%
	{\normalfont\normalsize\bfseries}}
\makeatother

\usepackage{hyperref} 
\makeatletter
\hypersetup{%
	colorlinks	=true,% Definition der Links im PDF File
	linkcolor	=blue,%
	citecolor	=red,%
	urlcolor	=blue
}
\makeatother

\begin{document}

\title{Convergence of a Robin boundary approximation for a Cahn--Hilliard system with dynamic boundary conditions}

\author{Patrik Knopf \footnotemark[1] \and Kei Fong Lam \footnotemark[2]}

\date{ }

\renewcommand{\thefootnote}{\fnsymbol{footnote}}
\footnotetext[1]{Department of Mathematics, University of Regensburg, 93053, Germany \tt(Patrik.Knopf@ur.de)}
\footnotetext[2]{Department of Mathematics, The Chinese University of Hong Kong, Shatin, N.T., Hong Kong \tt (kflam@math.cuhk.edu.hk)}

\maketitle

\begin{center}
	\textit{This is a preprint version of the paper. Please cite as:} \\  
	P. Knopf and K.F. Lam,
	Nonlinearity 33(8): 4191-4235, 2020  \\ 
	\url{https://doi.org/10.1088/1361-6544/ab8351}
\end{center}

\bigskip

\begin{abstract}
We prove the existence of unique weak solutions to an extension of a Cahn--Hilliard model proposed recently by C.~Liu and H.~Wu (2019), in which the new dynamic boundary condition is further generalised with an affine linear relation between the surface and bulk phase field variables.  As a first approach to tackle more general and nonlinear relations, we investigate the existence of unique weak solutions to a regularisation by a Robin boundary condition.  Included in our analysis is the case where there is no diffusion for the surface phase field, which causes new difficulties for the analysis of the Robin system.  Furthermore, for the case of affine linear relations, we show the weak convergence of solutions as the regularisation parameter tends to zero, and derive an error estimate between the two models.  This is supported by numerical experiments which also demonstrate some non-trivial dynamics for the extended Liu--Wu model that is not present in the original model. 
\end{abstract}



\section{Introduction}
Since its introduction by physicists \cite{Fis1,Fis2,Kenzler}, the Cahn--Hilliard equation with dynamic boundary conditions has been the subject of many analytical and numerical studies, see, for instance, \cite{LW,Miran} for a recent overview.  The modification from the standard zero Neumann boundary conditions to dynamic boundary conditions allows for the inclusion of effective short-range interactions between the boundary of the domain and the mixture contained within.  One application in which these interactions cannot be neglected lies in the modelling of moving contact lines in binary fluid mixtures \cite{GGM,GGW,Qian}, where instead of the fixed contact angle of $\frac{\pi}{2}$ imposed by the typical Neumann conditions for the Cahn--Hilliard component, it has been proposed to employ a dynamic boundary condition that allows some transport on the domain boundary to better emulate the dynamics of the moving contact line.

The term {\it dynamic boundary conditions} in fact encapsulates a class of boundary conditions, exhibiting a common structure involving a surface time-dependent partial differential equation (PDE) containing the normal derivative of the bulk quantity.  Many such dynamic boundary conditions can be derived as the variation of suitable surface free energies.  To fix ideas, let $T$ be an arbitrary positive real number and let $\Omega \subset \R^3$ be a bounded domain with sufficiently smooth boundary $\Gamma$.   We denote the normal derivative as $\pdnu f := \nabla f \cdot \bn$ with the unit outer normal vector $\bn$ on $\Gamma$.  For fixed $\eps > 0$, the bulk Ginzburg--Landau free energy functional is defined as
\begin{align*}
E_b(u) = \int_\Omega \frac{\eps}{2} \abs{\nabla u}^2 + \frac{1}{\eps} F(u) \dx
\end{align*}
where $F$ is a double well potential with two equal minima and $u$ denotes the phase field variable.  The Cahn--Hilliard equation is the $H^{-1}$ gradient flow of the above free energy:
\begin{align*}
\begin{cases}
u_t = \Lx \mu, \quad \mu = - \eps \Lx u + \eps^{-1} F'(u) & \text{ in } \Omega, \\
\pdnu u = \pdnu \mu =0 & \text{ on } \Gamma,
\end{cases}
\end{align*}
typically furnished with an initial condition $u(0) = u_0$ in $\Omega$.  The condition $\pdnu u = 0$ arises naturally from taking the variation of $E_b$, while the condition $\pdnu \mu = 0$ is chosen to guarantee conservation of mass.  To account for interactions with the domain boundary $\Gamma$, we prescribe a surface free energy functional
\begin{align*}
E_s(u) = \int_\Gamma \frac{\kappa}{2} \abs{\gradg u}^2 + G(u) \dG,
\end{align*}
where $\kappa \geq 0$ is a fixed constant, $G$ is a surface potential function, and $\gradg$ denotes the surface gradient on $\Gamma$.  By performing variation of $E_s$ with respect to suitable inner products, various dynamic boundary conditions replacing $\pdnu u = 0$ can be obtained.  For instance, the dynamic boundary condition of Allen--Cahn type
\begin{align*}
u_t = - \pdnu u + \kappa \LB u - G'(u) \text{ on } \Gamma
\end{align*}
proposed in \cite{Fis1,Fis2,Kenzler} can be regarded formally as the $L^2$-gradient flow of $E_s$, where $\LB$ is the Laplace--Beltrami operator and the bulk effects are captured by the term $\pdnu u$.  Another dynamic boundary condition, that was proposed in \cite{GMS} and studied in \cite{CherfilsGM,ColliF}, is of Cahn--Hilliard type,
\begin{align*}
	u_t = \sigma \LB \mu - \pdnu \mu, \quad \mu = - \kappa \LB u + G'(u) + \pdnu u \text{ on } \Gamma,
\end{align*}
for nonnegative constants $\sigma, \kappa \geq 0$, and arises from taking the subdifferential of the total energy $E = E_b + E_s$ with respect to a suitable scalar product.  In the recent work \cite{LW}, Liu and Wu employed the energetic variational approach to derive a new type of dynamic boundary condition for the Cahn--Hilliard equation.  This approach combines the least action principle and Onsager's principle of maximum energy dissipation to ensure that the resulting model fulfils the physical constraints of mass conservation, dissipation of energy and balance of forces.  Our present interest stems from the observation that we can extend the model and derivation of Liu and Wu to the following Cahn--Hilliard system with dynamic boundary conditions: 
\begin{subequations}\label{CHLW:lim}
\begin{alignat}{4}
&u_t = \Lx \mu, \quad &&\mu = -  \eps \Lx u + \eps^{-1} F'(u) && \text{ in } Q_T := \Omega \times (0,T), \\
&v_t = \LB \muG, \quad  && \muG = - \delta \kappa \LB v + \delta^{-1} G'(v) + \eps \alpha \pdnu u \ && \text{ on } \Sigma_T := \Gamma \times (0,T), \\
&u \vert_{\Sigma_T} = \alpha v + \beta, \quad  &&\pdnu \mu = 0 && \text{ on } \Sigma_T, \\
&u(0) = u_0 \text{ in } \Omega, \quad  && v(0) = \alpha^{-1} (u_0 \vert_\Gamma - \beta) && \text{ on } \Gamma.
\end{alignat}
\end{subequations}
Here, $\kappa \geq 0$, $\eps, \delta > 0$, $\alpha \neq 0$ and $\beta \in \R$ are fixed constants, $u$ and $v$ are the bulk and surface phase field variables, and $\mu$ and $\muG$ are the corresponding bulk and surface chemical potentials.  For applications, we may consider $u$ as the difference in volume fractions between two species of materials in the bulk domain, and $v$ as the difference in volume fractions between two species of materials restricted to the domain boundary.  Introducing the variable $\phi = \alpha \muG$ allows us to express \eqref{CHLW:lim} equivalently as
\begin{subequations}\label{lim:alt}
\begin{alignat}{4}
&u_t = \Lx \mu, \quad & &\mu = - \eps \Lx u + \eps^{-1} F'(u) & \text{ in } Q_T, \\
&u_t = \LB \phi,  \quad & &\pdnu \mu = 0, \quad \phi = - \delta \kappa  \LB u + \alpha \delta^{-1} G'(\alpha^{-1}(u - \beta)) + \eps \alpha^2 \pdnu u \;\;  & \text{ on } \Sigma_T,
\end{alignat}
\end{subequations}
which coincides with the model of \cite{CFW,LW} when $\eps = \delta = \alpha = 1$ and $\beta = 0$.  Our consideration of the affine linear transmission condition $u \vert_{\Sigma_T} = \alpha v + \beta$ as oppose to the classical relation $u \vert_{\Sigma_T} = v$ is motivated in part by observations of some non-trivial dynamics that is not present in the original model with $\alpha = 1$ and $\beta = 0$, see Section \ref{sec:num} below.   For example, if $\alpha = -1$ and $\beta = 0$, then the surface phase $v$ is the opposite of the bulk phase $u$. If now $F$ and $G$ have equal minima at $\pm 1$, we immediately see that $(u,v) = (\pm 1, \mp 1)$ is a solution to \eqref{lim:alt} which cannot be achieved for the usual dynamic boundary condition.

Associated to \eqref{lim:alt} is the energy functional
\begin{align*}
E_*(u) := \intO \frac{\eps}{2} \abs{\grad u}^2 + \frac{1}{\eps} F(u) \dx + \intG \frac{\delta \kappa}{2 \alpha^2} \abs{ \gradg u}^2 + \frac{1}{\delta} G(\alpha^{-1}(u - \beta)) \dG.
\end{align*}
In the appendix, we outline the derivation of \eqref{lim:alt} by the energetic variational approach with the help of $E_*$. In Section \ref{sec:GF}, we show that \eqref{lim:alt} (and hence \eqref{CHLW:lim}) can be interpreted as the gradient flow of $E_*$ with respect to a suitable inner product.  This allows us to employ a natural implicit time discretisation to prove, akin to the first author's previous work \cite{GK}, the existence of a unique global weak solution to \eqref{CHLW:lim}.

A further extension would be to replace the classical relation $u \vert_{\Sigma_T} = v$ with the more general situation $u \vert_{\Sigma_T} = H(v)$ for some continuous function $H: \R \to \R$, and we attempt to tackle this through the well known method of penalising this Dirichlet-like condition by a Robin boundary approximation.  More precisely, for $K > 0$, we consider the system
\begin{subequations}\label{CHLW}
\begin{alignat}{4}
&\K u_t = \Lx \K \mu, \quad &&\K \mu = - \eps \Lx \K u + \eps^{-1} F'(\K u) && \text{ in } Q_T, \\
&\K v_t = \LB \K  \muG, \quad &&\K \muG = - \delta \kappa \LB \K v + \delta^{-1} G'(\K v) + \eps H'(\K v) \pdnu \K  u \; && \text{ on } \Sigma_T, \\
&\eps K \pdnu \K u = H(\K v) - \K u, \quad &&\pdnu \K \mu = 0 && \text{ on } \Sigma_T, \\
&\K u(0) = \K u_0 \qquad \text{ in } \Omega, \quad &&\K v(0) = \K v_0 && \text{ on } \Gamma.
\end{alignat}
\end{subequations}
Formally, as $K \to 0$, we recover the desired relation $u \vert_{\Sigma_T} = H(v)$.  In the case $H(s) = \alpha s + \beta$ with $\alpha \neq 0$ and $\beta \in \R$, this strategy has been successful for the Allen--Cahn system with rigorous convergence results and error estimates derived in \cite{CFL}.  In this work we aim to deduce analogous assertions for the more difficult Cahn--Hilliard model \eqref{CHLW}.  

Associated to \eqref{CHLW} is the energy functional
\begin{align*}
E(\bu) := \intO \frac{\eps}{2} \abs{\grad u}^2 + \frac{1}{\eps} F(u) \dx + \intG \frac{\delta \kappa}{2} \abs{\gradg v}^2 + \frac{1}{\delta} G(v) + \frac{1}{2K} \abs{H(v) - u}^2 \dG,
\end{align*}
for $\bu = (u,v)$, and likewise, \eqref{CHLW} can be interpreted as the gradient flow of $E$ with respect to a suitable inner product.  Thus, a natural implicit time discretisation is used to show the existence of a unique global weak solution to \eqref{CHLW} for a large class of nonlinear relation $H$, see \eqref{ass:h}. While one may expect that, as a regularisation, solutions to \eqref{CHLW} exhibit better regularity properties than solutions to the limit system \eqref{CHLW:lim}, in fact the opposite is true. This means we encounter more difficulties, especially in the case $\kappa = 0$, than in the study of the limit model \eqref{CHLW:lim}.  The chief reason is that the surface variable $v$ is not just given by the trace of the bulk variable $u$, but rather coupled through the Robin boundary condition $\eps K \pdnu u = H(v) - u$. Therefore, the regularity of $v$ cannot be deduced by the trace theorem with the help of $u$ as it was done in \cite{GK,LW}.  We have to apply further techniques and restrictions to overcome these difficulties.

Aside from existence and uniqueness of global weak solutions to \eqref{CHLW:lim} and \eqref{CHLW}, we investigate the rigorous convergence as $K \to 0$ of solutions for the case $H(s) = \alpha s + \beta$ and we derive an error estimate of the form (under the same initial data for $\kappa > 0$)
\begin{equation}\label{intro:err}
\begin{aligned}
& \norm{\K u - u }_{L^4(0,T;L^2(\Omega))}^2 + \norm{\K u - u }_{L^2(0,T;H^1(\Omega))}^2 + \norm{\K v - v }_{L^4(0,T;L^2(\Gamma))}^2 \\
& \quad  + \norm{\K v - v }_{L^2(0,T;H^1(\Gamma))}^2 + K^{-1} \norm{\K u - (\alpha \K v + \beta)}_{L^2(\Sigma_T)}^2 \leq C K \norm{\pdnu u}_{L^2(\Sigma_T)}^2,
\end{aligned}
\end{equation}
which implies that the affine linear transmission condition $u \vert_{\Sigma_T} = \alpha v + \beta$ is attained from the approximation solutions $(\K u, \K v)$ to \eqref{CHLW} at a linear rate in $K$.  This is supported by numerical experiments performed for selected values of $(\alpha, \beta)$ in the affine linear case. Moreover, we have also observed the same rate of convergence in numerical experiments for two examples of nonlinear relations $H$.

The structure of this paper is as follows: in Section \ref{sec:main} we outline several preliminary results essential for our investigations and present the main results of this work.  The gradient flow structure for both Cahn--Hilliard systems is demonstrated in Section \ref{sec:GF}, while the existence and uniqueness of weak solutions to \eqref{CHLW:lim} and \eqref{CHLW} is established via an implicit time discretisation in Section \ref{sec:lim} and \ref{sec:CHLW}, respectively.  In Section \ref{sec:conv}, we discuss the convergence of \eqref{CHLW} to \eqref{CHLW:lim} for the case $H(s) = \alpha s + \beta$ and in Section \ref{sec:num} we present finite element discretisations of \eqref{lim:alt} and \eqref{CHLW}, numerical evidence in support of the error estimate \eqref{intro:err}, as well as a comparison between the dynamics of solutions to \eqref{lim:alt} for different values of $\alpha$ and $\beta$. 

\section{Preliminaries and Main results}\label{sec:main}
\paragraph{Notation.} Throughout this paper we use the following notation:  For any $1 \leq p \leq \infty$ and $k \geq 0$, the standard Lebesgue and Sobolev spaces defined on $\Omega$ are denoted as $L^p(\Omega)$ and $W^{k,p}(\Omega)$, along with the norms $\norm{\cdot}_{L^p(\Omega)}$ and $\norm{\cdot}_{W^{k,p}(\Omega)}$.  For the case $k = 2$, these spaces become Hilbert spaces and we use the notation $H^k(\Omega) = W^{k,2}(\Omega)$.  A similar notation is used for Lebesgue and Sobolev spaces on $\Gamma$.  For any Banach space $X$, we denote its dual space by $X'$ and the associated duality pairing by $\inn{\cdot}{\cdot}_X$.  If $X$ is a Hilbert space, we denote its inner product by $(\cdot, \cdot)_X$.  We define
\begin{align*}
\mean{u}_\Omega := \begin{cases}
\frac{1}{\abs{\Omega}} \inn{u}{1}_{H^1(\Omega)} & \text{ if } u \in H^1(\Omega)', \\
\frac{1}{\abs{\Omega}} \int_\Omega u \dx & \text{ if } u \in L^1(\Omega)
\end{cases}
\end{align*}
as the spatial mean of $u$, where $\abs{\Omega}$ denotes the $d$-dimensional Lebesgue measure of $\Omega$.  The spatial mean for $v \in H^1(\Gamma)'$ and $v \in L^1(\Gamma)$ can be defined analogously. For brevity, we also use the notation
\begin{align*}
	\CL^p := L^p(\Omega) \times L^p(\Gamma)
	\quad\text{and}\quad
	\H^k := H^k(\Omega) \times H^k(\Gamma)
\end{align*}
for all $p\in [1,\infty]$ and $k\ge 0$.

In the sections 2-6 we set the parameters $\eps = \delta = 1$ in both models \eqref{CHLW:lim} and \eqref{CHLW} as their values have no impact on the analysis we will carry out. 

\bigskip

\paragraph{Assumptions.}
\begin{enumerate}[label=$(\mathrm{A \arabic*})$, ref = $\mathrm{A \arabic*}$]
\item \label{ass:dom} We take $\Omega \subset \R^d$ with $d \in \{2,3\}$ to be a bounded domain with Lipschitz boundary $\Gamma$.  
\item  \label{ass:pot}
We assume that the potentials $F$ and $G$ are bounded from below by constants $C_F,C_G\in\R$ such that
\begin{align*}
F(s)  \ge C_F \quad\text{and}\quad G(s)  \ge C_G \quad \text{for all}\; s \in\R,
\end{align*}	
and exhibit a decomposition in the form
\begin{align*}
F(s)=F_1(s) + F_2(s)  \quad\text{and}\quad G(s)=G_1(s)  + G_2(s)  \quad\text{for all}\;  s \in\R
\end{align*}
satisfying
\begin{enumerate}[label=$(\mathrm{A 2 \roman*})$, ref = $\mathrm{A 2 \roman*}$]
\item \label{ass:pot:1} $F_1, F_2, G_1, G_2 \in C^2(\R)$.
\item \label{ass:pot:2} $F_1$ and $G_1$ are convex nonnegative functions.
\item \label{ass:pot:3} For any $\delta > 0$, there exist positive constants $A_{f,\delta}$ and $A_{g,\delta}$ such that 
\begin{align*}
\abs{F_1'(s)} \leq \delta F_1(s) + A_{f,\delta} \quad \text{ and } \quad 
\abs{G_1'(s)} \leq \delta G_1(s) + A_{g,\delta} \quad \text{ for all } s \in \R.
\end{align*}
\item \label{ass:pot:4}$F_2'$ and $G_2'$ are Lipschitz continuous.
\end{enumerate}
\item \label{ass:h} The function $H$ belongs to $C^2(\R)$ and there exist $q \in [1,\infty)$ and a positive constant $d$ such that for all $s,t \in \R$, 
\begin{align*}
& \abs{H'(s)} \leq d(1 + \abs{s}^{q}), \\
& \abs{H(s) - H(t)} \leq d(1 + \abs{s}^{q} + \abs{t}^{q} ) \abs{s-t}, \quad \abs{H'(s) - H'(t)} \leq d(1 + \abs{s}^{q-1} + \abs{t}^{q-1}) \abs{s-t}.
\end{align*}
\end{enumerate}

\begin{remark}	\label{REM:POT} The assumption \eqref{ass:pot:4} implies there exist positive constants $c_F$ and $c_G$ such that for all $s\in\R$,
\begin{gather*}
\abs{F_2'(s)} \le c_f(1+\abs{s}), \quad \abs{G_2'(s)} \le c_g(1+\abs{s}), \quad \abs{F_2(s)} \le c_F(1+\abs{s}^2), \quad \abs{G_2(s)} \le c_G(1+\abs{s}^2).
\end{gather*}
\end{remark}

\begin{remark}
For our notion of weak solutions to \eqref{CHLW:lim} and \eqref{CHLW}, cf.~\cite[Defn.~4.1]{GK}, it suffices to have a Lipschitz boundary $\Gamma$, where the definition of tangential gradient on a Lipschitz surface can be found in \cite[Defn.~3.1]{Buffa}.  For the derivation of the error estimate \eqref{intro:err}, we will assume a smooth boundary $\Gamma$ in order to employ classical elliptic regularity results, which also is consistent with the set-up of \cite{LW}.
\end{remark}
\bigskip

\paragraph{Preliminaries.} 
\begin{enumerate}[label=$(\mathrm{P \arabic*})$, ref = $\mathrm{P \arabic*}$]
\item \label{pre:Na} For $X = \Omega$ or $\Gamma$, let $H^1(X)'_0 := \{ \phi \in H^1(X)' \, : \, \mean{\phi}_X = 0 \}$ and $H^1(X)_0 := \{ \theta \in H^1(X) \, : \, \mean{\theta}_X = 0 \}$.  Then, the unique solvability of the problems
\begin{align*}
- \LB \theta_\Gamma = \phi_\Gamma \text{ on } \Gamma
, \qquad \begin{cases}
 - \Lx \theta = \phi &\text{ in } \Omega, \\
\pdnu \theta = 0 &\text{ on } \Gamma
\end{cases}
\end{align*}
induces the solution operators $\NNG : H^1(\Gamma)'_0 \to H^1(\Gamma)_0$, $\phi_\Gamma \mapsto \NNG(\phi_\Gamma) = \theta_\Gamma$ and ${\NN : H^1(\Omega)'_0 \to H^1(\Omega)_0}$, $ \phi \mapsto \NN(\phi) = \theta$, satisfying
\begin{align}
\label{SolnOp:b} \norm{\phi}_{L^2(\Omega)}^2 & = \int_\Omega \nabla \phi \cdot \nabla \NN(\phi) \dx \leq \norm{\nabla \phi}_{L^2(\Omega)} \norm{\nabla \NN(\phi)}_{L^2(\Omega)}, \\
\label{SolnOp:s} \norm{\phi_\Gamma}_{L^2(\Gamma)}^2 & = \int_\Gamma \gradg \phi_\Gamma \cdot \gradg \NNG(\phi_\Gamma) \dG \leq \norm{\gradg \phi_\Gamma}_{L^2(\Gamma)} \norm{\gradg \NNG(\phi_\Gamma)}_{L^2(\Gamma)}.
\end{align}
\item \label{pre:H1} In this paper, the dual spaces $H^1(\Omega)'$ and $H^1(\Gamma)'$ are endowed with the norms
\begin{align*}
\norm{\phi}_{H^1(\Omega)'}^2 &:= \norm{\grad \NN(\phi - \mean{\phi}_\Omega)}_{L^2(\Omega)}^2 + \abs{\mean{\phi}_\Omega}^2
&&\text{for all}\; \phi\in H^1(\Omega)',\\
\norm{\phi_\Gamma}_{H^1(\Gamma)'}^2 &:= \norm{\gradg \NN(\phi_\Gamma - \mean{\phi_\Gamma}_\Gamma)}_{L^2(\Gamma)}^2 + \abs{\mean{\phi_\Gamma}_\Gamma}^2
&&\text{for all}\; \phi_\Gamma\in H^1(\Gamma)',
\end{align*}
which are equivalent to the standard norms on these spaces. The spaces $H^1(\Omega)_0'$ and $H^1(\Gamma)_0'$ are endowed with the same norms, i.e.,\begin{align*}
\norm{\phi}_{H^1(\Omega)'} = \norm{\grad \NN(\phi)}_{L^2(\Omega)}, \quad
	\norm{\phi_\Gamma}_{H^1(\Gamma)'} = \norm{\gradg \NNG(\phi_\Gamma)}_{L^2(\Gamma)}
\end{align*}
for $\phi\in H^1(\Omega)_0'$ and $\phi_\Gamma \in H^1(\Gamma)_0'$.
\item \label{pre:Nb} 
Suppose that $\phi \in L^2(Q_T) \cap H^1(0,T;H^1(\Omega)'_0)$ along with solutions $\theta := \NN(\phi)$ and $\psi := \NN(\phi_t)$.  Then, it holds that $\theta_t = \psi$ in the sense of distributions, which can be inferred by testing the Neumann problem for $\theta$ with a test function $q_t \zeta$, where $q \in C^{\infty}_c(0,T)$ and $\zeta \in H^1(\Omega)$ are arbitrary, integrating by parts in time and employing the uniqueness of solutions.  With this characterisation one also arrives at
\begin{align}\label{SolnOp:u:time}
\frac{d}{dt} \frac{1}{2} \norm{\nabla \NN(\phi)}_{L^2(\Omega)}^2 = \int_\Omega \nabla \theta \cdot \nabla \psi \dx = \inn{\phi_t}{\theta}_{H^1(\Omega)} = \inn{\phi_t}{\NN(\phi)}_{H^1(\Omega)},
\end{align}
and analogously, for $\phi_\Gamma \in L^2(\Sigma_T) \cap H^1(0,T;H^1(\Gamma)'_0)$ with $\mean{\phi_\Gamma}_\Gamma = 0$, the solutions $\theta_\Gamma := \NNG(\phi_\Gamma)$ and $\psi_\Gamma := \NNG(\phi_{\Gamma,t})$ satisfy $\theta_{\Gamma,t} = \psi_\Gamma$ in the sense of distributions and the equality
\begin{align}\label{SolnOp:v:time}
\frac{d}{dt} \frac{1}{2} \norm{\gradg \NNG(\phi_\Gamma)}_{L^2(\Gamma)}^2 = \inn{\phi_{\Gamma,t}}{\NNG(\phi_\Gamma)}_{H^1(\Gamma)}.
\end{align}
\item \label{pre:H} For any $m=(\mo,\mg)\in\R^2$ we define
\begin{align*}
\Hm := \{ \bu=(u,v) \in \H^1 \,\big\vert\, \mean{u}_\Omega = \mo,\; \mean{v}_\Gamma = \mg \} , \quad \Ho := \H_{(0,0)} .
\end{align*}
For $\bu=(u,v) \in \Hod$, there exist $\theta_{\Omega,u} \in H^1(\Omega)_0$ and $\theta_{\Gamma,v} \in H^1(\Gamma)_0$ such that
\begin{align*}
\bu(\zeta,\xi) = \intO \grad\theta_{\Omega,u} \cdot \grad\zeta \dx + \intG \gradg\theta_{\Gamma,v} \cdot \gradg \xi \dG \quad\text{for all}\; (\zeta,\xi)\in\Ho.
\end{align*}
Hence, we can define the inner product
\begin{align}\label{defn:H0':inn}
(\bu_1,\bu_2)_\Hod &:= \intO \grad\theta_{\Omega,u_1} \cdot \grad\theta_{\Omega,u_2} \dx + \intG \gradg\theta_{\Gamma,v_1} \cdot\gradg\theta_{\Gamma,v_2} \dG
\end{align}
for functions $\bu_1=(u_1,v_1), \bu_2=(u_2,v_2) \in \Hod$. The induced norm is given by $\norm{\bu}_\Hod := (\bu,\bu)_\Hod^{1/2}$. If $\bu=(u,v)$ actually lies in $\Ho\subset\Hod$, the functions $\theta_{\Omega,u}$ and $\theta_{\Gamma,v}$ can be expressed by
\begin{align*}
\theta_{\Omega,u} = \NN(u) \quad\text{and}\quad \theta_{\Gamma,v} = \NNG(v).
\end{align*}
We remark that $\norm{\,\cdot\,}_\Hod$ is also a norm on $\Ho$. However, $\Ho$ is not complete with respect to this norm.
\item \label{pre:V} For fixed $\kappa \geq 0$ we define the Hilbert space
\begin{align*}
\Vk & := \begin{cases}
\{ \theta \in H^1(\Omega) \, : \, \theta \vert_\Gamma \in H^1(\Gamma) \}, & \kappa > 0, \\
H^1(\Omega), & \kappa = 0,
\end{cases}
\end{align*}
endowed with the inner product and induced norm
\begin{align*}
(\phi, \psi)_{\V^{\kappa}} & := \begin{cases}
(\phi, \psi)_{H^1(\Omega)} + (\phi \vert_\Gamma, \psi \vert_\Gamma)_{H^1(\Gamma)}, & \kappa > 0, \\
(\phi, \psi)_{H^1(\Omega)}, & \kappa = 0,
\end{cases}
\quad \norm{\phi}_{\V^{\kappa}} := (\phi, \phi)_{\V^{\kappa}}^{1/2}.
\end{align*}
For any $m=(\mo,\mg)\in\R^2$, we also define the subsets
\begin{align*}
\Vm := \{ \theta\in\Vk \,:\, \mean{\theta}_\Omega = \mo,\; \mean{\theta}_\Gamma = \mg \}, \quad \Vo := \V^\kappa_{(0,0)}.
\end{align*}
For $\phi \in \Vod$, there exist $\theta_{\Omega,\phi} \in H^1(\Omega)_0$ and $\theta_{\Gamma,\phi} \in H^1(\Gamma)_0$ such that
\begin{align*}
	\phi(\zeta) = \intO \grad\theta_{\Omega,\phi} \cdot \grad\zeta \dx + \intG \gradg\theta_{\Gamma,\phi} \cdot \gradg \zeta \dG \quad\text{for all}\; \zeta\in\Vo.
\end{align*}
Thus, we can define an inner product on $\Vod$ by
\begin{align*}
(\phi,\psi)_\Vod = \intO \grad\theta_{\Omega,\phi} \cdot \grad\theta_{\Omega,\psi} \dx + \intG \gradg\theta_{\Gamma,\phi} \cdot \gradg\theta_{\Gamma,\psi} \dG \quad\text{for all}\; \phi,\psi\in\Vod.
\end{align*}
Its induced norm is given by $\norm{\phi}_\Vod = (\phi,\phi)_\Vod^{1/2}$. If $\phi$ lies in $\Vo \subset \Vod$, the functions $\theta_{\Omega,\phi}$ and $\theta_{\Gamma,\phi}$ can be expressed by
\begin{align*}
\theta_{\Omega,\phi} = \NN(\phi) \quad\text{and}\quad \theta_{\Gamma,\phi} = \NNG(\phi).
\end{align*}
It holds that $\norm{\,\cdot\,}_\Vod$ is also a norm on $\Vo$ but $\Vo$ is not complete with respect to this norm.
\item \label{pre:trace} A closer examination of the proof of the trace theorem (cf.\,\cite{Evans}) shows that there exists a constant $C$ depending only on $\Omega$ and $p \in [1,\infty)$ such that for all $f\in W^{1,p}(\Omega)$, 
\begin{align*}
\norm{f}_{L^p(\Gamma)}^p \leq C \big (\norm{f}_{L^p(\Omega)}^p + \norm{\nabla f}_{L^p(\Omega)} \norm{f}_{L^p(\Omega)}^{p-1} \big ).
\end{align*}
\end{enumerate}

\bigskip

The main results of this paper are listed as follows.

\begin{thm}[Global existence and uniqueness for \eqref{lim:alt}]\label{thm:lim}
Let $T > 0$, $\kappa \geq 0$, $\alpha \neq 0$, $\beta \in \R$ and $m = (m_\Omega, m_\Gamma) \in \R^2$ be arbitrary and suppose that \eqref{ass:dom} and \eqref{ass:pot} hold. Then, for any initial condition $u_0 \in \V_m^\kappa$ satisfying $F(u_0) \in L^1(\Omega)$ and $G(\alpha^{-1}(u_0 - \beta)) \in L^1(\Gamma)$, there exists a unique weak solution $(u, \mu, \phi)$ to \eqref{CHLW:lim} in the following sense:
\begin{enumerate}
\item[(i)] The solution has the following regularity
\begin{align*}
u & \in C^{0,\frac{1}{4}}([0,T];L^2(\Omega)) \cap L^\infty(0,T;H^1(\Omega)) \cap H^1(0,T;H^1(\Omega)'), \\
u & \in H^1(0,T;H^1(\Gamma)') \cap \begin{cases}
L^\infty(0,T;H^1(\Gamma)) \cap C^{0,\frac{1}{4}}(0,T;L^2(\Gamma)), & \kappa > 0, \\
L^\infty(0,T;H^{1/2}(\Gamma)) \cap C^{0,\frac{1}{4}}(0,T;H^1(\Gamma)'), & \kappa = 0,
\end{cases} \\
\mu & \in L^2(0,T;H^1(\Omega)), \quad \phi \in L^2(0,T;H^1(\Gamma)).
\end{align*}
\item[(ii)] For all $\zeta \in H^1(\Omega)$, $\theta \in H^1(\Gamma)$ and $\psi \in \V^{\kappa} \cap L^{\infty}(\Omega)$ with $\psi\vert_\Gamma \in L^\infty(\Gamma)$, and for a.e.~$t \in (0,T)$, the following weak formulation 
\begin{align*}
0 & = \inn{u_t}{\zeta}_{H^1(\Omega)} + \int_\Omega \nabla \mu \cdot \nabla \zeta \dx, \\
0 & = \inn{u_t}{\theta}_{H^1(\Gamma)} + \int_\Gamma  \gradg  \phi \cdot \gradg \theta \dG, \\
0 & = \int_\Omega \nabla u \cdot \nabla \psi + (F'(u) - \mu) \psi \dx + \frac{1}{\alpha^2} \int_\Gamma \kappa\gradg u \cdot \gradg \psi +   (\alpha G'  ( \tfrac{u-\beta}{\alpha} )  - \phi) \psi \dG.
\end{align*}
and the initial condition $u(0) = u_0$ are satisfied.
\item[(iii)] For all $t \in [0,T]$, the energy inequality is satisfied 
\begin{align*}
E_*(u(t)) + \frac{1}{2} \int_0^t \norm{\nabla \mu(s)}_{L^2(\Omega)}^2 + \norm{\gradg \muG(s)}_{L^2(\Gamma)}^2 \ds \leq E_*(u_0).
\end{align*}
\end{enumerate}
In addition, if the boundary $\Gamma$ is smooth and there exist positive constants $c_0$ and $c_1$ such that for all $s \in \R$, 
\begin{align}\label{GK:f1}
\begin{cases}
\abs{F_1'(s)} \leq c_0 \Big (1 + \abs{s}^3 \Big), &\text{ if } \kappa > 0, \\
\abs{F_1'(s)} \leq c_0 \Big (1 + \abs{s}^3 \Big), \quad \abs{G_1''(s)} \leq c_1 , & \text{ if } \kappa = 0,
\end{cases}
\end{align}
then for any $\kappa \geq 0$ it holds that $u \in L^2(0,T;H^2(\Omega))$. Moreover, if $\kappa > 0$, then $u\vert_{\Sigma_T} \in L^2(0,T;H^2(\Gamma))$.
\end{thm}

\begin{thm}[Global existence and uniqueness for \eqref{CHLW}]\label{thm:CHLW}
Let $T > 0$, $\kappa \geq 0$, $K > 0$ and $m = (m_\Omega, m_\Gamma) \in \R^2$ be arbitrary.  Suppose that \eqref{ass:dom}-\eqref{ass:h} hold, and in the case $\kappa = 0$, the following assumptions hold additionally 
\begin{enumerate}[label=$(\mathrm{B \arabic*})$, ref = $\mathrm{B \arabic *}$]
\item \label{ass:pot:5:k} $G_2 \equiv 0$ and there exist a real number $p \geq 2$, and constants $a, c > 0$ and $b \geq 0$ such that 
\begin{align*}
a \abs{s}^{p-2} - b \leq \abs{G_1''(s)} \leq c(1 + \abs{s}^{p-2}) \quad \text{ for all } s \in \R.
\end{align*}
\item \label{ass:h:k=0} $H(s) = \alpha s + \beta$ for all $s \in \R$ where $\alpha \neq 0$ and $b \in \R$ are constants.
\end{enumerate}
For the exponent $p$ in \eqref{ass:pot:5:k}, let
\begin{align}
\label{DEF:XY}
\Xk := 
\begin{cases}
L^2(\Gamma) \cap L^p(\Gamma) &\text{ if } \kappa=0,\\
H^1(\Gamma) &\text{ if } \kappa>0,
\end{cases}
\qquad
\Yk := 
\begin{cases}
H^1(\Gamma)' &\text{ if } \kappa=0,\\
L^2(\Gamma) &\text{ if } \kappa>0.
\end{cases}
\end{align}
Then, for any initial condition $\bu_0 = (u_0, v_0) \in H^1(\Omega) \times \Xk$ with $\mean{u_0}_\Omega = m_\Omega$ and $\mean{v_0}_\Gamma = m_\Gamma$ satisfying $F(u_0) \in L^1(\Omega)$ and $G(v_0) \in L^1(\Gamma)$, there exists a unique weak solution $(u, v, \mu, \muG)$ to \eqref{CHLW} in the following sense:
\begin{enumerate}
\item[(i)] The solution has the following regularity 
\begin{align*}
u & \in C^{0,\frac{1}{4}}([0,T];L^2(\Omega)) \cap L^\infty(0,T;H^1(\Omega)) \cap H^1(0,T;H^1(\Omega)'), \\
v & \in C^{0,\frac{1}{4}}([0,T];\Yk) \cap L^{\infty}(0,T;\Xk) \cap H^1(0,T;H^1(\Gamma)'), \\
\mu & \in L^2(0,T;H^1(\Omega)), \quad \muG \in L^2(0,T;H^1(\Gamma)).
\end{align*}
\item[(ii)] For all $\zeta \in H^1(\Omega)$, $\theta \in H^1(\Gamma)$, $\eta \in H^1(\Omega) \cap L^{\infty}(\Omega)$, $\psi \in \Xk \cap L^{\infty}(\Gamma)$, and for a.e.~$t \in (0,T)$, the following weak formulation 
\begin{align}\label{DEF:CHLW1}
\inn{u_t}{\zeta}_{H^1(\Omega)} &= -\intO \grad\mu \cdot \grad\zeta \dx,\\
\label{DEF:CHLW2}
\inn{v_t}{\theta}_{H^1(\Gamma)} &= 	-\intG \gradg \muG \cdot \gradg \theta \dG,\\
\label{DEF:CHLW3}
\intO \grad u\cdot\grad \eta + F'(u)\eta -  \mu\eta \dx &= \intG K^{-1}\big(H(v)-u\big)  \eta  \dG,\\
\label{DEF:CHLW4}	\intG  \kappa \gradg v\cdot \gradg \psi + G'(v)\psi - \muG \psi \dG& = - \intG K^{-1} \big( H(v) - u \big) H'(v) \psi  \dG,
\end{align}
and the initial conditions $u(0)=u_0$ and $v(0)=v_0$ are satisfied.
\item[(iii)] For all $t\in[0,T]$, the energy inequality is satisfied 
\begin{align}\label{DEF:CHLW5}
E(u(t),v(t)) + \frac{1}{2} \int_0^t \|\grad\mu(s)\|_{L^2(\Omega)}^2 + \|\grad\mu_\Gamma(s)\|_{L^2(\Gamma)}^2 \ds \le E(u_0,v_0).
\end{align}
\end{enumerate}
\end{thm}

\begin{thm}[Convergence for $K\to 0$ and error estimates]\label{thm:conv}
Suppose that \eqref{ass:dom}, \eqref{ass:pot} and \eqref{ass:pot:5:k} hold and let $\Xk$ be as defined in \eqref{DEF:XY}. For arbitrary $T > 0$, $K > 0$, $\kappa \geq 0$, $\alpha \neq 0$, $\beta \in \R$ and $m = (m_\Omega, m_\Gamma) \in \R^2$, let $H(s) = \alpha s + \beta$, and $(\K u, \K v, \K \mu, \K \muG)$ denote the weak solution to \eqref{CHLW} on $[0,T]$ with initial data $(\K u_0, \K v_0) \in H^1(\Omega) \times \Xk$ satisfying
\begin{align}\label{conv:weak:ini}
E(\K u_0, \K v_0) \leq C, \quad \mean{\K u_0}_\Omega = m_\Omega, \quad \mean{\K v_0}_\Gamma = m_\Gamma \quad \text{ for all } K > 0.
\end{align}
Then, along a nonrelabelled subsequence, 
\begin{align*}
& \K u_0 \wto u_0 \text{ in } H^1(\Omega), \quad \K v_0 \to v_0:= \alpha^{-1} (u_0 \vert_\Gamma - \beta) \text{ in } L^2(\Gamma), \\
& \K u \wsto u \text{ in } L^{\infty}(0,T;H^1(\Omega)) \cap H^1(0,T;H^1(\Omega)'), \\
& \K v \wsto v:= \alpha^{-1}(u - \beta) \text{ in } L^{\infty}(0,T;\Xk) \cap H^1(0,T;H^1(\Gamma)'), \\
& \K \mu \wto \mu \text{ in } L^2(0,T;H^1(\Gamma)), \quad \alpha \K \muG \wto \phi \text{ in } L^2(0,T;H^1(\Gamma)),
\end{align*}
where $(u, \mu, \phi)$ is the weak solution to \eqref{lim:alt} on $[0,T]$ with initial data $u_0 \in \V_m^\kappa$ in the sense of items (i) and (ii) of Theorem \ref{thm:lim}.

Moreover, if \eqref{GK:f1} holds, then for any $\kappa > 0$ there exists a positive constant $C$ independent of $K$ and the solution variables such that 
\begin{equation}\label{CHLW:rate}
\begin{aligned}
& \sup_{t \in [0,T]} \norm{(\K u - u, \K v - v)(t)}_{(\Ho)'}^2 + \norm{\K u - u}_{L^4(0,T;L^2(\Omega))}^2 + \norm{\K u - u }_{L^2(0,T;H^1(\Omega))}^2  \\
& \qquad + \norm{\K v - v}_{L^4(0,T;L^2(\Gamma))}^2 + \norm{\K v - v }_{L^2(0,T;H^1(\Gamma))}^2 + K^{-1} \norm{\K u - (\alpha \K v + \beta)}_{L^2(\Sigma_T)}^2 \\
& \quad \leq C \Big ( K \norm{\pdnu u}_{L^2(\Sigma_T)}^2 + \norm{(\K u_0 - u_0, \K v_0 - v_0)}_{(\Ho)'}^2\Big ).
\end{aligned}
\end{equation}
For the case $\kappa = 0$, the above estimate \eqref{CHLW:rate} holds with the modification where we remove the $L^4(0,T;L^2(\Gamma))$ and $L^2(0,T;H^1(\Gamma))$ estimates for $\K v - v$. \end{thm}

\section{Gradient flow structure}\label{sec:GF}
In this section we show that the Cahn--Hilliard systems \eqref{CHLW:lim} and \eqref{CHLW} can be expressed as gradient flows of the energy functionals $E_*$ and $E$ with respect to suitable inner products.

\subsection{Limit system \eqref{lim:alt}}
The first variation of $E_*$ with respect to $u$ in the direction $\eta \in \V^{\kappa}_0 \cap L^\infty(\Omega)$ with $\eta\vert_\Gamma\in L^\infty(\Gamma)$, where we recall the subspace $\V^{\kappa}_0$ defined in \eqref{pre:V} with dual space $(\V^{\kappa}_0)'$, reads as
\begin{align*}
\frac{\delta E_*}{\delta u}(u)[\eta] = \int_\Omega \nabla u \cdot \nabla \eta + F'(u) \eta \dx + \int_\Gamma \frac{\kappa}{\alpha^2} \gradg u \cdot \gradg \eta + \frac{1}{\alpha} G'(\alpha^{-1}(u-\beta)) \eta \dG.
\end{align*}
We introduce the solution operator $\NNG^\alpha : H^1(\Gamma)'_0 \to H^1(\Gamma)$ which is defined by
\begin{align*}
- \LB \NNG^\alpha(\phi_\Gamma) = \frac{1}{\alpha} \phi_\Gamma \text{ on } \Gamma,
\end{align*}
for $\phi_\Gamma \in H^1(\Gamma)'_0$.  Then, for any $u, \eta \in \V^{\kappa}_0$, the functions $\NN(u)$, $\NNG^\alpha(u)$, $\NN(\eta)$ and $\NNG^\alpha(\eta)$ are well-defined, and so we can define the inner product
\begin{align*}
(u, \eta)_{(\V^{\kappa}_0)'_*} := \int_\Omega \nabla \NN(u) \cdot \nabla \NN(\eta) \dx + \int_\Gamma \gradg \NNG^\alpha (u) \cdot \gradg \NNG^\alpha(\eta) \dG \quad \text{ for } u, \eta \in \V^{\kappa}_0,
\end{align*}
with norm $\norm{u}_{(\V^{\kappa}_0)'_*} = (u, u)_{(\V^{\kappa}_0)'_*}^{1/2}$.

We seek to identify the gradient $\nabla E_* \vert_{u}$ of $E_*$ at $u$ associated to this inner product, which satisfies $(\nabla E_* \vert_{u}, \eta)_{(\V^{\kappa}_0)'_*} = \frac{\delta E_*}{\delta u}(u) [\eta]$ for all $\eta \in \V^{\kappa}_0 \cap L^\infty(\Omega)$ with $\eta\vert_\Gamma\in L^\infty(\Gamma)$.  In particular, this yields for $\theta = \nabla E_* \vert_{u}$ the identity
\begin{equation}\label{LW:lim:GF:1}
\begin{aligned}
& \int_\Omega \nabla \NN(\theta) \cdot \nabla \NN(\eta) \dx + \int_\Gamma \gradg \NNG^\alpha(\theta) \cdot \gradg \NNG^\alpha(\eta) \dG \\
& \quad = \int_\Omega F'(u) \eta + \nabla u \cdot \nabla \eta  \dx + \int_\Gamma \frac{\kappa}{\alpha^2} \gradg u \cdot \gradg \eta + \frac{1}{\alpha} G'(\alpha^{-1}(u-\beta)) \eta \dG.
\end{aligned}
\end{equation}
Following ideas in \cite[Sec.~4.2]{GK}, for arbitrary $\tilde \eta \in \V^{\kappa} \cap L^\infty(\Omega)$ with $\tilde \eta\vert_\Gamma\in L^\infty(\Gamma)$ and arbitrary nonnegative $\zeta \in C^{\infty}_c(\Omega)$ not identically zero, we define
\begin{align}\label{LW:GF:c1c2}
c_1 = \mean{\tilde \eta}_{\Gamma} = \frac{1}{\abs{\Gamma}} \int_\Gamma \tilde \eta \dG, \quad c_2 = \frac{\mean{\tilde \eta}_\Omega - \mean{\tilde \eta}_\Gamma}{\mean{\zeta}_\Omega}
\end{align}
so that
\begin{align}\label{eta:test}
\eta := \tilde \eta - c_1 - c_2 \zeta \in \V^{\kappa}_0  \cap L^\infty(\Omega) \quad\text{with}\quad \eta\vert_\Gamma \in L^\infty(\Gamma).
\end{align} 
Defining the constants 
\begin{align}
C_\Omega & := \frac{1}{\abs{\Omega} \mean{\zeta}_\Omega} \Big ( \int_\Omega -\NN(\theta) \zeta + \nabla u \cdot \nabla \zeta + F'(u) \zeta \dx \Big ), \label{LW:cOme}\\
 C_\Gamma & :=  \frac{1}{\abs{\Gamma}}  \Big ( \int_\Omega F'(u) \dx + \int_\Gamma \frac{1}{\alpha}G'(\alpha^{-1}(u - \beta)) \dG \Big ) - C_\Omega \frac{\abs{\Omega}}{\abs{\Gamma}}, \label{LW:cGam}
\end{align}
that are independent of $\tilde \eta$ and plugging the choice \eqref{eta:test} into \eqref{LW:lim:GF:1} yields  the following identifications after a short calculation via the fundamental lemma of calculus of variations:
\begin{equation}\label{LW:lim:GF:2}
\begin{aligned}
\mu  := \NN(\theta) + C_\Omega &= F'(u) - \Lx u &&\text{ in } \Omega, \\
\frac{\phi}{\alpha^2}  := \frac{1}{\alpha} \NNG^\alpha(\theta) + C_\Gamma 
&= \frac{1}{\alpha} G'(\alpha^{-1}(u-\beta)) + \pdnu u -\frac{\kappa}{\alpha^2} \LB u &&\text{ on } \Gamma.
\end{aligned}
\end{equation}
Thus, the gradient flow
\begin{align*}
(u_t,\eta)_{(\V^{\kappa}_0)'} = - \frac{\delta E_*}{\delta u}(u)[\eta] = - (\nabla E_* \vert_{u}, \eta)_{(\V^{\kappa}_0)'} \quad \text{ for all } \eta \in \V^{\kappa}_0 \cap L^\infty(\Omega),\,\eta\vert_\Gamma \in L^\infty(\Gamma)
\end{align*}
can now be expressed as
\begin{align*}
& \int_\Omega \NN(u_t) \eta \dx + \int_\Gamma \frac{1}{\alpha} \NNG^\alpha(u_t) \eta \dG =  \int_\Omega \nabla \NN(u_t) \cdot \nabla \NN(\eta) \dx + \int_\Gamma \gradg \NNG^\alpha(u_t) \cdot \gradg \NNG^\alpha(\eta) \dG \\
& \quad = - \int_\Omega \nabla (\NN(\theta) - C_\Omega) \cdot \nabla \NN(\eta) \dx - \int_\Gamma \gradg (\NNG^\alpha(\theta) - \alpha C_\Gamma) \cdot \gradg \NNG^\alpha(\eta) \dG \\
& \quad = -\int_\Omega \mu \eta \dx - \int_\Gamma \frac{\phi}{\alpha^2} \eta \dG.
\end{align*}
Choosing $\eta = \tilde \eta - c_1 - c_2 \zeta$ as in \eqref{eta:test} we obtain the identification
\begin{align}\label{LW:lim:GF:3}
\NN(u_t) - b_\Omega = - \mu \text{ in } \Omega, \quad \frac{1}{\alpha} \NNG^\alpha(u_t) - b_\Gamma = - \frac{\phi}{\alpha^2} \text{ on } \Gamma,
\end{align}
where the constants $b_\Omega$ and $b_\Gamma$ independent of $\tilde \eta$ are defined as
\begin{align*}
b_\Omega := \frac{\mean{(\mu - \NN(u_t))\zeta}_\Omega}{\mean{\zeta}_\Omega}, \quad b_\Gamma := \frac{1}{\abs{\Gamma}} \Big ( \int_\Omega \mu \dx + \int_\Gamma \frac{\phi}{\alpha^2} \dG \Big ) - b_\Omega \frac{\abs{\Omega}}{\abs{\Gamma}}.
\end{align*}
Using the definition of $\NN$ and $\NNG^{\alpha}$, the system comprising of \eqref{LW:lim:GF:2} and \eqref{LW:lim:GF:3} is exactly the Cahn--Hilliard system \eqref{lim:alt}.

\subsection{Robin system \eqref{CHLW}}
The first variation of $E$ with respect to $\bu = (u, v)$ in direction $\bm{\eta} = (\eta, \phi) \in \Ho\cap \CL^\infty$ reads as
\begin{equation}\label{firstvar}
\begin{aligned}
\frac{\delta E}{\delta \bm{u}}(\bm{u}) [\bm{\eta}] & = \int_\Omega \nabla u \cdot \nabla \eta + F'(u) \eta \dx  \\
& \quad + \int_\Gamma \kappa \gradg v \cdot \gradg \phi + G'(v) \phi + K^{-1}(u - H(v))(\eta - H'(v)\phi) \dG.
\end{aligned}
\end{equation}
Recalling the inner product $(\cdot,\cdot)_{(\Ho)'}$ defined in \eqref{defn:H0':inn}, we seek the gradient $\nabla E \vert_{\bm{u}}$ of $E$ at $\bm{u}$ associated to this inner product, which satisfies
$(\nabla E \vert_{\bm{u}}, \bm{\eta})_{(\Ho)'} = \frac{\delta E}{\delta \bm{u}}(\bm{u})[\bm{\eta}]$ for all $\bm{\eta} \in \Ho  \cap \CL^\infty$.  In particular, this yields for $(\theta, \zeta) = \nabla_{(\Ho)'} E \vert_{\bm{u}}$ the identity 
\begin{align}\label{CHLW:GF:1}
\notag \int_\Omega \NN(\theta) \eta \dx + \int_\Gamma  \NNG(\zeta) \phi \dG
 & = \int_\Omega (F'(u) - \Lx u) \eta \dx  + \int_\Gamma (K^{-1}(u - H(v)) + \pdnu u) \eta \dG \\
& \quad  + \int_\Gamma (G'(v) + K^{-1} H'(v)(H(v) - u) - \kappa \LB v) \phi \dG.
\end{align}
For arbitrary test functions $\tilde \eta \in H^1(\Omega) \cap L^\infty(\Omega)$, $\tilde \phi \in H^1(\Gamma) \cap L^\infty(\Gamma)$, we observe that $\bm{\eta} := (\tilde \eta - \mean{\tilde \eta}_\Omega, \tilde \phi - \mean{\tilde \phi}_\Gamma) \in \Ho \cap \CL^\infty$.  Defining constants 
\begin{align}
\label{DEF:CONS}
c_\Omega := \mean{F'(u)}_\Omega + \frac{\abs{\Gamma}}{\abs{\Omega}} \mean{ K^{-1}(u - H(v))}_\Gamma, \quad c_\Gamma := \mean{G'(v) + K^{-1} H'(v)(H(v) - u)}_\Gamma,
\end{align}
independent of $(\tilde \eta, \tilde \phi)$ and plugging the choice of $\bm{\eta}$ into \eqref{CHLW:GF:1} yields the following identifications after a short calculation:
\begin{equation}\label{LW:GF:2}
\begin{aligned}
\mu := \NN(\theta) + c_\Omega &= F'(u) - \Lx u &&\text{ in } \Omega, \\
\muG := \NNG(\zeta) + c_\Gamma &= G'(v) + H'(v) \pdnu u - \kappa \LB v,
\quad K \pdnu u = H(v) - u &&\text{ on } \Gamma.
\end{aligned}
\end{equation}
Thus, the gradient flow
\begin{align}\label{LW:GF}
(\bm{u}_t, \bm{\eta})_{(\Ho)'} = - \frac{\delta E}{\delta \bm{u}}(\bm{u})[\bm{\eta}] = - (\nabla E \vert_{\bm{u}},\bm{\eta})_{(\Ho)'} \quad \text{ for all } \bm{\eta} \in \Ho \cap \CL^\infty
\end{align}
can now be expressed as
\begin{align*}
& \int_\Omega \NN(u_t) \eta \dx + \int_\Gamma \NNG(v_t) \phi \dG = \int_\Omega \nabla \NN(u_t) \cdot \nabla \NN(\eta) \dx + \int_\Gamma \gradg \NNG(v_t) \cdot \gradg \NNG(\phi) \dG \\
& \quad = - \int_\Omega \nabla \mu \cdot \nabla \NN(\eta) \dx - \int_\Gamma \gradg \muG \cdot \gradg \NNG(\phi) \dG = - \int_\Omega \mu \eta \dx - \int \muG \phi \dG.
\end{align*}
Choosing $\eta := \tilde \eta - \mean{\tilde \eta}_\Omega$ and $\phi := \tilde \phi - \mean{\tilde \phi}_\Gamma$ for arbitrary $\tilde \eta \in H^1(\Omega) \cap L^\infty(\Omega)$ and $\tilde \phi \in H^1(\Gamma) \cap L^\infty(\Gamma)$, we obtain the identifications
\begin{align}\label{LW:GF:3}
\NN(u_t) = - \mu + c_\Omega \text{ in } \Omega, \quad \NNG(v_t) = - \muG + c_\Gamma \text{ on } \Gamma.
\end{align}
Using the definition of $\NN$ and $\NNG$, the system comprising of \eqref{LW:GF:2} and \eqref{LW:GF:3} is exactly the Cahn--Hilliard system \eqref{CHLW}.

\section{Proof of Theorem \ref{thm:lim}}\label{sec:lim}
We define
\begin{align*}
\tilde G(s) = G(\alpha^{-1}(s - \beta)) \quad \text{ for all } s \in \R,
\end{align*}
and it is easy to see that the assumption \eqref{ass:pot} implies that $\tilde G$ satisfies \cite[(P2)]{GK}, namely, $\tilde G$ is bounded from below and admits a decomposition $\tilde G = \tilde G_1 + \tilde G_2$ where $\tilde G_1$ is convex,  nonnegative and satisfies for any $\delta > 0$,
\begin{align*}
\abs{\tilde G_1'(s)} = \frac{1}{\abs{\alpha}} \abs{G_1'(\alpha^{-1}(s-\beta))} \leq \frac{\delta}{\abs{\alpha}} \abs{G_1(\alpha^{-1}(s-\beta))} + \frac{C_{G,\delta}}{\abs{\alpha}} = \frac{\delta}{\abs{\alpha}} \tilde G_1(s) + \tilde C_{G,\delta},
\end{align*}
for some positive constant $\tilde C_{G,\delta}$ depending only on $\delta$,
while $\tilde G_2'$ is Lipschitz continuous
\begin{align*}
\abs{\tilde G_2'(s) - \tilde G_2'(r)} = \frac{1}{\abs{\alpha}} \abs{G_2'(\alpha^{-1}(s-\beta)) - G_2'(\alpha^{-1}(r-\beta))} \leq \frac{C_g}{\alpha^2} \abs{s-r},
\end{align*}
where $C_g$ is the Lipschitz constant of $G_2$.

Since the main ideas have been outlined in \cite{GK}, we simply sketch the crucial steps of the proof.  
\paragraph{Step 1: Implicit time discretisation.}
Let $N \in \N$ be arbitrary and $\tau := T/N$ denote a time step size.  For $n = 0, 1, 2, \dots, N$, we define functions $u^n$ recursively by the following construction.  The zeroth order iterate is taken to be the initial data, i.e., $u^0 := u_0$.  If $u^n$ is already constructed, then we choose $u^{n+1}$ to be a minimiser of the functional
\begin{align*}
J_n(u)  := \frac{1}{2 \tau} \norm{u - u^n}_{(\V^\kappa_0)'_*}^2 + E_*(u)
\end{align*}
on the set $\V^{\kappa}_m$. We point out that $J_n$ may attain the value $+\infty$, since the $\V^{\kappa}$-regularity for $u$ may not guarantee that $E_*(u)$ is finite if the potentials $F$ and $G$ grow sufficiently fast.  However, any minimiser to $J_n$ will have finite energy due to its construction.  By following the proof of \cite[Lem.~5]{GK} we find that $J_n$ is bounded from below for all $u \in \V^{\kappa}_m$ and via the direct method of the calculus of variations, we infer the existence of a global minimiser to $J_n$ on the space $\V^{\kappa}_m$.  Hence, the existence of $u^{n+1}$ is guaranteed.  Employing the convexity of $F_1$ and $G_1$, we proceed as in \cite{garckeelas} to obtain the Euler--Lagrange equation
\begin{align*} 
0 & = \inn{\tau^{-1}(u^{n+1} - u^n)}{\eta}_{(\V^{\kappa}_0)'} + \int_\Omega \nabla u^{n+1} \cdot \nabla \eta + F'(u^{n+1}) \eta \dx \\
& \quad + \int_\Gamma \frac{\kappa}{\alpha^2} \gradg u^{n+1} \cdot \gradg \eta + \frac{1}{\alpha} G'(\alpha^{-1}(u^{n+1}-\beta)) \eta \dG
\end{align*}
for any $\eta \in \V^{\kappa}_0 \cap L^\infty(\Omega)$  with $\eta\vert_\Gamma \in L^\infty(\Gamma)$.  Setting
\begin{align*}
\mathring{\mu}^{n+1} := \NN \Big ( - \tau^{-1}(u^{n+1} - u^n) \Big ), \quad \mathring{\phi}^{n+1} = \NNG^\alpha \Big ( - \tau^{-1}(u^{n+1} - u^n ) \Big )
\end{align*}
so that by \eqref{pre:Na}, $\mathring{\phi}^{n+1} \in H^1(\Gamma)_0$ and $\mathring{\mu}^{n+1} \in H^1(\Omega)_0$ are, respectively, the unique solutions to the Poisson problems
\begin{align*}
- \LB \mathring{\phi}^{n+1} = -\frac{u^{n+1} - u^n}{\alpha \tau} \;\text{ on } \Gamma, \qquad 
\begin{cases}
- \Lx \mathring{\mu}^{n+1}  = -\frac{u^{n+1} - u^n}{\tau} & \text{ in } \Omega, \\ \pdnu \mathring{\mu}^{n+1} = 0 & \text{ on } \Gamma,
\end{cases}
\end{align*}
it holds that
\begin{align*}
(\tau^{-1}(u^{n+1} - u^n), \eta)_{(\V^{\kappa}_0)'_*} = - \int_\Omega \mathring{\mu}^{n+1} \eta \dx - \int_\Gamma \frac{1}{\alpha} \mathring{\phi}^{n+1} \eta \dG \quad \text{ for any } \eta \in \V^{\kappa}_0.
\end{align*}
Furthermore, for any arbitrary nonnegative $\zeta \in C^{\infty}_c(\Omega)$ not identically zero, and any arbitrary $\tilde \eta \in \Vk \cap L^\infty(\Omega)$ with $\tilde \eta\vert_\Gamma\in L^\infty(\Gamma)$, we define $\eta = \tilde \eta - c_1 - c_2 \zeta \in \V^{\kappa}_0 \cap L^\infty(\Omega)$ with $\eta\vert_\Gamma\in L^\infty(\Gamma)$ where $c_1$ and $c_2$ are the constants defined in \eqref{LW:GF:c1c2}, as well as constants
\begin{align*}
C_\Omega^{n+1} & := \frac{1}{\abs{\Omega} \mean{\zeta}_\Omega} \Big (\int_\Omega ( F'(u^{n+1}) - \mathring{\mu}^{n+1} )\zeta- \nabla u^{n+1} \cdot \nabla \zeta \dx \Big ), \\
 C_\Gamma^{n+1} & := \frac{1}{\abs{\Gamma}} \Big ( \int_\Omega F'(u^{n+1}) \dx + \int_\Gamma \frac{1}{\alpha} G'(\alpha^{-1}(u^{n+1} - \beta)) \dG \Big ) - C_\Omega^{n+1} \frac{\abs{\Omega}}{\abs{\Gamma}}
\end{align*}
that are independent of $\tilde \eta$. Then the functions
\begin{align*}
\mu^{n+1} := \mathring{\mu}^{n+1} + C_\Omega^{n+1}, \quad \phi^{n+1} := \alpha \mathring{\phi}^{n+1} + \alpha^2 C_\Gamma^{n+1}
\end{align*}
satisfy 
\begin{subequations}\label{EQ:DWFLIM}
\begin{alignat}{2}
0 & = \int_\Omega \nabla \mu^{n+1} \cdot \nabla \zeta + \frac{u^{n+1} - u^n}{\tau} \zeta \dx , \\
0 & = \int_\Gamma \gradg \phi^{n+1} \cdot \gradg \theta + \frac{u^{n+1} - u^n}{\tau} \theta \dG , \\
0 & = \int_\Omega (F'(u^{n+1}) - \mu^{n+1})\tilde \eta + \nabla u^{n+1} \cdot \nabla \tilde \eta \dx &&\notag \\
& \quad + \int_\Gamma \left( \frac{1}{\alpha} G' \Big ( \frac{u^{n+1}-\beta}{\alpha} \Big)  - \frac{\phi^{n+1}}{\alpha^2} \right) \tilde \eta + \frac{\kappa}{\alpha^2} \gradg u^{n+1} \cdot \gradg \tilde \eta \dG ,
\end{alignat}
\end{subequations}
for all $\zeta \in H^1(\Omega)$, $\theta \in H^1(\Gamma)$ and $\tilde \eta \in \V^{\kappa} \bl \cap L^\infty(\Omega)$ with $\eta\vert_\Gamma\in L^\infty(\Gamma)$.  This leads to an implicit time discretisation of \eqref{lim:alt}.

Hence, the collection $\{(u^n, \mu^n, \phi^n)\}_{n = 1,\dots, N}$ describes a time-discrete approximate solution.  For $t \in [0,T]$ and $n \in \{1,2, \dots, N\}$, we define the piecewise constant extension
\begin{align*}
(u_N, \mu_N, \phi_N)(\cdot, t) :=  (u^n, \mu^n, \phi^n) 
\end{align*}
for $t \in ((n-1)\tau, n\tau]$, and the piecewise linear extension
\begin{align*}
(\bar u_N, \bar \mu_N, \bar \phi_N)(\cdot, t)   := \gamma(u^n, \mu^n, \phi^n) + (1-\gamma)(u^{n-1}, \mu^{n-1}, \phi^{n-1}) 
\end{align*}
for any $\gamma \in [0,1]$ and $t = \gamma n \tau + (1-\gamma)(n-1)\tau$.

\paragraph{Step 2: Uniform estimates.} Following \cite[Sec.~4.4 and 4.5]{GK} we immediately infer the following uniform estimates:  There exists a positive constant $C$ independent of $N$, $n$ and $\tau$, such that for any $t_1, t_2 \in [0,T]$,
\begin{align*}
& \norm{u_N}_{L^\infty(0,T;H^1(\Omega))} + \norm{\mu_N}_{L^2(0,T;H^1(\Omega))} + \norm{\phi_N}_{L^2(0,T;H^1(\Gamma))} \leq C, \\
& \norm{u_N}_{L^\infty(0,T;H^p(\Gamma))} \leq C, \quad \text{ where }\; p = \begin{cases}
\frac 1 2, & \text{ if } \kappa = 0, \\
1, & \text{ if } \kappa > 0,
\end{cases} \\
&\norm{\bar u_N(t_1) - \bar u_N(t_2)}_{L^2(\Omega)} \leq C \abs{t_1 - t_2}^{1/4}, \\
&\norm{\bar u_N(t_1) - \bar u_N(t_2)}_{\Yk} \leq C \abs{t_1 - t_2}^{1/4},
\end{align*}
where $\Yk$ is defined in \eqref{DEF:XY}.  Then, from the uniform boundedness of $\nabla \mu_N$ and $\gradg \phi_N$, we also deduce the boundedness of $\bar u_{N, t}$ in $L^2(0,T;H^1(\Omega)')$ and $\bar u_{N,t}$ in $L^2(0,T;H^1(\Gamma)')$.  

\paragraph{Step 3: Passing to the limit.}  Proceeding as in \cite[Sec.~4.6]{GK}, we obtain the existence of limit functions $(u, \mu, \phi)$ in the limit $N \to \infty$ possessing the regularities stated in item (i) of Theorem \ref{thm:lim}.  Items (ii) and (iii) of Theorem \ref{thm:lim} can also be shown analogously to \cite[Sec.~4.7]{GK}, and the uniqueness of the weak solution follows as in \cite[Sec.~5]{GK}.

\paragraph{Step 4: Higher regularity.} We now establish the $H^2(\Omega)$-regularity assertion by following the arguments outlined in \cite[Sec.~4]{CC}.  The new hypothesis \eqref{GK:f1} implies that $F'(u) \in L^{\infty}(0,T;L^2(\Omega))$, and so the third equality in item (ii) of Theorem \ref{thm:lim} can be made to hold for all $\psi \in \V^{\kappa}$.  Now, we choose $\psi \in H^1_0(\Omega)$ and obtain the elliptic equation
\begin{align}\label{ell:b}
- \Lx u = \mu - F'(u)
\end{align}
holding in the sense of distributions in $\Omega$, with a right-hand side belonging to $L^2(Q_T)$.  Elliptic regularity theory (see for instance \cite[Thm.~A.2]{CFL} with $g_0 = u \vert_{\Sigma_T}$, $p = 2$, $A = - \Lx$, $f = \mu - F'(u)$, $r = 0$, $t = 1$) yields that for $\kappa > 0$,
\begin{align*}
\norm{u}_{L^2(0,T;H^{\frac{3}{2}}(\Omega))} \leq C \Big ( \norm{\Lx u}_{L^2(0,T;L^2(\Omega))} + \norm{u \vert_{\Sigma_T}}_{L^2(0,T;H^1(\Gamma))} \Big ) \leq C.
\end{align*}
Since $\Lx u \in L^2(0,T;L^2(\Omega))$ and thus $u \in L^2(0,T;H^{3/2}(\Omega))$, by a variant of the trace theorem involving elliptic operators (see for instance \cite[Thm.~A.1]{CFL} with $k = 1$, $p = 2$, $s = \frac{3}{2}$), we have $\pdnu u \in L^2(\Sigma_T)$.  Then, using \eqref{ell:b} and the third equality in item (ii) of Theorem \ref{thm:lim} we infer the surface elliptic equation
\begin{align}\label{ell:s}
- \frac{\kappa}{\alpha^2} \LB u = \frac{\phi}{\alpha^2} - \pdnu u - \frac{1}{\alpha} G' \Big ( \frac{u - \beta}{\alpha} \Big )
\end{align}
holding in the sense of distributions on $\Gamma$, with the right-hand side belonging to $L^2(\Sigma_T)$.  
Hence, by elliptic regularity theory, we have $u \vert_{\Sigma_T} \in L^2(0,T;H^2(\Gamma))$, and when revisiting the elliptic equation \eqref{ell:b} now with a more regular trace $u \vert_{\Sigma_T}$, we have for $\kappa > 0$,
\begin{align*}
\norm{u}_{L^2(0,T;H^2(\Omega))} \leq C \Big ( \norm{\Lx u}_{L^2(0,T;L^2(\Omega))} + \norm{u \vert_{\Sigma_T}}_{L^2(0,T;H^2(\Gamma))} \Big ) \leq C.
\end{align*}
For the case $\kappa = 0$ we take note of the Sobolev embedding $H^{1/2}(\Gamma) \emb L^4(\Gamma)$, so that it follows from \eqref{ell:s} that $\pdnu u \in L^2(\Sigma_T)$.  Then, by elliptic regularity theory (see \cite[Thm.~A.2]{CFL} with $A = - \Lx$, $f = \mu - F'(u)$, $r = 0$, $p = 2$, $g_1 = \pdnu u$, $t = 0$), we obtain $u \in L^2(0,T;H^{\frac{3}{2}}(\Omega))$.  By the trace theorem, $u \vert_{\Sigma_T} \in L^2(0,T;H^1(\Gamma))$ and hence, due to \eqref{GK:f1}, it holds that $G'((u-\beta)/\alpha) \in L^2(0,T;H^1(\Gamma))$.  From \eqref{ell:s} we now have $\pdnu u \in L^2(0,T;H^1(\Gamma))$ and by \cite[Thm.~A.2]{CFL} (with $t = 1$) it holds that $u \in L^2(0,T;H^2(\Omega))$.  This finishes the proof of Theorem \ref{thm:lim}.

\section{Proof of Theorem \ref{thm:CHLW}}\label{sec:CHLW}
The proof is divided into seven steps. Throughout this section, the symbol $C$ will denote generic positive constants independent of $N$, $n$ and $\tau$ that may change their value from line to line.

\paragraph{Step 1: Implicit time discretisation.} 	Let $N \in \N$ be arbitrary and $\tau := T/N$ denote a time step size.  For $n = 0, 1, 2, \dots, N$, we define functions $\bu^n = (u^n, v^n)$ recursively by the following construction.  The zeroth order iterate is taken to be the initial data, i.e., $\bu^0 = (u^0,v^0) := (u_0,v_0)$.  If $\bu^n$ is already constructed, we choose $\bu^{n+1}$ to be a minimiser of the functional
\begin{align*}
J_n(\bu) := \frac{1}{2 \tau} \norm{\bu - \bu^n}_\Hod^2 + E(\bu)
\end{align*}
on the set 
\begin{align}
\Zmk := \left\{ (u,v)\in H^1(\Omega) \times \Xk \;\middle|\; \mean{u}_\Omega = \mo, \,
\mean{v}_\Gamma = \mg \right\}.
\end{align}
As in section \ref{sec:lim}, if the potentials $F$ and $G$ grow sufficiently fast, the functional $J_n$ may attain the value $+\infty$, but minimisers, whose existence are addressed in Step 2, will have finite energy. At this moment we briefly describe the construction of time-discrete solutions. Convexity of $F_1$ and $G_1$ allow us to proceed as in \cite{garckeelas} to infer that $\bu^{n+1}$ satisfies 
\begin{equation}\label{EQ:DGF}
\begin{aligned}
0 &=  (\tau^{-1} (\bu^{n+1} - \bu^n) , \feta)_\Hod + \intO \grad u^{n+1} \cdot \grad \eta + F'(u^{n+1}) \eta \dx \\
&\quad + \intG \kappa \gradg v^{n+1} \cdot \gradg \psi + G'(v^{n+1}) \psi + K^{-1}(u^{n+1}-H(v^{n+1}))(\eta - H'(v^{n+1})\psi)\dG \\
& = \intO  \nabla u^{n+1} \cdot \nabla \eta + F'(u^{n+1}) \eta - \mathring \mu^{n+1} \eta \dx + \intG \kappa \gradg v^{n+1} \cdot \gradg \psi \dG   \\
& \quad + \intG G'(v^{n+1}) \psi - \mathring \mu_\Gamma^{n+1}  \psi  + K^{-1}(u^{n+1} - H(v^{n+1}))(\eta - H'(v^{n+1}) \psi) \dG,
\end{aligned}
\end{equation}
for all $\feta=(\eta,\psi) \in\Ho \cap \CL^\infty$, where
\begin{align*}
\mathring \mu^{n+1} := \NN\big(-\tau^{-1}(u^{n+1} - u^n)\big) \in H^1(\Omega)_0,
\quad \mathring \mu_\Gamma^{n+1} :=\NNG\big(-\tau^{-1}(v^{n+1} - v^n)\big) \in H^1(\Gamma)_0.
\end{align*}
For arbitrary $\tilde \eta \in H^1(\Omega) \cap L^{\infty}(\Omega)$ and $\tilde \psi \in H^1(\Gamma) \cap L^{\infty}(\Gamma)$, we define 
\begin{align*}
\bm{\eta} = (\eta, \psi) := (\tilde \eta - \mean{\tilde \eta}_\Omega, \tilde \psi - \mean{\tilde \psi}_{\Gamma}) \in \Ho  \cap \CL^\infty
\end{align*}
as well as functions
\begin{align*}
\mu^{n+1} := \mathring \mu^{n+1} + c_\Omega^{n+1}, \qquad \mu_\Gamma^{n+1} := \mathring \mu_\Gamma^{n+1} + c_\Gamma^{n+1}
\end{align*}
where the constants are given by
\begin{align*}
c_\Omega^{n+1} &:= \mean{F'(u^{n+1})}_\Omega + \frac{\abs{\Gamma}}{\abs{\Omega}} \mean{K^{-1}(u^{n+1} - H(v^{n+1}))}_\Gamma, \\
c_\Gamma^{n+1} &:= \mean{G'(v^{n+1}) + K^{-1} H'(v^{n+1})(H(v^{n+1}) - u)}_\Gamma.
\end{align*}
Then, following the previous section, we infer from \eqref{EQ:DGF} that the quadruplet $(u^{n+1}$, $v^{n+1}$, $\mu^{n+1}$, $\mu_\Gamma^{n+1})$ satisfies the equations
\begin{align}
\label{EQ:DWF1}
\intO \frac{u^{n+1}-u^n}{\tau} \,\zeta \dx &= -\intO \grad\mu^{n+1} \cdot \grad\zeta \dx,\\
\label{EQ:DWF2}
\intG \frac{v^{n+1}-v^n}{\tau} \, \theta \dx &= -\intG \gradg \mu^{n+1} \cdot \gradg \theta \dG,\\
\label{EQ:DWF3}
\intO \grad u^{n+1}\cdot\grad \eta + F'(u^{n+1})\eta -  \mu^{n+1}\eta \dx &= \intG \tfrac{1}{K} (H(v^{n+1})-u^{n+1})  \eta  \dG,\\
\label{EQ:DWF4}
\intG  \kappa \gradg v^{n+1}\cdot \gradg \psi + G'(v^{n+1})\psi - \mu_\Gamma^{n+1} \psi \dG &= - \intG \tfrac{1}{K} \big( H(v^{n+1}) - u^{n+1} \big) H'(v^{n+1}) \psi  \dG.
\end{align}
for all $\zeta\in H^1(\Omega)$, $\theta\in H^1(\Gamma)$, $\eta\in H^1(\Omega) \cap L^\infty(\Omega)$, $\psi\in \Xk \cap L^\infty(\Gamma)$.  This system is an implicit time discretisation of \eqref{CHLW}. Hence, the collection $\{u^{n}, v^{n}, \mu^{n}, \mu_\Gamma^{n}\}_{n =1, \dots, N}$ is a time-discrete approximate solution.  For $t \in [0,T]$ and $n \in \{1,2, \dots, N\}$, we define the piecewise constant extension
\begin{align}\label{pi:cons}
(u_N, v_N, \mu_N, \muGN)(\cdot, t) :=  (u_N^n, v_N^n, \mu_N^n, \muGN^n) :=  (u^n, v^n, \mu^n, \muG^n),
\end{align}
for $t \in ((n-1)\tau, n\tau]$ and the piecewise linear extension
\begin{align}\label{pi:lin}
(\lu_N, \lv_N, \lmu_N, \lmuGN)(\cdot, t)  := \gamma(u_N^n, v_N^n, \mu_N^n, \muGN^n) + (1-\gamma)(u_N^{n-1}, v_N^{n-1}, \mu_N^{n-1}, \muGN^{n-1})
\end{align}
for any $\gamma \in [0,1]$ and $t = \gamma n \tau + (1-\gamma)(n-1)\tau$.

\paragraph{Step 2: Existence of a minimiser.}
We can use the direct method in the calculus of variations to show that the functional $J_n$ has at least one minimiser in the set $\Zmk$. First, we observe that $J_n$ is bounded from below:
\begin{align}\label{BND:JN}
J_n(\bu) \ge \intO F(u) \dx + \intG G(v) \dG \ge C^* \quad\text{ with } C^* := C_F\abs{\Omega} + C_G\abs{\Gamma},
\end{align}
for all $\bu=(u,v)\in \Zmk$. Hence, the infimum $M:=\inf_{\Zmk} J_n$ exists and we can find a minimising sequence $(\bu_k)_{k\in\N}\subset \Zmk$ such that
\begin{align*}
\lim_{k \to \infty} J_n(\bu_k) = M, \quad 	J_n(\bu_k) \le M+1 \quad \text{ for all } k\in\N,
\end{align*}
which leads to the estimates
\begin{align}\label{BND:MC}
\frac 1 2 \norm{\grad u_k}_{L^2(\Omega)}^2 + \frac \kappa 2 \norm{\gradg v_k}_{L^2(\Gamma)}^2 + \frac{1}{2K} \norm{u_k-H(v_k)}_{L^2(\Gamma)}^2 + \intG G_1(v_k) \dG \leq C.
\end{align}
for all $k\in\N$. As $\mean{u_k}_\Omega=\mo$, we can use \eqref{BND:MC} and Poincar\'e's inequality to conclude that $(u_k)_{k \in \N}$ is a bounded sequence in $H^1(\Omega)$. Thus, there exists a function $\bar u\in H^1(\Omega)$ such that $u_k \wto \bar u$ in $H^1(\Omega)$ as $k\to\infty$ along a nonrelabelled subsequence. This implies that $\mo = \mean{u_k}_\Omega\to\mean{\bar u}_\Omega$ and $\mean{\bar u}_\Omega = \mo$ directly follows. 

For the surface variable we have two cases: If $\kappa>0$ we proceed similarly as above and infer the existence of a function $\bar v\in H^1(\Gamma)$ with $v_k\wto \bar v$ in $H^1(\Gamma)$ and $\mean{\bar v}_\Gamma = \mg$. This implies that the limit $\bar\bu=(\bar u,\bar v)$ belongs to $\Zmk$.  In the case $\kappa=0$ we recall that $H(s)=\alpha s+\beta$ for all $s\in\R$. By \eqref{BND:MC} and the continuous embedding $H^1(\Omega)\emb L^2(\Gamma)$ we get
\begin{align*} 
\norm{v_k}_{L^2(\Gamma)} \le C + C\norm{H(v_k)}_{L^2(\Gamma)} \le C + C\norm{u_k-H(v_k)}_{L^2(\Gamma)} +  C\norm{u_k}_{H^1(\Omega)} \le C.
\end{align*} 
Furthermore, from \eqref{ass:pot:5:k}, we can use \eqref{BND:MC} to obtain the estimate
\begin{align*}
\norm{v_k}_{L^p(\Gamma)}^p \le C + C \intG G_1(v_k) \dG \le C.
\end{align*} 
Hence, there exists a function $\bar v\in \Xk$ such that $v_k\wto \bar v$ in $\Xk$. The weak convergence implies that $\mg = \mean{v_k}_\Gamma \to \mean{\bar v}_\Gamma$ and $\mean{\bar v}_\Gamma = \mg$ directly follows. This means that $\bar\bu=(\bar u,\bar v)\in \Zmk$.
	
It remains to prove that $\bar\bu$ is a minimiser of the functional $J_n$.  The convergence $u_k\wto \bar u$ in $H^1(\Omega)$, continuity and nonnegativity of $F(\cdot) - C_F$, and Fatou's lemma yield
\begin{align*}
\intO F(\bar u) \dx &\le \underset{k \to \infty}{\lim\inf}\intO F(u_k) \dx.
\end{align*}
If $\kappa>0$ we can proceed similarly and show that
\begin{equation}\label{LIMSUP}
\begin{aligned}
\intG G(\bar v) \dG & \le \underset{k\to \infty}{\lim\inf}\intO G(v_k) \dG, \\
 \intG \abs{\bar u-H(\bar v)}^2 \dG &\le \underset{k\to \infty}{\lim\inf} \intG \abs{u_k-H(v_k)}^2 \dG.
\end{aligned}
\end{equation}
In the case $\kappa=0$ we recall that $G=G_1$ is continuous and convex (hence weakly lower semicontinuous) and that
	\begin{align*}
	\intG \abs{\bar u-H(\bar v)}^2 \dG = \intG (\bar u-\beta)^2 - 2\alpha(\bar u-\beta)\bar v + \alpha^2\bar v^2 \dG.
	\end{align*}
Hence, \eqref{LIMSUP} holds true thanks to weak lower semicontinuity, and ensures that in all cases we have
\begin{align*}
J_n(\bar \bu) \leq \underset{k \to \infty}{\lim\inf} J_n(\bu_k) = M.
\end{align*}

\paragraph{Step 3: Uniform estimates.}
In the following, we exploit the gradient flow structure to establish uniform bounds on the piecewise constant extensions, i.e.,
\begin{equation}
	\label{EST:PCE}
\begin{aligned}
& \norm{u_N}_{L^\infty(0,T;H^1(\Omega))} + \norm{v_N}_{L^\infty(0,T;\Xk)} \\
&\qquad + \norm{\mu_N}_{L^2(0,T;H^1(\Omega))}  + \norm{\muGN}_{L^2(0,T;H^1(\Gamma))} \leq C,
\end{aligned}
\end{equation}
by following the approach in \cite[Sec.~4.4]{GK}.  Firstly, since $\buN = (u_N,v_N)$ is a minimiser of $J_n$ on the set $\Zmk$, we infer the a priori estimate
\begin{align}
\label{EST:APR}
\frac{1}{2\tau} \norm{\bu^{n+1}-\bu^n}_\Hod^2 + E(\bu^{n+1}) \le E(\bu^n) \quad \text{ for } n\in\{0,1, \dots ,N-1\},
\end{align}
and subsequently, for all $n \in \{0, 1, \dots, N-1\}$,
\begin{align}\label{EST:APR2}
\frac{1}{2} \norm{\grad u^{n+1}}_{L^2(\Omega)} + \frac \kappa 2 \norm{\gradg v^{n+1}}_{L^2(\Gamma)} + \frac{1}{2K} \norm{u-H(v)}_{L^2(\Gamma)}^2 + \intG G_1(v^{n+1}) \dG \leq C.
\end{align}
Since $\mean{u^{n+1}}_\Omega=\mo$ and $\mean{v^{n+1}}_\Gamma=\mg$, we can proceed similarly as in Step 2 and use the definition \eqref{pi:cons} to obtain the uniform bounds
\begin{align}
\label{EST:UVN}
\norm{u_N}_{H^1(\Omega)} + \norm{v_N}_{\Xk} \le C.
\end{align}
Next, testing \eqref{EQ:DWF3} with arbitrary $\eta\in C_c^\infty(\Omega)$ such that $0\le\eta\le 1$ and using \eqref{EST:UVN} gives
\begin{align*}
\abs{ \intO \mu^{n+1} \eta \dx } \le C \norm{\grad \eta}_{L^2(\Omega)} + \intO |F'(u^{n+1})| \dx 
\end{align*}
From assumption \eqref{ass:pot:3} with $\delta=1$ and Remark \ref{REM:POT} we infer that
\begin{align}
\label{EST:fL1}
&\intO |F'(u^{n+1})| \dx \le E(\bu_0) + \abs{C^*} + A_{f,1}\abs{\Omega} + c_f\abs{\Omega} 
+ c_f \abs{\Omega}^{\frac 1 2} \norm{u^{n+1}}_{L^2(\Omega)} \leq C.
\end{align}
Consequently, there exists some nonnegative constant $C(\eta)$ depending on $\eta$ but not on $N$, $\tau$ or $n$, such that
\begin{align}
\label{EST:CETA}
\abs{ \intO \mu^{n+1} \eta \dx } \le C(\eta). 
\end{align}
Proceeding as in \cite[Sec.~4.3]{GK} we can use a generalised Poincar\'e inequality (see \cite[p.~242]{Alt}) to conclude that 
\begin{align}
\label{EST:MUN}
\norm{\mu^{n+1}}_{L^2(\Omega)} \le C\, \left(1+\norm{\grad \mu^{n+1}}_{L^2(\Omega)}\right) \quad \text{ for all }  n \in \{0, 1,...,N-1\}.
\end{align}
To establish a uniform bound on $\norm{\grad \mu^{n+1}}_{L^2(\Omega)}$ we recall from \eqref{pi:cons} that
\begin{align*}
\big(u_N,v_N,\mu_N,\muGN)(s) = \big(u_N,v_N,\mu_N,\muGN)(t) = \big(u_N^n,v_N^n,\mu_N^n,\muGN^n)
\end{align*}
for all $s\in (t-\tau,t]$, $n\in\{1,...,N-1\}$ where $t=\tau n$ is fixed.  From the definition of $\mu_N$ and $\muGN$, as well as the a priori estimate \eqref{EST:APR}, we derive that
\begin{align*}
&E\big(\bu_N(t)\big) + \frac 1 2 \int_{t-\tau}^t \norm{\grad\mu_N(s)}_{L^2(\Omega)}^2  + \norm{\gradg\muGN(s)}_{L^2(\Gamma)}^2 \ds \\
&\quad = E\big(\bu_N(t)\big) + \frac{1}{2\tau^2} \int_{t-\tau}^t \norm{\bu_N(s)-\bu_N(s-\tau)}_\Hod^2 \ds\\
&\quad = E\big(\bu_N(t)\big) + \frac{1}{2\tau} \norm{\bu_N(t)-\bu_N(t-\tau)}_\Hod^2 
\le E\big(\bu_N(t-\tau)\big).
\end{align*}
It follows inductively that
\begin{align}
\label{IEQ:EUN}
E\big(\bu_N(t)\big) + \frac 1 2 \int_{0}^t \norm{\grad\mu_N(s)}_{L^2(\Omega)}^2 
+ \norm{\gradg\muGN(s)}_{L^2(\Gamma)}^2 \ds \le E(\bu_0),
\end{align}
and in particular, with $t=N\tau \equiv T$, we have
\begin{align}
\label{EST:GMUN}
\norm{\grad\mu_N}_{L^2(0,T;L^2(\Omega))}^2 +  \norm{\gradg\muGN}_{L^2(0,T;L^2(\Gamma))}^2 \le 2E(\bu_0) + 2\abs{C^*} \le C.
\end{align}
Together with \eqref{EST:MUN} we conclude that $\mu_N$ is uniformly bounded in $L^2(0,T;H^1(\Omega))$.

To establish a uniform bound on $\muGN$ we test \eqref{EQ:DWF4} with $\psi\equiv 1$. This gives
\begin{align*}
\abs{\intG \muG^{n+1} \dG} \leq  \intG \abs{G'(v^{n+1})} \dG + \tfrac{1}{K} \norm{H(v^{n+1}) - u^{n+1}}_{L^2(\Gamma)} \norm{H'(v^{n+1})}_{L^2(\Gamma)} \leq C
\end{align*}
upon using \eqref{EST:APR2}, an analogous estimate to \eqref{EST:fL1} for $G'(v^{n+1})$, the fact that $H'(v^{n+1}) = \alpha$ for $\kappa = 0$ and 
\begin{align*}
\norm{H'(v^{n+1})}_{L^2(\Gamma)} \leq C \big (1 + \norm{v^{n+1}}_{H^1(\Gamma)}^{2p} \big ) \leq C \quad \text{ for } \kappa > 0.
\end{align*}
Thus, $\mean{\muGN}_\Gamma$ is uniformly bounded in $L^\infty(0,T)$, and by Poincar\'e's inequality we conclude that $\muGN$ is uniformly bounded in $L^2(0,T;H^1(\Gamma))$. 

\paragraph{Step 4: H\"older-in-time estimates.}
Via interpolation type arguments we can show that the piecewise linear extension is H\"older continuous in time, i.e., 
\begin{align}
\label{EST:HLD}
\norm{\lu(t_1)-\lu(t_2)}_{L^2(\Omega)} + \norm{\lv(t_1)-\lv(t_2)}_{\Yk} \le C \abs{t_1-t_2}^{\frac 1 4}.
\end{align}
To prove this assertion let $t_1,t_2\in[0,T]$ be arbitrary where $t_1 > t_2$ without loss of generality. Since $\lu_N$ and $\lv_N$ are continuous and piecewise linear in time, they are weakly differentiable and we infer from \eqref{EQ:DWF1} and \eqref{EQ:DWF2} that 
\begin{align}
\label{EQ:IDT}
\intO \pd_t \lu_N \zeta \dx = -\intO \grad \mu_N \cdot \grad\zeta \dx, \quad 
\intG \pd_t \lv_N \theta \dx = -\intG \gradg \muGN \cdot \gradg \theta \dG
\end{align}
for all $t\in[0,T]$ and all test functions $\zeta\in H^1(\Omega)$ and $\theta\in H^1(\Gamma)$. Choosing $\zeta := \lu_N(t_2)-\lu_N(t_1)$, integrating with respect to $t$ from $t_2$ to $t_1$ and using the uniform bounds from Step 3, we obtain that
\begin{align*}
&\norm{\lu_N(t_1)-\lu_N(t_2)}_{L^2(\Omega)}^2 = -\int_{t_2}^{t_1} \intO \grad\mu_N\cdot\big(\lu_N(t_1)-\lu_N(t_2)\big) \dx \dt\\
&\quad\le 2\norm{u_N}_{L^\infty(0,T;H^1(\Omega))} \norm{\mu_N}_{L^2(0,T;H^1(\Omega))} \abs{t_1-t_2}^{\frac 1 2} \leq C \abs{t_1-t_2}^{\frac 1 2}.
\end{align*}
In the case $\kappa>0$ an analogous estimate can be established for $\lv$. If $\kappa=0$, we recall that $\mean{\lv(t_1)-\lv(t_2)}_\Gamma=0$ and choose $\theta:= \NNG\big(\lv(t_1)-\lv(t_2)\big)$. Proceeding similar as above and using integration by parts yields
\begin{align*}
&\norm{\lv_N(t_1)-\lv_N(t_2)}_{H^1(\Gamma)'}^2 = \intG \big ( \lv_N(t_1) - \lv_N(t_2) \big ) \theta \dG \\
& \quad = - \int_{t_2}^{t_1} \intG \gradg \muGN (t) \cdot \gradg \theta \dG \dt = \int_{t_2}^{t_1} \intG \muGN(t) \big ( \lv_N(t_1) - \lv_N(t_2) \big ) \dG \dt  \\
& \quad \leq 2 \norm{v_N}_{L^{\infty}(0,T;L^2(\Gamma))} \norm{\muGN}_{L^2(\Sigma_T)} \abs{t_1 - t_2}^{\frac{1}{2}} \leq C \abs{t_1 - t_2}^{\frac{1}{2}}.
\end{align*}
For any $t \in [0,T]$, we can find $n \in \{1, \dots, N\}$ and $\gamma \in [0,1]$ such that $t = \gamma n \tau + (1-\gamma) (n-1) \tau$.  Then, an immediate consequence of \eqref{EST:HLD} and the definitions \eqref{pi:cons} and \eqref{pi:lin} are the uniform estimates
\begin{equation}\label{DIFF}
\begin{aligned}
\norm{\lu_N(t) - u_N(t)}_{L^2(\Omega)}  &\leq (1-\gamma) \norm{u_N^n - u_N^{n-1}}_{L^2(\Omega)} \\
& = (1-\gamma) \norm{\lu_N(n \tau) - \lu_N((n-1)\tau)}_{L^2(\Omega)} \leq C \tau^{1/4}, \\
\norm{\lv_N(t) - v_N(t)}_{\Yk} & \leq C \tau^{1/4}.
\end{aligned}
\end{equation}
Moreover, by integrating \eqref{EQ:IDT} with respect to $t$ from $0$ to $T$ for test functions $\zeta \in L^2(0,T;H^1(\Omega))$ and $\theta \in L^2(0,T;H^1(\Gamma))$, it is easy to see that
\begin{align}\label{EST:pdt}
\norm{\pd_t \lu_N}_{L^2(0,T;H^1(\Omega)')} + \norm{\pd_t \lv_N}_{L^2(0,T;H^1(\Gamma)')} \leq C.
\end{align}

\paragraph{Step 5: Convergence assertions.} 
Thanks to the uniform estimates \eqref{EST:PCE}, \eqref{EST:HLD}, \eqref{DIFF} and \eqref{EST:pdt}, we deduce the existence of functions
\begin{align}
\label{REG:LIM}
\left\{\;\;
\begin{aligned}
&u\in C^{0,\frac 1 4}([0,T];L^2(\Omega)) \cap L^\infty(0,T;H^1(\Omega)) \cap H^1(0,T;H^1(\Omega)'), \\
&v \in C^{0,\frac 1 4}([0,T];\Yk) \cap L^\infty(0,T;\Xk) \cap H^1(0,T;H^1(\Gamma)'), \\ 
&\mu\in L^2(0,T;H^1(\Omega)), \qquad \mu_\Gamma\in L^2(0,T;H^1(\Gamma)),
\end{aligned}
\right.
\end{align}
such that for any $\delta \in (0,\frac{1}{4})$, the following convergences hold:
\begin{align*}
& u_N \to u \text{ weakly-* in } L^{\infty}(0,T;H^1(\Omega)), \text{ strongly in } L^{\infty}(0,T;L^2(\Omega)) \text{ and a.e.~in } Q_T, \\
& \lu_N \to u \text{ weakly in } H^1(0,T;H^1(\Omega)'), \text{ strongly in } C^{0,\delta}([0,T];L^2(\Omega)) \text{ and a.e.~in } Q_T,  \\
& v_N \to v \text{ weakly-* in } L^{\infty}(0,T;\Xk), \text{ strongly in } L^{\infty}(0,T;\Yk), \\
& \lv_N \to v \text{ weakly in } H^1(0,T;H^1(\Gamma)'), \text{ strongly in } C^{0,\delta}([0,T];\Yk), \\
& \mu_N \to \mu \text{ weakly in } L^2(0,T;H^1(\Omega)), \\
& \muGN \to \muG \text{ weakly in } L^2(0,T;H^1(\Gamma)), 
\end{align*}
and in the case $\kappa > 0$ we additionally have $v_N \to v$ and $\lv_N \to v$ a.e.~on $\Sigma_T$.  Since the proof follows almost analogously from the arguments in \cite[Lem.~8]{GK}, we omit the details.

Additionally, from the trace inequality in \eqref{pre:trace} and \eqref{EST:PCE}, we find that
\begin{align*}
\norm{u_N - u}_{L^2(\Gamma)} \leq C \big ( \norm{u_N - u}_{L^2(\Omega)} + \norm{u_N - u}_{L^2(\Omega)}^{1/2} \norm{\nabla (u_N - u)}_{L^2(\Omega)}^{1/2} \big ) \leq C \norm{u_N - u}_{L^2(\Omega)}^{1/2},
\end{align*}
and so it holds that
\begin{align*}
u_N \to u \text{ strongly in } L^{\infty}(0,T;L^2(\Gamma)) \text{ and a.e.~on } \Sigma_T.
\end{align*}

\paragraph{Step 6: Existence of weak solutions.}
Due to Step 5, the limit $(u, v, \mu, \muG)$ possesses the regularity stated in item (i) of Theorem \ref{thm:CHLW}.  To verify item (ii), we pass to the limit in the weak formulation.  It is easy to see that \eqref{DEF:CHLW1} and \eqref{DEF:CHLW2} can be deduced from integrating \eqref{EQ:IDT} over $[0,T]$ for test functions $\zeta \in L^2(0,T;H^1(\Omega))$ and $\theta \in L^2(0,T;H^1(\Gamma))$, and passing to the limit $N \to \infty$.  The recovery of terms depending linearly on $(u, v, \mu, \muG)$ in \eqref{DEF:CHLW3} and \eqref{DEF:CHLW4} is straightforward, while the recovery of nonlinear terms involving $F'(u)$ and $G'(v)$ (if $\kappa > 0$) follows analogously as in \cite[Sec.~4.7]{GK} via Vitali's convergence theorem.  In the following we sketch the details for the nonlinear terms involving $H$ and also for $G'(v)$ in the case $\kappa = 0$.

For the case $\kappa > 0$, for any $\eps > 0$, let $\Sigma_\eps$ denote a measurable subset of $\Sigma_T$ with $\abs{\Sigma_\eps} < \eps$.  In virtue of \eqref{ass:h}, \eqref{EST:PCE}, H\"older's inequality and the Sobolev embedding $H^1(\Gamma) \emb L^r(\Gamma)$ for all $r < \infty$, it is easy to see that 
\begin{align*}
\int_{\Sigma_\eps} \abs{H'(v_N(t))}^2 \dG \dt &\leq C( \eps + \norm{v_N}_{L^{4q}(\Sigma_\eps)}^{2q} \eps^{1/2}) \leq C\eps^{1/2}, \\
\int_{\Sigma_\eps} \abs{H(v_N(t))}^2 \dG \dt &\leq C( \eps + \norm{v_N}_{L^{4q+4}(\Sigma_\eps)}^{2q+2} \eps^{1/2}) \leq C \eps^{1/2}.
\end{align*}
Then Vitali's convergence theorem ensures that $H(v_N) \to H(v)$ and $H'(v_N) \to H'(v)$ strongly in $L^2(\Sigma_T)$, which is sufficient to recover the right-hand sides of \eqref{DEF:CHLW3} and \eqref{DEF:CHLW4}. 

For the case $\kappa = 0$, assumption \eqref{ass:h:k=0} implies that $H'(v_N) = \alpha$ is a constant, and hence the right-hand sides of \eqref{DEF:CHLW3} and \eqref{DEF:CHLW4} can be obtained easily.  Meanwhile, to recover the term involving $G'(v)$ on the left-hand side of \eqref{DEF:CHLW4}, we employ Minty's trick.  By assumption \eqref{ass:pot:5:k}, the surface potential $G$ is convex and nonnegative.  We associate its derivative $G'$ with the operator
\begin{align*}
\CG:L^r(\Sigma_T) \to \big(L^r(\Sigma_T)\big)'\cong L^{r^*}(\Sigma_T),\quad w \mapsto \CG(w):=G'(w),
\end{align*}
where $r=\max\{2,p-1\}$ and $r^*=\frac{r}{r-1}$. This mapping is well-defined since, 
\begin{align*}
\norm{\CG(w)}_{L^{r^*}(\Sigma_T)}^{r^*} & = \intS \abs{G'(w)}^{\frac{r}{r-1}} \dG dt \le C + C \intS \abs{w}^r \dG dt < \infty
\end{align*}
for all $w\in L^r(\Sigma_T)$. In particular, since $(v_N)$ is bounded in $L^\infty(0,T;L^r(\Gamma)) \subset L^r(\Sigma_T)$, this shows that $\CG(v_N)$ is uniformly bounded in $L^{r^*}(\Sigma_T)$, and by the Banach--Alaoglu theorem  there exists $\CG^*\in L^{r^*}(\Sigma_T)$ such that $\CG(v_N) \wto \CG^*$ in $L^{r^*}(\Sigma_T)$ as $N \to \infty$ along a nonrelabelled subsequence. Then, for any arbitrary test function $\psi \in L^r(\Sigma_T)$, integrating \eqref{EQ:DWF4} over $[0,T]$ and passing to the limit leads to
\begin{align}\label{EQ:G*}
\int_{\Sigma_T} \muG \psi - \tfrac{\alpha}{K} (\alpha v + \beta - u) \psi - \CG^* \psi \dG \dt = 0.
\end{align}
The idea is to use Minty's trick to show that $\CG^*=\CG(v)$.  Convexity of $G$ implies that the operator $\CG$ is monotone.  Furthermore, for arbitrary $w, y, z \in L^r(\Sigma_T)$,
\begin{align*}
\inn{\CG(w) - \CG(y)}{z}_{L^r(\Sigma_T)} & = \int_{\Sigma_T} \int_0^1 G''(\lambda w + (1-\lambda)y)(w - y) z\, d \lambda  \dG \dt \\
& \leq C \big ( 1 + \norm{w}_{L^r(\Sigma_T)}^{r-2} + \norm{y}_{L^r(\Sigma_T)}^{r-2}\big ) \norm{w-y}_{L^r(\Sigma_T)} \norm{z}_{L^r(\Sigma_T)},
\end{align*}
which implies that $s \mapsto \inn{\CG(y+ sz)}{z}_{L^r(\Sigma_T)}$ is continuous for arbitrary $y, z \in L^r(\Sigma_T)$, i.e., $\CG$ is a hemicontinuous operator.  Together with monotonicity, we infer that $\CG$ is maximal monotone.  Hence, in order to prove the identification $\CG^* = \CG(v)$ it suffices to show that 
\begin{align}\label{IEQ:MINTY}
\int_{\Sigma_T} (\CG^* - \CG(w))(v - w) \dG \dt \geq 0 \quad \text{ for all } w \in L^r(\Sigma_T).
\end{align}
Let $w \in L^r(\Sigma_T)$ be arbitrary.  By the monotonicity of $G'$ we see that
\begin{equation}\label{IEQ:MG}
\begin{aligned}
0 & \leq \int_{\Sigma_T} (G'(v_N) - G'(w))(v_N - w) \dG \dt \\
& = \int_{\Sigma_T} \muGN (v_N - \lv_N) + \muGN \lv_N - \tfrac{\alpha^2}{K} \abs{v_N}^2 - \tfrac{\alpha}{K}( \beta - u_N) v_N  \dG \dt \\
& \quad - \int_{\Sigma_T} G'(v_N) w + G'(w) (v_N - w) \dG \dt.
\end{aligned}
\end{equation}
In light of the convergence results stated in Step 5, as well as 
\begin{align*}
v_N - \lv_N \to 0 \text{ strongly in } L^2(0,T;H^1(\Gamma)'), \quad \limsup_{N \to \infty} - \alpha^2 \norm{v_N}_{L^2(\Sigma_T)}^2 \leq - \alpha^2 \norm{v}_{L^2(\Sigma_T)}^2
\end{align*}
which follows from \eqref{DIFF} and weak lower semicontinuity, respectively, we infer from taking the limit superior in \eqref{IEQ:MG} the inequality
\begin{align*}
0 \leq \int_{\Sigma_T} \muG v - \tfrac{\alpha}{K} (\alpha v + \beta - u) v - \CG^* w + G'(w)(v - w) \dG \dt.
\end{align*}
Using the identity \eqref{EQ:G*} yields \eqref{IEQ:MINTY} and in particular, $G'(v_N) \to G'(v)$ weakly in $L^{r^*}(\Sigma_T)$.  Thus, we have derived item (ii) of Theorem \ref{thm:CHLW}.

For item (iii), proceeding similarly to the approach in Step 2 we find that
\begin{gather*}
\intO F(u(t)) \dx \le \underset{N\to \infty}{\lim\inf}\;\intO F(u_N(t)) \dx,\quad \intG G(v(t)) \dG \le \underset{N\to \infty}{\lim\inf}\;\intO G(v_N(t)) \dG,\\
\intG \abs{u(t)-H(v(t))}^2 \dG \le \underset{N\to \infty}{\lim\inf}\; \intG \abs{u_N(t)-H(v_N(t))}^2 \dG.
\end{gather*}
for almost all $t\in [0,T]$. These results can be used to prove the energy inequality \eqref{DEF:CHLW5} by taking the limit inferior in \eqref{IEQ:EUN}.  

\paragraph{Step 7: Uniqueness.}
Let $(u_1, v_1, \mu_1, \mu_{\Gamma,1})$ and $(u_2, v_2, \mu_2, \mu_{\Gamma,2})$ denote weak two solutions to \eqref{CHLW} corresponding to the same initial data. We define the difference of these solutions as
\begin{align*}
 (\lu,\lv,\lmu,\lmuG) := (u_1,v_1,\mu_1,\mu_{\Gamma,1}) -(u_2,v_2,\mu_2,\mu_{\Gamma,2}).
\end{align*}
For convenience, we use the notation
\begin{align*}
\bar{h} := H(v_1) - H(v_2), \quad \bar{h}' := H'(v_1) - H'(v_2).
\end{align*}
For arbitrary $t_0 \in (0,T]$, test functions $\xi \in L^2(0,T;H^1(\Omega))$ and $\varphi \in L^2(0,T;H^1(\Gamma))$, we define $Q_{t_0} = \Omega \times (0, t_0)$, $\Sigma_{t_0} = \Gamma \times (0,t_0)$, 
\begin{align}\label{uniq:test}
\zeta(\cdot,t) :=	\begin{cases}
\int_t^{t_0} \xi (\cdot,s) \ds, &\text{ if } t\le t_0,\\
0 &\text{ if } t>t_0
\end{cases}
\quad \text{ and } \quad \theta(\cdot,t) :=
\begin{cases}
\int_t^{t_0} \varphi (\cdot,s) \ds, &\text{ if } t\le t_0,\\
0 &\text{ if } t>t_0.
\end{cases}
\end{align}
Then, $\zeta \in L^2(0,T;H^1(\Omega)) \cap H^1(0,T;L^2(\Omega))$ and $\theta \in L^2(0,T;H^1(\Gamma)) \cap H^1(0,T;L^2(\Gamma))$.  From \eqref{DEF:CHLW1} with this particular $\zeta$, we have
\begin{align*}
\int_{Q_{t_0}} \lu \xi \dx \dt = \int_0^{t_0} \inn{\lu_t}{\zeta}_{H^1(\Omega)} \dt = - \int_{Q_{t_0}} \nabla \lmu \cdot \nabla \zeta \dx \dt = - \int_{Q_{t_0}} \nabla \left ( \int_{0}^t \lmu \ds \right ) \cdot \nabla \xi \dx \dt.
\end{align*}
This implies that we can write
\begin{align*}
\NN (\lu) = - \int_0^t \lmu \ds + c \quad \text{ and } \quad \pd_t \NN(\lu) = - \lmu
\end{align*}
for some constant $c \in \R$.  Choosing $\xi = \lmu$ in turn gives
\begin{align*}
\int_{Q_{t_0}} \lu \lmu \dx \dt = - \int_{Q_{t_0}} \nabla \NN(\lu) \cdot \nabla \pd_t \NN(\lu) \dx \dt = - \frac{1}{2} \norm{\nabla \NN(\lu(t_0))}_{L^2(\Omega)}^2,
\end{align*}
where we used $\lu(0) = 0$ and $\NN(0) = 0$.  Proceeding similarly, we infer that
\begin{align*}
\int_{\Sigma_{t_0}} \lv \lmuG \dG \dt = - \frac{1}{2} \norm{\gradg \NNG(\lv(t_0))}_{L^2(\Gamma)}^2,
\end{align*}
and thus,
\begin{align}\label{EQ:HOD}
\frac{1}{2} \norm{\lbu(t_0)}_{(\Ho)'}^2 + \int_{Q_{t_0}} \lu \lmu \dx \dt + \int_{\Sigma_{t_0}} \lv \lmuG \dG \dt = 0.
\end{align}
For any $M > 0$, we define the projection $\PM : \R \to \R$ as
\begin{align*}
\PM(s) = \begin{cases}
s & \text{ if } \abs{s} < M, \\
\frac{s}{\abs{s}} M & \text{ if } \abs{s} \geq M
\end{cases}
\end{align*}
satisfying $\abs{\PM(s)} \leq \abs{s}$ and $\PM'(s) = \chi_{[-M,M]}(s)$, where $\chi_A$ denotes the characteristic function of the set $A$.  Then, in \eqref{DEF:CHLW3} we can consider the test function $\eta = \chi_{[0,t_0]} \PM(\lu) \in L^2(0,T;H^1(\Omega)) \cap L^{\infty}(Q_T)$.  Monotonicity of $F_1'$ yields
\begin{align*}
\int_{Q_{t_0}} \nabla \lu \cdot \nabla \eta + (F_2'(u_1) - F_2'(u_2)) \eta - \lmu \eta \dx \dt - \int_{\Sigma_{t_0}} K^{-1}(\bar{h} - \lu) \eta \dG \dt \leq 0.
\end{align*}
By the dominated convergence theorem, we can pass to the limit $M \to \infty$ rigorously and obtain an analogous inequality where $\eta$ is replaced by $\lu$.  Lipschitz continuity of $F_2'$ then yields
\begin{align}\label{uniq:1}
\int_{Q_{t_0}} \abs{\nabla \lu}^{2} - \lmu \lu \dx \dt - \int_{\Sigma_{t_0}} K^{-1}(\bar{h} - \lu) \lu \dG \dt \leq C \norm{\lu}_{L^2(Q_{t_0})}^2.
\end{align}
For the case $\kappa > 0$, proceeding similarly, in \eqref{DEF:CHLW4} we consider the test function $\psi = \chi_{[0,t_0]} \PM(\lv) \in L^2(0,T;\Xk) \cap L^{\infty}(\Sigma_T)$ and obtain
\begin{equation}\label{uniq:2}
\begin{aligned}
& \int_{\Sigma_{t_0}} \kappa \abs{\gradg \lv}^2 - \lmuG \lv + K^{-1}H'(v_1)(\bar{h}- \lu) \lv \dG \dt \\
& \quad + \int_{\Sigma_{t_0}} K^{-1}(H(v_2) - u_2) \bar{h}' \lv \dG \leq C \norm{\lv}_{L^2(\Sigma_{t_0})}^2.
\end{aligned}
\end{equation}
Adding these two inequalities to \eqref{EQ:HOD} leads to
\begin{align*}
& \frac{1}{2} \norm{\lbu(t_0)}_{(\Ho)'}^2 + \norm{\nabla \lu}_{L^2(Q_{t_0})}^2 + \kappa \norm{\gradg \lv}_{L^2(\Sigma_{t_0})}^2 + K^{-1} \norm{\lu}_{L^2(\Sigma_{t_0})}^2 - C \big (\norm{\lv}_{L^2(\Sigma_{t_0})}^2 + \norm{\lu}_{L^2(Q_{t_0})}^2 \big )
 \\
& \quad \leq \int_{\Sigma_{t_0}} K^{-1}[\bar{h}(\lu - H'(v_1) \lv) + H'(v_1) \lu \lv - (H(v_2) - u_2)\bar{h}' \lv]  \dG \dt  \\
& \quad \leq C \sum_{i=1}^2 \int_{\Sigma_{t_0}}(1 + \abs{v_i}^{q}) \abs{\lv} \abs{\lu} + (1+\abs{v_i}^{2q}) \abs{\lv}^2 + (1+ \abs{v_i}^q + \abs{u_2})(1+ \abs{v_i}^{q-1}) \abs{\lv}^2 \dG \dt \\
& \quad \leq \frac{1}{2K} \norm{\lu}_{L^2(\Sigma_{t_0})}^2 + C  \norm{\lv}_{L^2(0,t_0;L^4(\Gamma))}^2,
\end{align*}
where we have used the embeddings $H^1(\Omega)\emb L^4(\Gamma)$ and $H^1(\Gamma) \emb L^r(\Gamma)$ for all $r < \infty$ and the regularities of the solutions.  By the Gagliardo--Nirenberg inequality in two dimensions we find that
\begin{align*}
\norm{\lv}_{L^4(\Gamma)}^2 \leq C \big ( \norm{\gradg \lv}_{L^2(\Gamma)} \norm{\lv}_{L^2(\Gamma)} + \norm{\lv}_{L^2(\Gamma)}^2 \big ) .
\end{align*}
Hence, applying Young's inequality to absorb all terms involving $\gradg \lv$ by the left-hand side, and using \eqref{SolnOp:b}, \eqref{SolnOp:s} and \eqref{defn:H0':inn}, we arrive at 
\begin{align*}
& \frac{1}{2} \norm{\lbu(t_0)}_{(\Ho)'}^2 + \norm{\nabla \lu}_{L^2(Q_{t_0})}^2 + \frac{3\kappa}{4} \norm{\gradg \lv}_{L^2(\Sigma_{t_0})}^2 + \frac{1}{2K} \norm{\lu}_{L^2(\Sigma_{t_0})}^2 \\
& \quad \leq C \big (\norm{\lu}_{L^2(Q_{t_0})}^2 + \norm{\lv}_{L^2(\Sigma_{t_0})}^2 \big ) \\
& \quad \leq \frac{\kappa}{4} \norm{\gradg \lv}_{L^2(\Sigma_{t_0})}^2 + \frac{1}{2} \norm{\nabla \lu}_{L^2(Q_{t_0})}^2 + C \int_{0}^{t_0} \norm{\lbu}_{(\Ho)'}^2 \dt.
\end{align*}
Uniqueness thus follows from Gronwall's inequality.

For the case $\kappa = 0$, recall that $H(s) = \alpha s + \beta$ is affine linear and that $G = G_1$ is convex.  Thus, the right-hand side of \eqref{uniq:2} arising from the Lipschitz continuity of $G_2'$ is now zero.  Adding \eqref{uniq:1} and \eqref{uniq:2} to \eqref{EQ:HOD} now gives
\begin{align*}
& \frac{1}{2} \norm{\lbu(t_0)}_{(\Ho)'}^2 + \norm{\nabla \lu}_{L^2(Q_{t_0})}^2 + K^{-1} \norm{\lu}_{L^2(\Sigma_{t_0})}^2 + \alpha^2 K^{-1} \norm{\lv}_{L^2(\Sigma_{t_0})}^2 - C \norm{\lu}_{L^2(Q_{t_0})}^2 \\
& \quad \leq \int_{\Sigma_{t_0}} K^{-1} [2 \alpha \lv \lu] \dG \dt \leq \frac{\alpha^2}{2K} \norm{\lv}_{L^2(\Sigma_{t_0})}^2 + \frac{2}{K} \norm{\lu}_{L^2(\Sigma_{t_0})}^2,
\end{align*}
owning to the fact $\bar{h} = \alpha \lv$ and $\bar{h}' = 0$.  By the trace inequality in \eqref{pre:trace}, Young's inequality, and arguing as in the above case for terms involving the $L^2(Q_{t_0})$-norm of $\lu$, we see that
\begin{align*}
& \frac{1}{2} \norm{\lbu(t_0)}_{(\Ho)'}^2 + \frac{1}{2} \norm{\nabla \lu}_{L^2(Q_{t_0})}^2 + \frac{\alpha^2}{2K} \norm{\lv}_{L^2(\Sigma_{t_0})}^2 \leq C \int_0^{t_0} \norm{\lbu}_{(\Ho)'}^2 \dt.
\end{align*}
Hence, uniqueness follows from Gronwall's inequality.  This completes the proof of Theorem \ref{thm:CHLW}.

\section{Proof of Theorem \ref{thm:conv}}\label{sec:conv}
In this section, we take $H(s) = \alpha s + \beta$ for some $\alpha \neq 0$ and $\beta \in \R$.  The letter $C$ will denote generic positive constants independent of $K$ that may change their value from line to line. Moreover, the symbol $\Xk$ will denote the space introduced in \eqref{DEF:XY}.

\paragraph{Weak convergence.} We first address the convergence of the initial data.  From \eqref{conv:weak:ini}, \eqref{ass:pot:5:k} and the Poincar\'e inequality we see that 
\begin{align*}
\norm{\K u_0}_{H^1(\Omega)} \leq C, \quad  \norm{\K v_0}_{\Xk} \leq C, \quad \frac{1}{K}\norm{H(\K v_0) - \K u_0}_{L^2(\Gamma)}^2 \leq C.
\end{align*}
Hence, there exists a nonrelabelled subsequence such that
\begin{align*}
& \K u_0 \wto u_0 \text{ in } H^1(\Omega), \quad \K u_0 \to u_0 \text{ in } L^2(\Omega) \text{ and a.e.~in } \Omega, \\ 
& \K v_0 \wto v_0 \text{ in } \Xk, \quad \K u_0 \to u_0 \text{ in } L^2(\Gamma) \text{ and a.e.~on } \Gamma, \\
&H(\K v_0) - \K u_0 \to 0 \text{ in } L^2(\Gamma) \text{ and a.e.~on } \Gamma.
\end{align*}
The last compactness statement and affine linearity of $H$ imply that $\K v_0 \to \alpha^{-1}(u_0 - \beta)$ strongly in $L^2(\Gamma)$ and a.e.~on $\Gamma$, and so $u_0 \vert_{\Gamma} \in \Xk$ by uniqueness of strong and weak limits.  In particular, we have that the limit of the initial data satisfies $u_0 \in \V_m^{\kappa}$.

Similarly, from the energy inequality \eqref{DEF:CHLW5}, the assumption \eqref{conv:weak:ini}, and Poincar\'e's inequality, we immediate infer the following uniform estimates:
\begin{align*}
& \norm{\K u}_{L^{\infty}(0,T;H^1(\Omega))} + \norm{F(\K u)}_{L^{\infty}(0,T;L^1(\Omega))} + \norm{\nabla \K \mu}_{L^2(Q_T)}^2 + \norm{\K v}_{L^{\infty}(0,T;\Xk)}^2 \\
& \quad + \norm{G(\K v)}_{L^{\infty}(0,T;L^1(\Gamma))} + \frac{1}{K} \norm{H(\K v) - \K u}_{L^{\infty}(0,T;L^2(\Gamma))}^2  + \norm{\gradg \K \muG}_{L^2(\Sigma_T)}^2 \leq C.
\end{align*}
Arguing as in Step 3 of the proof of Theorem \ref{thm:CHLW}, we use a test function $\eta \in C^{\infty}_c(\Omega)$ with $0 \leq \eta \leq 1$ in the weak formulation \eqref{DEF:CHLW3} to obtain, in conjunction with \eqref{ass:pot:3} and the previous uniform estimates, that 
\begin{align*}
\abs{\int_\Omega \K \mu \eta \dx} \leq C \norm{\nabla \eta}_{L^2(\Omega)} + \int_\Omega \abs{F'(\K u)} \dx \leq C(\eta).
\end{align*}
Similar to the derivation of \eqref{EST:MUN}, we use the generalised Poincar\'e inequality (cf. \cite[p.~242]{Alt}) to infer that
\begin{align*}
\norm{\K \mu}_{L^2(Q)} \leq C.
\end{align*}
Consequently, choosing $\eta = \alpha$ in \eqref{DEF:CHLW3} we have 
\begin{align*}
\abs{\int_\Gamma \alpha K^{-1} (H(\K v) - \K u) \dG} \leq C \norm{\K \mu}_{L^2(\Omega)} + C \norm{F'(\K u)}_{L^1(\Omega)},
\end{align*}
and together with \eqref{DEF:CHLW4} this yields
\begin{align*}
\abs{\int_\Gamma \K \muG \dG} \leq \norm{G'(\K v)}_{L^1(\Gamma)} + C \norm{\K \mu}_{L^2(\Omega)} + C \norm{F'(\K u)}_{L^1(\Omega)}.
\end{align*}
As the right-hand side is bounded in $L^2(0,T)$ we deduce by the Poincar\'e inequality that
\begin{align*}
\norm{\K \muG}_{L^2(\Sigma_T)} \leq C.
\end{align*}
Lastly, the uniform boundedness of $\nabla \K \mu$ and $\gradg \K \muG$ implies that
\begin{align*}
\norm{\K u_t}_{L^2(0,T;H^1(\Omega)')} + \norm{\K v_t}_{L^2(0,T;H^1(\Gamma)')} \leq C.
\end{align*}
These uniform estimates allows us to deduce that, along a nonrelabelled subsequence, there exists a limit quadruplet $(u_*, v_*, \mu_*, \mu_{\Gamma,*})$ such that
\begin{align*}
\K u & \to u_* &&\text{ weakly* in } L^{\infty}(0,T;H^1(\Omega))\cap H^1(0,T;H^1(\Omega)')\\
\K u & \to u_* &&\text{ strongly in } C([0,T];L^2(\Omega)) \cap L^{\infty}(0,T;L^2(\Gamma)) \text{ and a.e.~in } \overline{Q_T}, \\
\K v & \to v_* &&\text{ weakly* in } L^{\infty}(0,T;\Xk) \cap H^1(0,T;H^1(\Gamma)'), \\
\K v & \to v_* &&\text{ strongly in } C([0,T];\Yk) \text{ and a.e.~on } \Sigma_T, \\
\K \mu & \to \mu_* &&\text{ weakly in } L^2(0,T;H^1(\Omega)), \\
\K \muG & \to \mu_{\Gamma,*} &&\text{ weakly in } L^2(0,T;H^1(\Gamma)), \\
H(\K v) - \K u & \to 0 &&\text{ strongly in } L^{2}(0,T;L^2(\Gamma)) \text{ and a.e.~on } \Sigma_T.
\end{align*} 
Here we point out that the a.e.~convergence of $\K v$ to $v_*$ on $\Sigma_T$ in the case $\kappa = 0$ is due to the fact that $H$ is affine linear and $\K u \to u_*$ a.e.~on $\Sigma_T$.  In particular, we have $H(\K v) \to u_*$ a.e.~on $\Sigma_T$ and hence $\K v \to v_* = \alpha^{-1}(u_* - \beta)$ a.e.~on $\Sigma_T$.  Taking arbitrary test functions $\zeta \in L^2(0,T;H^1(\Omega))$ in \eqref{DEF:CHLW1} and $\theta \in L^2(0,T;H^1(\Gamma))$ in \eqref{DEF:CHLW2}, integrating over $[0,T]$ and passing to the limit $K \to 0$ yields
\begin{align*}
0 & = \int_0^T \hspace{-3pt}\Big ( \inn{u_{*,t}}{\zeta}_{H^1(\Omega)} + \intO \nabla \mu_* \cdot \nabla \zeta \dx \Big ) \dt, \; 0  = \int_0^T\hspace{-3pt} \Big (\inn{v_{*,t}}{\theta}_{H^1(\Gamma)} + \intG \gradg \mu_{\Gamma,*} \cdot \gradg \theta \dG \Big ) \dt.
\end{align*}
We now take an arbitrary test function $\eta \in L^2(0,T;\Vk) \cap L^{\infty}({Q_T})$ with $\eta\vert_{\Sigma_T}\in L^\infty(\Sigma_T)$ in \eqref{DEF:CHLW3} and choose $\psi = \alpha^{-1} \eta \vert_{\Sigma_T}$ in \eqref{DEF:CHLW4}.  Adding these two equalities and integrating over $[0,T]$ leads to
\begin{align*}
0 & = \int_{Q_T} \nabla \K u \cdot \nabla \eta + F'(\K u) \eta - \K \mu \eta \dx \dt  + \int_{\Sigma_T} \frac{\kappa}{\alpha} \gradg \K v \cdot \gradg \eta + \frac{1}{\alpha} G'(\K v) \eta - \frac{1}{\alpha} \K \muG \eta \dG \dt.
\end{align*}
For the case $\kappa > 0$, we can follow the arguments in Step 6 of the proof of Theorem \ref{thm:CHLW} to pass to the limit $K \to 0$ in order to obtain
\begin{align*}
0 & = \int_{Q_T} \nabla u_* \cdot \nabla \eta + F'(u_*) \eta - \mu_* \eta \dx \dt  + \int_{\Sigma_T} \frac{\kappa}{\alpha} \gradg v_* \cdot \gradg \eta + \frac{1}{\alpha} G'(v_*) \eta - \frac{1}{\alpha} \mu_{\Gamma,*} \eta \dG \dt.
\end{align*}
Hence, substituting $\phi_* = \alpha \mu_{\Gamma,*}$ and $v_* = \alpha^{-1}(u_* - \beta)$ in the above equalities shows that the triplet $(u_*, \mu_*, \phi_*)$ is a weak solution to \eqref{lim:alt} in the sense of item (ii) of Theorem \ref{thm:lim}.  

For the case $\kappa = 0$, we can employ exactly the same argument.  Indeed, we have $\K v \to v_*$ a.e.~on $\Sigma_T$ and by Vitali's convergence theorem $G'(\K v) \to G'(v)$ strongly in $L^1(\Sigma_T)$.  Let us point out that in Step 6 of the proof of Theorem \ref{thm:CHLW}, we do not have a pointwise a.e~convergence of $v_N$ for the case $\kappa = 0$, and hence we have to use Minty's method to identify the limit of $G'(v_N)$. 

Lastly, to recovery the H\"older-in-time regularity of the limit functions, we fix arbitrary $t_1, t_2 \in [0,T]$ with $t_1 > t_2$, and consider the (time-independent) test function $\zeta = u_*(t_1) - u_*(t_2)$ in \eqref{DEF:CHLW1} for $(\K u, \K \mu)$. Then integrating over $[t_2, t_1]$ and passing to the limit $K \to 0$ gives
\begin{align*}
\norm{u_*(t_1) - u_*(t_2)}_{L^2(\Omega)}^2 & = - \int_{t_2}^{t_1} \int_\Omega \nabla \mu_* \cdot \nabla (u_*(t_1) -u_*(t_2)) \dx \dt \\
& \leq 2\norm{\nabla \mu_*}_{L^2(\Sigma_T)} \norm{u_*}_{L^{\infty}(0,T;H^1(\Omega))} \abs{t_1 - t_2}^{1/2} \leq C \abs{t_1 - t_2}^{1/2},
\end{align*}
which implies that $u_* \in C^{0,\frac{1}{4}}(0,T;L^2(\Omega))$.  For $\kappa > 0$, we use the same argument to obtain $\alpha^{-1}(u_* \vert_{\Sigma_T} - \beta) = v_* \in C^{0,\frac{1}{4}}(0,T;L^2(\Gamma))$.  Furthermore, for $\kappa = 0$, we take $\theta =  \NNG(v_*(t_1) - v_*(t_2))$ in \eqref{DEF:CHLW2} for $(\K v, \K \muG)$. Then integrating over $[t_2, t_1]$ and passing to the limit $K \to 0$ gives
\begin{align*}
\norm{v_*(t_1) - v_*(t_2)}_{H^1(\Gamma)'}^2 & = \inn{v_*(t_1) - v_*(t_2)}{\NNG(v_*(t_1) - v_*(t_2))}_{H^1(\Gamma)} \\
& = - \int_{t_2}^{t_1} \int_\Gamma \gradg \mu_{\Gamma,*} \cdot \gradg \NNG(v_*(t_1) - v_*(t_2)) \dG \dt\\
& \leq \norm{\gradg \mu_{\Gamma,*}}_{L^2(\Sigma_T)} \norm{v_*(t_1) - v_*(t_2)}_{H^1(\Gamma)'} \abs{t_1 - t_2}^{1/2}.
\end{align*}
Hence, $(u_*, \mu_*, \phi_*)$ fulfills also item (i) of Theorem \ref{thm:lim}.

\paragraph{Error estimates.} We reuse the notation $(\lu, \lv, \lmu, \lmuG)$ to denote the differences between $(\K u, \K v, \K \mu, \K \muG)$ and $(u, v, \mu, \muG)$, where $v := \alpha^{-1}(u \vert_{\Sigma_T} - \beta)$ and $\muG := \alpha^{-1} \phi$.  With the regularity assumed on $u$, which is guaranteed by Theorem \ref{thm:lim} under the condition \eqref{GK:f1}, we see that the weak formulation of \eqref{CHLW:lim} in item (ii) of Theorem \ref{thm:lim} can be expressed as
\begin{align*}
0 & = \inn{u_t}{\zeta}_{H^1(\Omega)} + \int_\Omega \nabla \mu \cdot \nabla \zeta \dx, \\
0 & = \inn{v_t}{\theta}_{H^1(\Gamma)} + \int_\Gamma \gradg \muG \cdot \gradg \theta \dG, \\
0 & = \int_\Omega \nabla u \cdot \nabla \eta + (F'(u) - \mu) \eta \dx - \int_\Gamma \pdnu u \eta \dG, \\
0 & = \int_\Gamma \kappa \gradg v \cdot \gradg \psi  + (G'(v) - \muG + \alpha \pdnu u) \psi \dG
\end{align*}
holding for all $\zeta \in H^1(\Omega)$, $\theta \in H^1(\Gamma)$, $\eta \in H^1(\Omega) \cap L^{\infty}(\Omega)$, $\psi \in \Xk \cap L^{\infty}(\Gamma)$, and for a.e.~$t \in (0,T)$.  Then, it holds that
\begin{align}
\label{conv:1} \inn{\lu_t}{\zeta}_{H^1(\Omega)} &= -\intO \grad \lmu \cdot \grad\zeta \dx,\\
\label{conv:2}
\inn{\lv_t}{\theta}_{H^1(\Gamma)} &= -\intG \gradg \lmuG \cdot \gradg \theta \dG,\\
\label{conv:3} 
\intO \grad \lu \cdot\grad \eta + (F'(\K u) - F'(u)) \eta -  \lmu \eta \dx &= \intG \Big ( K^{-1}\big(\alpha \lv - \lu \big) - \pdnu u  \Big )\eta  \dG,\\
\label{conv:4} \intG  \kappa \gradg \lv \cdot \gradg \psi + (G'(\K v) - G'(v))\psi - \lmuG \psi \dG& = - \intG \alpha \Big (K^{-1} \big( \alpha \lv - \lu \big) - \pdnu u \Big ) \psi  \dG,
\end{align}
where we have used the relation $u \vert_{\Sigma_T} = \alpha v + \beta$ to obtain the differences in the right-hand sides of \eqref{conv:3} and \eqref{conv:4}.  We note that the hypothesis \eqref{conv:weak:ini} implies $\mean{\lu(t)}_\Omega = 0$ and $\mean{\lv(t)}_\Gamma = 0$ for all $t \in [0,T]$.  Then, choosing $\zeta = \NN(\lu)$ and $\theta = \NN(\lv)$ yields by application of \eqref{pre:Nb},
\begin{align}\label{conv:est:-1}
\frac{1}{2} \frac{d}{dt} \norm{\lbu}_{(\Ho)'}^2 + \int_\Omega \lu \lmu \dx + \int_\Gamma \lv \lmuG \dG = 0.
\end{align}
As in the uniqueness proof of Theorem \ref{thm:CHLW}, we consider $\eta = \PM(\lu)$ and $\psi = \PM(\lv)$ and employ the monotonicity of $F_1'$ and $G_1'$.  Then passing to the limit $M \to \infty$ in the resulting inequalities and use the Lipschitz continuity of $F_2'$ and $G_2'$ leads to
\begin{equation}\label{conv:est:0}
\begin{aligned}
& \int_\Omega \abs{\nabla \lu}^2 - \lmu \lu \dx + \int_\Gamma \kappa \abs{\gradg \lv}^2 - \lmuG \lv \dG + K^{-1} \abs{\alpha \lv - \lu}^2 \dG \\
& \quad \leq \intG \pdnu u (\alpha \lv - \lu) \dG + C \big ( \norm{\lu}_{L^2(\Omega)}^2 + \norm{\lv}_{L^2(\Gamma)}^2 \big ) \\
& \quad \leq \intG \pdnu u (\alpha \lv - \lu) \dG + \frac{1}{2} \norm{\nabla \lu}_{L^2(\Omega)}^2 + \frac{\kappa}{2} \norm{\gradg \lv }_{L^2(\Gamma)}^2 + C \norm{\lbu}_{(\Ho)'}^2
\end{aligned}
\end{equation}
where for the last inequality we employed the interpolation estimates in \eqref{pre:Na}, as well as the identification of the $(\Ho)'$-norm in \eqref{pre:H}.  Adding \eqref{conv:est:0} to \eqref{conv:est:-1} leads to the differential inequality for the case $\kappa > 0$
\begin{equation}\label{conv:est:1}
\begin{aligned}
& \frac{d}{dt} \norm{\lbu}_{(\Ho)'}^2 + \norm{\nabla \lu}_{L^2(\Omega)}^2 + \kappa \norm{\gradg \lv}_{L^2(\Gamma)}^2 + K^{-1} \norm{\alpha \lv - \lu}_{L^2(\Gamma)}^2 \\
& \quad \leq CK \norm{\pdnu u}_{L^2(\Gamma)}^2 + C \norm{\lbu}_{(\Ho)'}^2,
\end{aligned}
\end{equation}
with a positive constant $C$ independent of $K$.  Then, the desired estimate \eqref{CHLW:rate} follows from first applying Gronwall's inequality to \eqref{conv:est:1} and then applying the interpolation inequalities from \eqref{pre:Na}.

For the case $\kappa = 0$, the convexity of $G$ from \eqref{ass:pot:5:k} implies that the term involving $\norm{\lv}_{L^2(\Gamma)}^2$ on the right-hand side of the first inequality in \eqref{conv:est:0} is not present.  Although the surface gradient term on the left-hand side of \eqref{conv:est:0} is lost, nevertheless we can proceed as before to infer the estimate \eqref{CHLW:rate} without the terms for $\K v - v$.

\section{Numerical approximation}\label{sec:num}
In this section we present fully discrete finite element approximations of the systems \eqref{CHLW} and \eqref{lim:alt} based on the implicit time discretisations in Section \ref{sec:CHLW} and Section \ref{sec:lim}, respectively. Moreover, we display several two-dimensional numerical simulations.  Let us mention that there are numerous contributions on the numerical approximation of the Cahn--Hilliard equation with Neumann/periodic boundary conditions (of which we only cite \cite{FP,FP2,GT,SXY,SY} and refer the reader to the references cited therein), and with dynamic boundary conditions \cite{CP,CPP,FYW,IMP,Trau}.  As the focus of this paper is on the analytical aspects, we leave the investigation of efficient numerical algorithms for \eqref{CHLW} to future research, and a recent contribution on the numerical analysis for the original model of \cite{LW} can be found in \cite{Metzger}.

\subsection{Finite element scheme}
Let $\Omega_h$ be a polygonal $(d=2)$ or a polyhedral $(d=3)$ approximation of the domain $\Omega$ with $\Gamma_h = \pd \Omega_h$.  The faces of $\Gamma_h$ are $(d-1)$-simplices with vertices lying on $\Gamma$ in the sense that $\Gamma_h$ is an interpolation of $\Gamma$.  We consider a quasi-uniform triangulation $\Th$ of $\Omega_h$ consisting of closed simplices $R$ with $h = \max \{ \mathrm{diam}(R) \, : \, R \in \Th\}$.  Associated to $\Th$ is the set of nodal points $\{x_i\}_{i=1}^{K_h}$ for $K_h \in \N$ and the set of continuous, piecewise linear nodal basis functions $\{\Lambda_i\}_{i=1}^{K_h}$ satisfying $\Lambda_i(x_j) = \delta_{ij}$ for $j \in \{1, \dots, K_h\}$.  Let $P_1(D)$ denote the space of affine linear polynomials on a closed, connected set $D\subset\R^d$.  Then, we define the finite element space
\begin{align*}
V^h := \mathrm{span} \{ \Lambda_1, \dots, \Lambda_{K_h}\} = \{  \Lambda \in C(\overline{\Omega_h}) \, : \, \Lambda \vert_{R} \in P_1(R) \, \forall R \in \Th \}.
\end{align*}
As observed in \cite[Sec.~2]{CP} and \cite[Sec.~3.2.1]{KL}, the set of restrictions $\{\Lambda_i \vert_{\Gamma_h}\}_{i=1}^{K_h}$ of the basis functions forms a basis for the space of piecewise linear polynomials on $\Gamma_h$, which allows us to define
\begin{align*}
V^h_\Gamma := \mathrm{span} \{ \Lambda_1 \vert_{\Gamma_h}, \dots, \Lambda_{K_h} \vert_{\Gamma_h} \} = \{ \Lambda \in C(\Gamma_h) \, : \, \Lambda \vert_{D} \in P_1(D) \; \forall   D = R \cap \Gamma_h \neq \emptyset \text{ s.t.}\, R \in \Th\}.
\end{align*}
We further introduce the characteristic function $\chi$ for the boundary nodal points as
\begin{align*}
\chi : \{ 1,\dots, K_h \} \to \{0,1\}, \quad \chi_k := \chi(k) = \begin{cases}
1 & \text{ if } x_k \in \Gamma_h, \\
0 & \text{ if } x_k \notin \Gamma_h,
\end{cases}
\end{align*}
and the Lagrange interpolation operators $I^h : C(\overline{\Omega_h}) \to V^h$ and $I^h_\Gamma : C(\Gamma_h) \to V^h_\Gamma$ as 
\begin{align*}
I^h \zeta = \sum_{i=1}^{K_h} \zeta(x_i) \Lambda_i, \quad I^h_\Gamma \xi = \sum_{i=1}^{K_h} \chi_i \xi(x_i) \Lambda_i \vert_{\Gamma_h} \quad \text{ for all } \zeta \in C(\overline{\Omega_h}), \, \xi \in C(\Gamma_h).
\end{align*}
For an arbitrary but fixed $T > 0$, let $N \in \N$ denote the number of time steps and $\tau := T/N >0$ the step size.  Based on the implicit time discretisation in Section \ref{sec:CHLW}, we propose the following fully discrete finite element scheme: Given $u_h^k \in V^h$ and $v_h^k \in V_\Gamma^h$, find $u_h^{k+1}, \mu_h^{k+1} \in V^h$ and $v_h^{k+1}, \mu_{\Gamma,h}^{k+1} \in V_\Gamma^h$ satisfying
\begin{subequations}\label{FE}
\begin{alignat}{2}
0 & = (u_h^{k+1} - u_h^k, \zeta)_\Omega^h + \tau (\nabla \mu_h^{k+1}, \nabla \zeta)_\Omega, \\
0 & = (v_h^{k+1} - v_h^k, \theta)_\Gamma^h + \tau (\gradg \mu_{\Gamma,h}^{k+1}, \gradg \theta)_\Gamma, \\
0 & = \eps (\nabla u_h^{k+1}, \nabla \eta)_\Omega + \tfrac{1}{\eps} (F'(u_h^{k+1}), \eta)_\Omega^h - (\mu_h^{k+1}, \eta)_\Omega^h - \tfrac{1}{K} ((H(v_h^{k+1}) - u_h^{k+1}), \eta)_\Gamma^h, \\
0 & = \kappa \delta (\gradg v_h^{k+1}, \gradg \psi)_\Gamma + \tfrac{1}{\delta} (G'(v_h^{k+1}), \psi)_\Gamma^h - (\mu_{\Gamma,h}^{k+1}, \psi)_\Gamma^h \\
\notag & \quad + \tfrac{1}{K} ((H(v_h^{k+1}) - u_h^{k+1}) H'(v_h^{k+1}), \psi)_\Gamma^h,
\end{alignat}
\end{subequations}
for all $\zeta, \theta \in V^h$ and $\eta, \psi \in V_\Gamma^h$.  In the above, $(\cdot, \cdot)_\Omega$ and $(\cdot, \cdot)_\Gamma$ denote the $L^2$-inner product on $\Omega_h$ and on $\Gamma_h$, respectively.  Meanwhile, $(\cdot, \cdot)_\Omega^h$ and $(\cdot, \cdot)_\Gamma^h$ are the lumped integrations defined as
\begin{align*}
(f, g)_\Omega^h := \int_{\Omega_h} I^h(fg) \dx \quad \forall f, g \in C(\overline{\Omega_h}),  \quad 
(p,q)_{\Gamma}^h := \int_{\Gamma_h} I^h_\Gamma(pq) \dG \quad \forall p,q \in C(\Gamma_h).
\end{align*}
Furthermore, for $k = 0$ we set $u_h^0 = \Pi^h u_0$ where $\Pi^h$ is the $L^2$-projection onto $V^h$ and analogously, $v_h^0 = \Pi_\Gamma^h v_0$ where $\Pi_\Gamma^h$ is the $L^2$-projection onto $V_\Gamma^h$.  The use of lumped integration in \eqref{FE} leads to diagonal mass matrices:
\begin{align*}
M_{ij} = (\Lambda_i, \Lambda_j)_{\Omega}^h = \mathrm{diag} \Big ( (1, \Lambda_i)_\Omega \Big )_{i=1, \dots, K_h}, \; M^\Gamma_{ij} = (\chi_i \Lambda_i , \Lambda_j)_{\Gamma}^h = \mathrm{diag} \Big ( \chi_i (1,\Lambda_i )_\Gamma \Big )_{i=1, \dots, K_h}.
\end{align*}
Together with the symmetric stiffness matrices 
\begin{align*}
A_{ij} = ((\nabla \Lambda_i, \nabla \Lambda_j)_\Omega)_{i,j=1, \dots, K_h}, \quad A^\Gamma_{ij} = ((\gradg \Lambda_i, \gradg \Lambda_j)_\Gamma)_{i,j=1, \dots, K_h}
\end{align*}
and the vectors of nonlinearities (using the Einstein summation convention in the arguments)
\begin{align*}
\mathcal{F}(U)_i & = M_{ii} F' ( U_j \Lambda_j  )(x_i), \quad \mathcal{G}(V)_i = M^\Gamma_{ii} G' ( V_j \Lambda_j  )(x_i), \\
\mathcal{H}(V)_i & = M^\Gamma_{ii} H ( V_j \Lambda_j  )(x_i), \quad \mathcal{J}(U,V)_i = M^\Gamma_{ii}  \Big [ H'  (  V_j \Lambda_j  ) \Big ( H (V_j \Lambda_j ) - U_j \Lambda_j \Big ) \Big ](x_i)
\end{align*}
for $U = (U_1, \dots, U_{K_h})^{\top}$ and $V = (V_1, \dots, V_{K_h})^{\top}$
we find that \eqref{FE} with $\zeta = \eta = \Lambda_j$ and $\phi = \psi = \Lambda_j \vert_{\Gamma_h}$ leads to the following system of nonlinear equations to solve
\begin{align}\label{discrete:FE}
\bm{L} \begin{pmatrix} 
U^{k+1} \\[1ex] V^{k+1} \\[1ex]\Xi^{k+1} \\[1ex] \Phi^{k+1}
\end{pmatrix} := \begin{pmatrix} 
M U^{k+1} + \tau A \Xi^{k+1} - M U^k \\[1ex]
M^\Gamma V^{k+1} + \tau A^\Gamma \Phi^{k+1} - M^\Gamma V^k \\[1ex]
\eps A U^{k+1} + \eps^{-1} \mathcal{F}(U^{k+1}) - M \Xi^{k+1} - \tfrac{1}{K} \mathcal{H}(V^{k+1}) + \tfrac{1}{K} M^\Gamma U^{k+1} \\[1ex]
\kappa \delta A^\Gamma V^{k+1} + \delta^{-1} \mathcal{G}(V^{k+1}) - M^\Gamma \Phi^{k+1} + \tfrac{1}{K} \mathcal{J}(U^{k+1}, V^{k+1})
\end{pmatrix} = \begin{pmatrix} 0 \\[1ex] 0 \\[1ex] 0 \\[1ex] 0 \end{pmatrix}
\end{align}
for our finite element solutions
\begin{align*}
u_h^{k+1} = \sum_{i=1}^{K_h} U_i^{k+1} \Lambda_i, \; \, v_h^{k+1} = \sum_{i=1}^{K_h} \chi_i V_i^{k+1} \Lambda_i \vert_{\Gamma_h}, \; \,  \mu_h^{k+1} = \sum_{i=1}^{K_h} \Xi_i^{k+1} \Lambda_i, \; \,  \mu_{\Gamma,h}^{k+1} = \sum_{i=1}^{K_h} \chi_i \Phi_i^{k+1} \Lambda_i \vert_{\Gamma_h}.
\end{align*}
For the case $H(s) = \alpha s + \beta$, we have
\begin{align*}
\mathcal{H}(V) = \alpha M^\Gamma V + \beta M^\Gamma E, \quad \mathcal{J}(U,V) = \alpha M^\Gamma \Big ( \alpha V + \beta E - U \Big ),
\end{align*}
where $E = (1, \dots, 1)^{\top}$ is the vector with all entries equal to one, and the discrete nonlinear system \eqref{discrete:FE} simplifies considerably.

In view of the choice of our finite element basis functions, we can derive a finite element scheme for the limit system \eqref{lim:alt} based on the implicit time discretisation \eqref{EQ:DWFLIM}: Given $u_h^k \in V^h$, find $u_h^{k+1}, \mu_h^{k+1} \in V^h$ and $\phi_h^{k+1} \in V_\Gamma^h$ satisfying
\begin{subequations}\label{FE:lim}
\begin{alignat}{2}
0 & = (u_h^{k+1} - u_h^k, \zeta)_\Omega^h + \tau (\nabla \mu_h^{k+1}, \nabla \zeta)_\Omega, \\
0 & = (u_h^{k+1} - u_h^k, \theta)_\Gamma^h + \tau (\gradg \phi_h^{k+1}, \gradg \theta)_\Gamma, \\
0 & = \eps (\nabla u_h^{k+1}, \nabla \tilde \eta)_\Omega +(  \tfrac{1}{\eps} F'(u_h^{k+1}) -\mu_h^{k+1}, \tilde \eta)_\Omega^h + \tfrac{\kappa \delta}{\alpha^2} (\gradg u_h^{k+1}, \gradg \tilde \eta)_\Gamma \\
\notag & \quad +( \tfrac{1}{\alpha \delta} G' (\tfrac{1}{\alpha}(u_h^{k+1} - \beta)) - \tfrac{1}{\alpha^2} \phi_h^{k+1}, \tilde \eta)_\Gamma^h 
\end{alignat}
\end{subequations}
for all $\zeta, \tilde \eta \in V^h$ and $\theta \in V_\Gamma^h$.  Introducing the vector 
\begin{align*}
\widetilde{\mathcal{G}}(U)_i = M_{ii}^\Gamma G'( \tfrac{1}{\alpha}(U_j \Lambda_j - \beta) )(x_i)
\end{align*}
and recalling the matrices $M$, $M^\Gamma$, $A$ and $A^\Gamma$, as well as the vector $\mathcal{F}(U)$,  we find that \eqref{FE:lim} with $\zeta = \tilde \eta = \Lambda_j$ and $\theta = \Lambda_j \vert_{\Gamma_h}$ leads to the following system of nonlinear equations to solve
\begin{align}\label{discrete:FE:lim}
\begin{pmatrix} 
M U^{k+1} + \tau A \Xi^{k+1} - M U^k \\
M^\Gamma U^{k+1} + \tau A^\Gamma \Phi^{k+1} - M^\Gamma U^k \\
(\eps A + \tfrac{\kappa \delta}{\alpha^2} A^\Gamma) U^{k+1} + \tfrac{1}{\eps} \mathcal{F}(U^{k+1}) - M \Xi^{k+1} + \tfrac{1}{\alpha \delta} \widetilde{\mathcal{G}}(U^{k+1}) - \tfrac{1}{\alpha^2} M^\Gamma \Phi^{k+1}
\end{pmatrix} = \begin{pmatrix}
0 \\ 0 \\ 0
\end{pmatrix}
\end{align}
for our finite element solutions
\begin{align*}
u_h^{k+1} = \sum_{i=1}^{K_h} U_i^{k+1} \Lambda_i, \; \,  \mu_h^{k+1} = \sum_{i=1}^{K_h} \Xi_i^{k+1} \Lambda_i, \; \,  \phi_h^{k+1} = \sum_{i=1}^{K_h} \chi_i \Phi_i^{k+1} \Lambda_i \vert_{\Gamma_h}.
\end{align*}

The subsequent nonlinear schemes \eqref{discrete:FE} and \eqref{discrete:FE:lim} are solved using a Newton method with MATLAB built on previous work by D.~Trautwein \cite{Trau}.   In the following, we always take $\Omega = \Omega_h = [0,1] \times [0,1]$ with a standard Friedrichs--Keller triangulation and spatial step size $h=0.01$. This yields 101 degrees of freedom along each axis.

\subsection{Qualitative experiments}

\subsubsection{Dynamics of the limit system \eqref{CHLW:lim}}
In the first set of experiment, we demonstrate the difference in dynamics for the limit system \eqref{CHLW:lim} for the case $\alpha = 1, \beta = 0$, i.e., the original model proposed in \cite{LW}, and for nontrivial values of $\{\alpha, \beta \}$.  We begin with the initial condition
\begin{align*}
u_0(x,y) = \begin{cases} 1 & \text{ if } x > 1/2, \\ -1 & \text{ if } x \leq 1/2, \end{cases} 
\end{align*}
for $(x,y) \in [0,1]^2$.  It is clear that $\mean{u}_\Omega = 0$ and $\mean{u}_\Gamma = 0$.  

\begin{figure}[h]
\captionsetup[subfigure]{labelformat=empty}
\centering 
\subfloat[]{\includegraphics[width=0.4\textwidth]{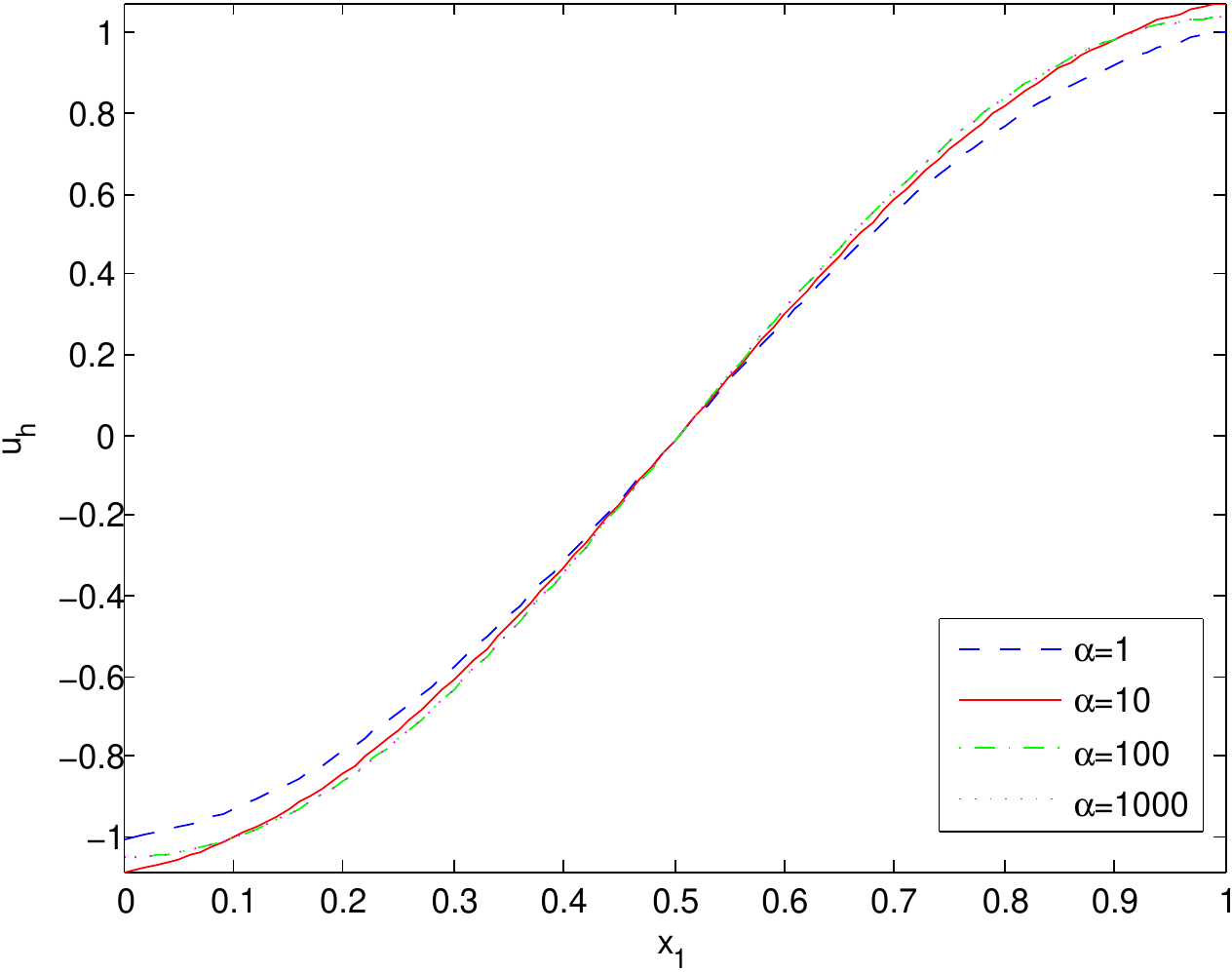}}
\hspace{2cm}
\subfloat[]{\includegraphics[width=0.4\textwidth]{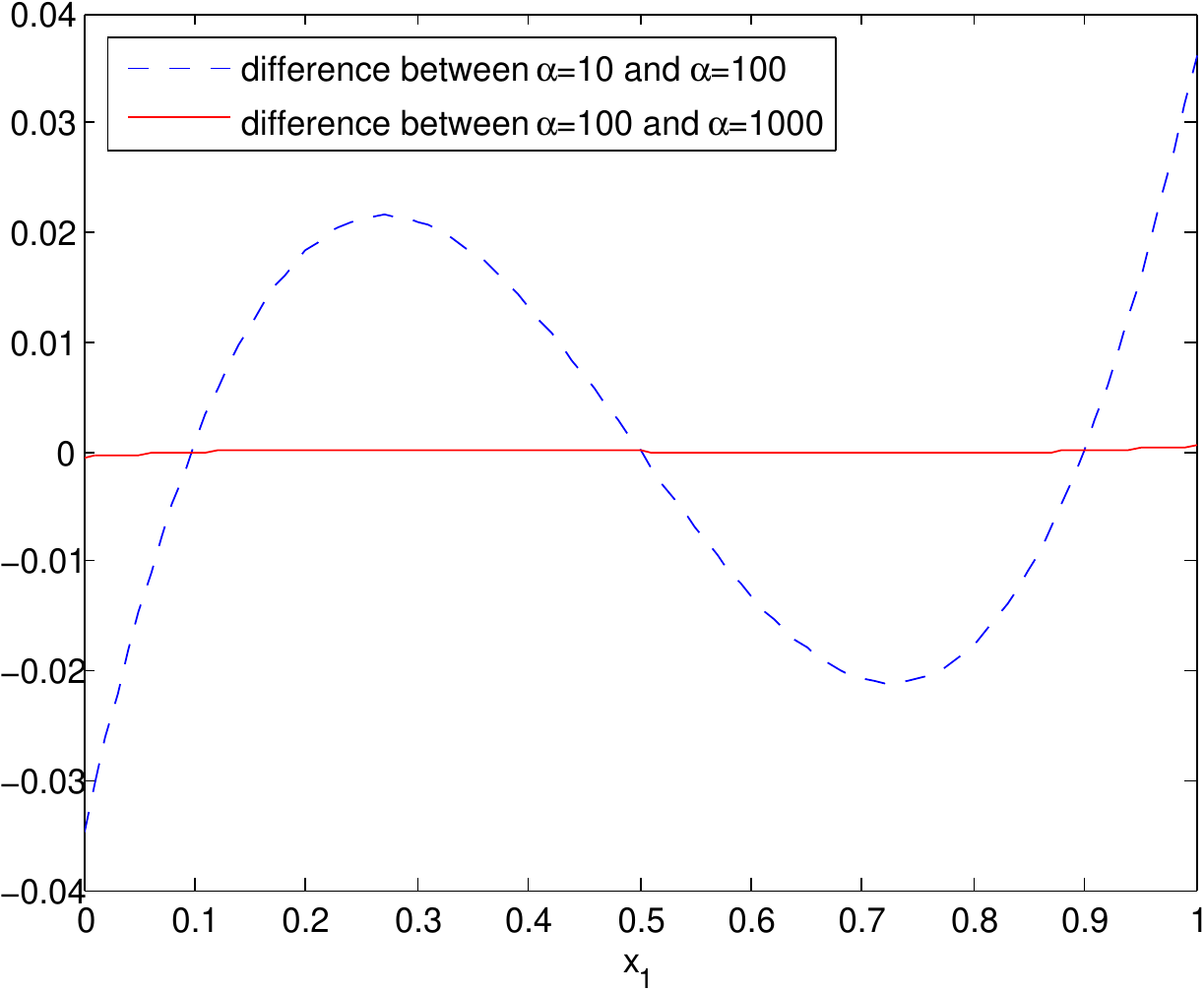}}
\caption{(Left) Profiles of the discrete solution $u_h^k$ to \eqref{CHLW:lim} at iteration $k = 200$ for $\alpha \in \{1,10,100,1000\}$ and $\beta = 0$.  (Right) The differences between the profiles for $\alpha \in \{10, 100\}$ and for $\alpha \in \{100,1000\}$.}
\label{fig:1:varyalpha}
\end{figure}

\paragraph{Effects of varying $\alpha$ with $\beta = 0$.}  We then solve \eqref{discrete:FE:lim} with the fixed parameters and potentials
\begin{align}\label{Ex1:para}
\tau = 10^{-5}, \quad \eps = 1, \quad \delta = 0.1, \quad \beta = 0, \quad \kappa = 1, \quad F(s) = G(s) = \tfrac{1}{4}(s^2 - 1)^2.
\end{align}
In Figure \ref{fig:1:varyalpha} we plot the projection of the discrete solution $u_h^k$ on the line $\{y = 0.5 \}$ at iteration $k = 200$ for $\alpha \in \{1, 10, 100, 1000\}$, as the solution seems to be constant in the orthogonal direction, as well as the differences between the profiles for $\alpha = 10$ and $\alpha = 100$, and for $\alpha = 100$ and $\alpha = 1000$.  Identical profiles are observed for the case of negative $\alpha$, and so we do not report them here.
It appears that the overall solution profile is rather robust with respect to the size of $\alpha$.  We also report that the magnitude of the difference between the profiles for $\alpha = 10$ and $\alpha = 100$ is of order $3 \times 10^{-2}$, while the magnitude of the difference between the profiles for $\alpha = 100$ and $\alpha = 1000$ is of order $5 \times 10^{-4}$.  Further comparison between $\alpha = 1000$ and $\alpha = 10000$ yields a difference of order $5 \times 10^{-6}$, which seems to suggest some form of convergence as $\alpha \to \infty$. 

\begin{figure}[h]
\captionsetup[subfigure]{labelformat=empty}
\centering
\subfloat[$\beta=0$]{\includegraphics[width=0.3\textwidth]{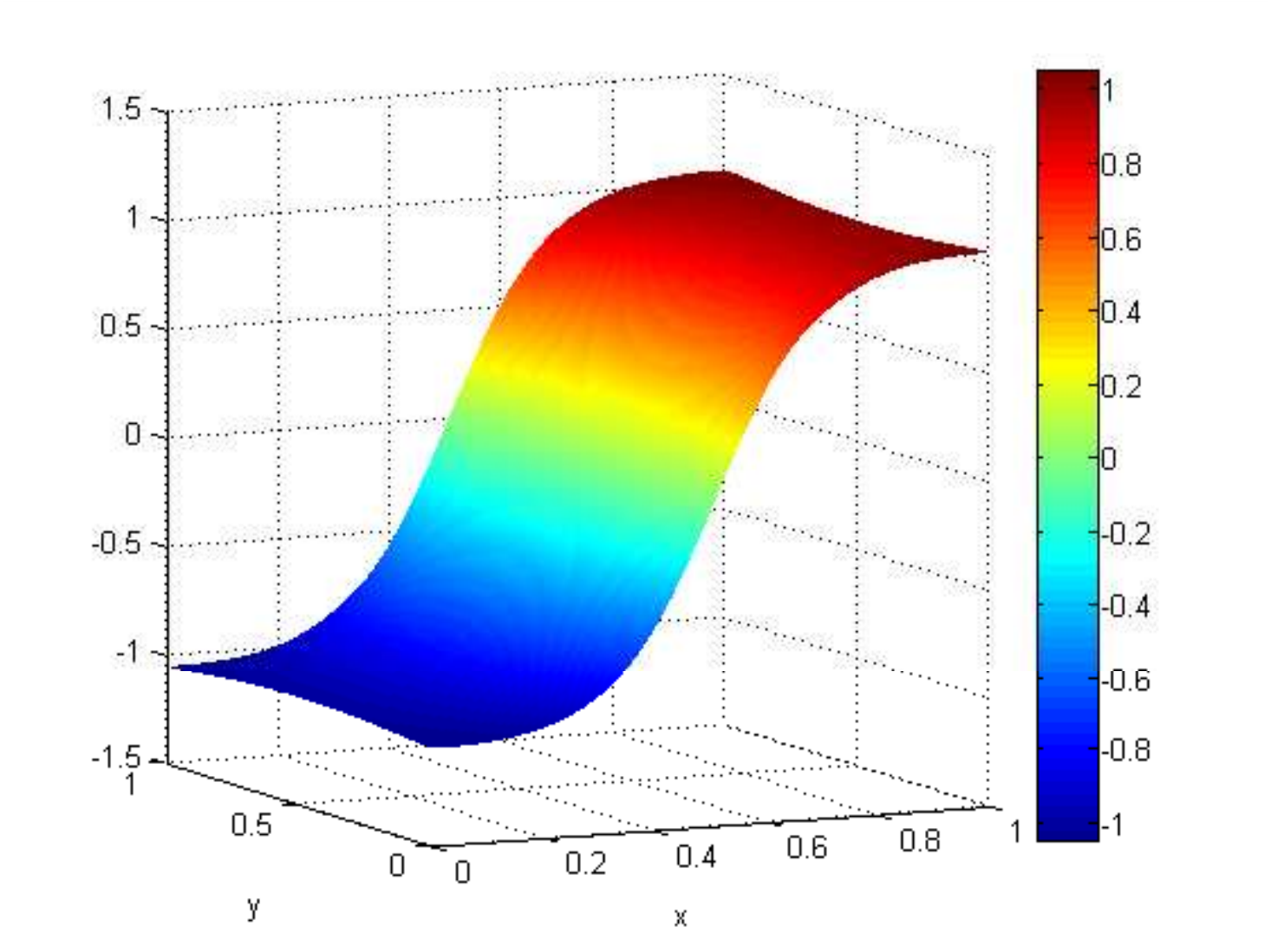}}
\hfill
\subfloat[$\beta=1.5$]{\includegraphics[width=0.3\textwidth]{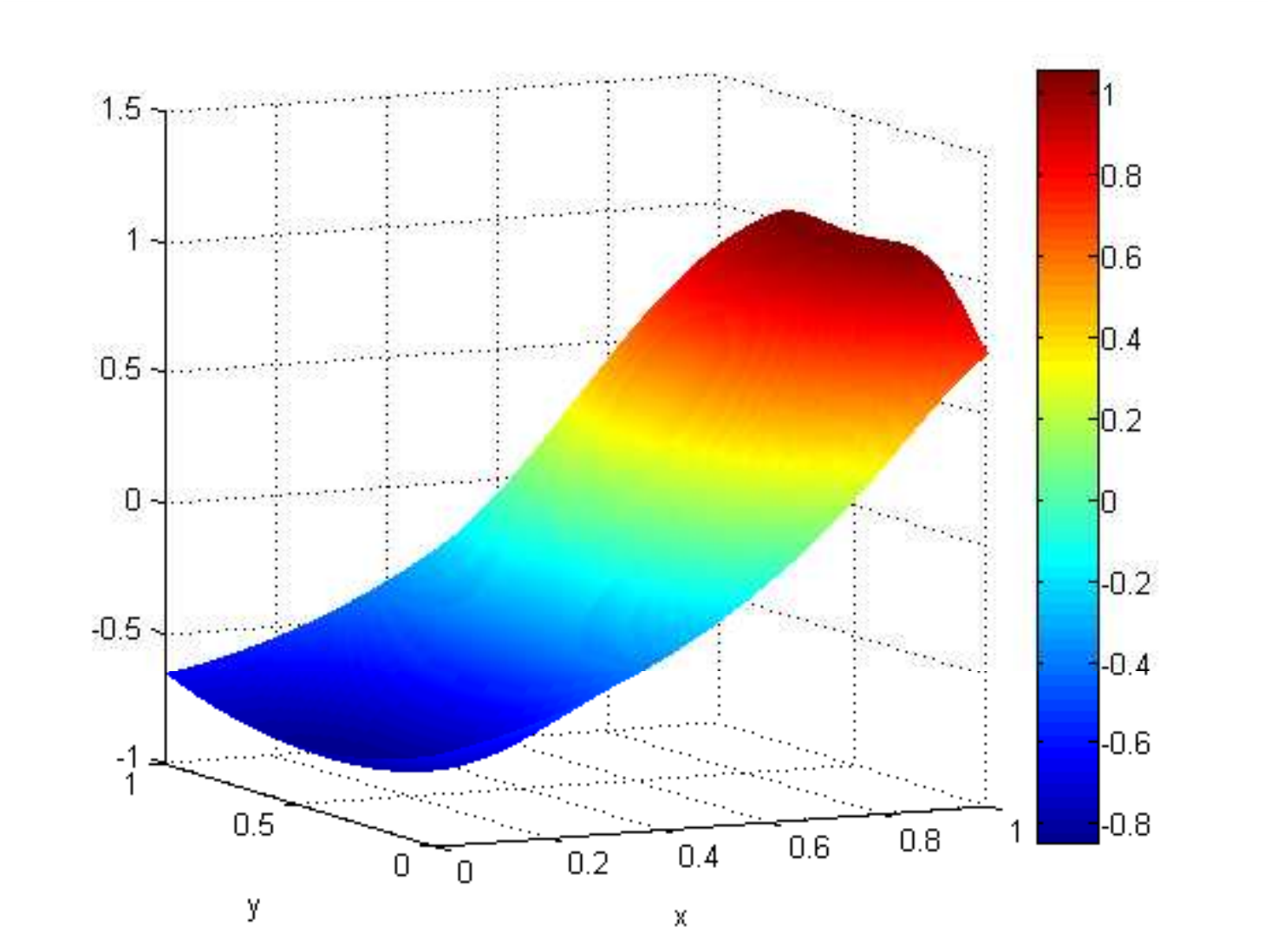}}
\hfill
\subfloat[$\beta=2$]{\includegraphics[width=0.3\textwidth]{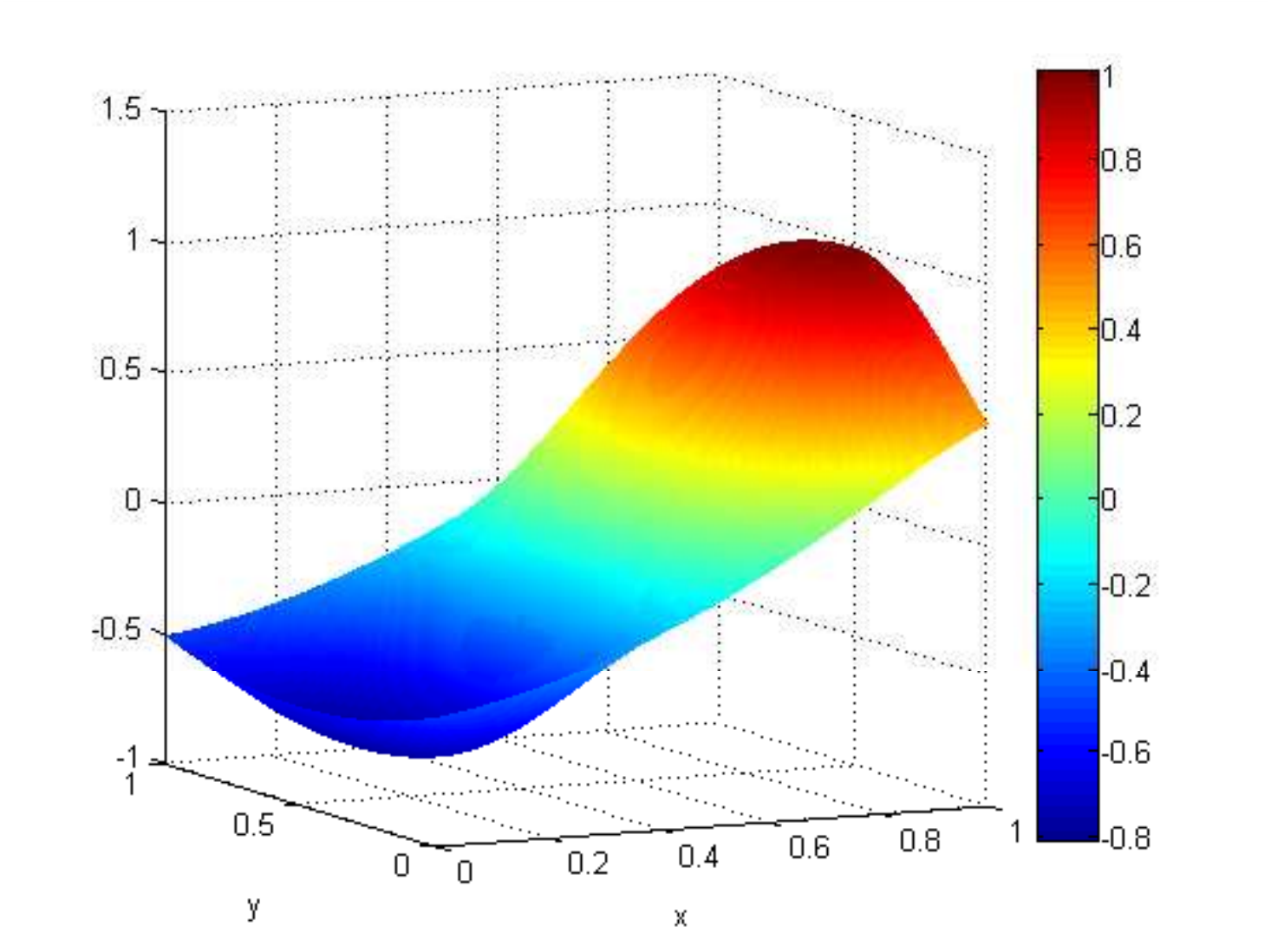}} \\
\subfloat[$\beta=3$]{\includegraphics[width=0.3\textwidth]{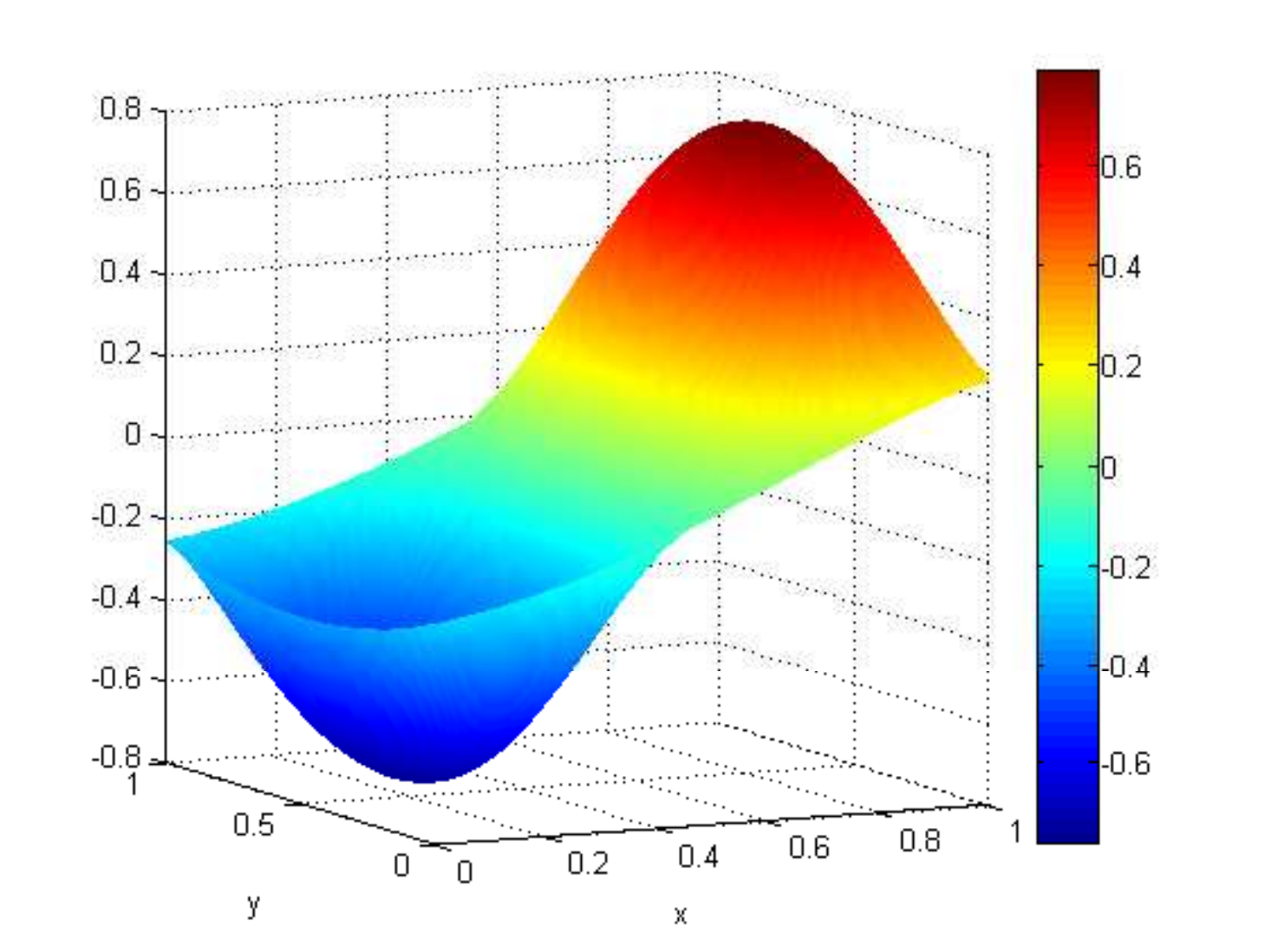}}
\hfill
\subfloat[$\beta=5$]{\includegraphics[width=0.3\textwidth]{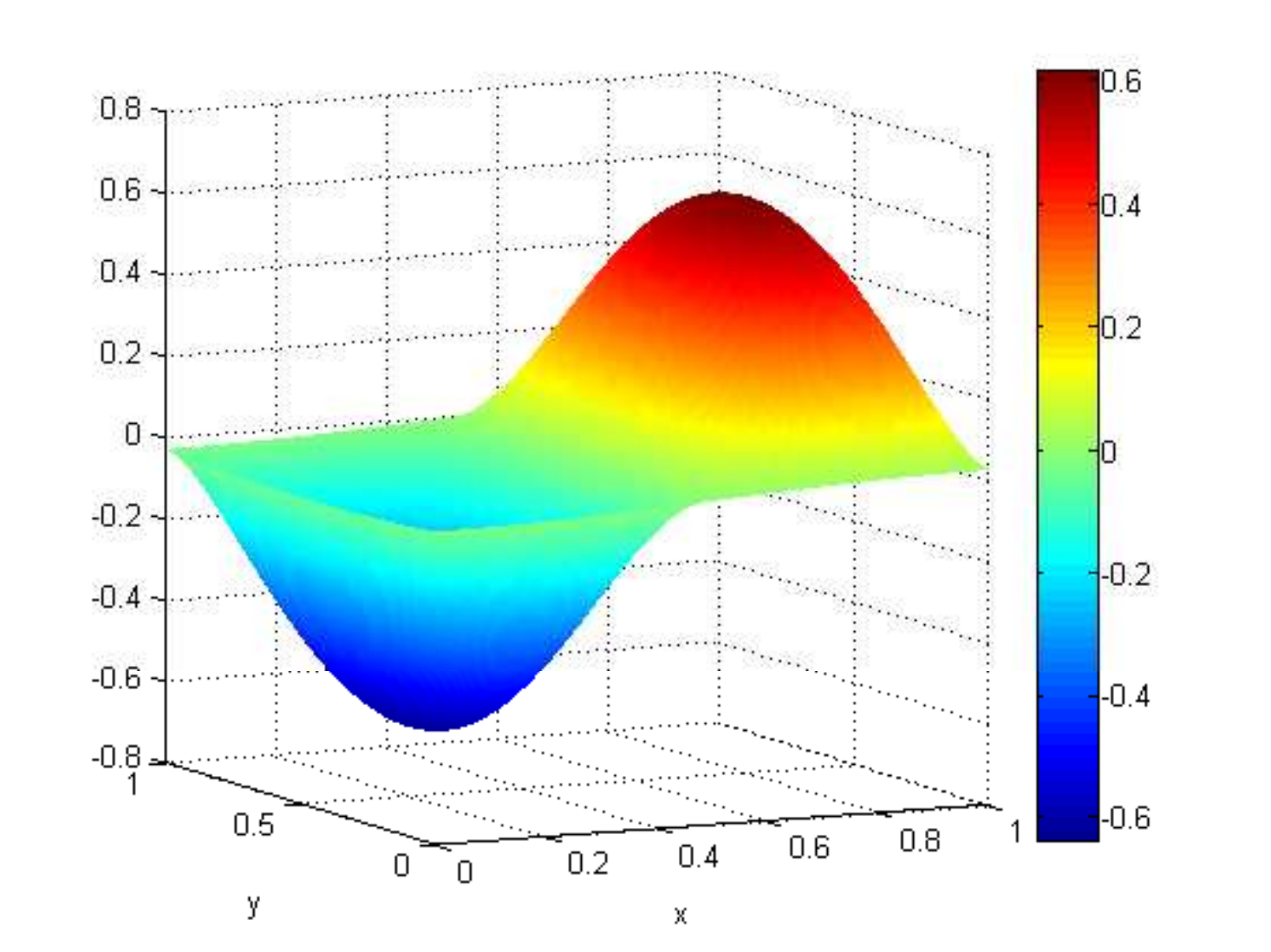}}
\hfill
\subfloat[$\beta=100$]{\includegraphics[width=0.3\textwidth]{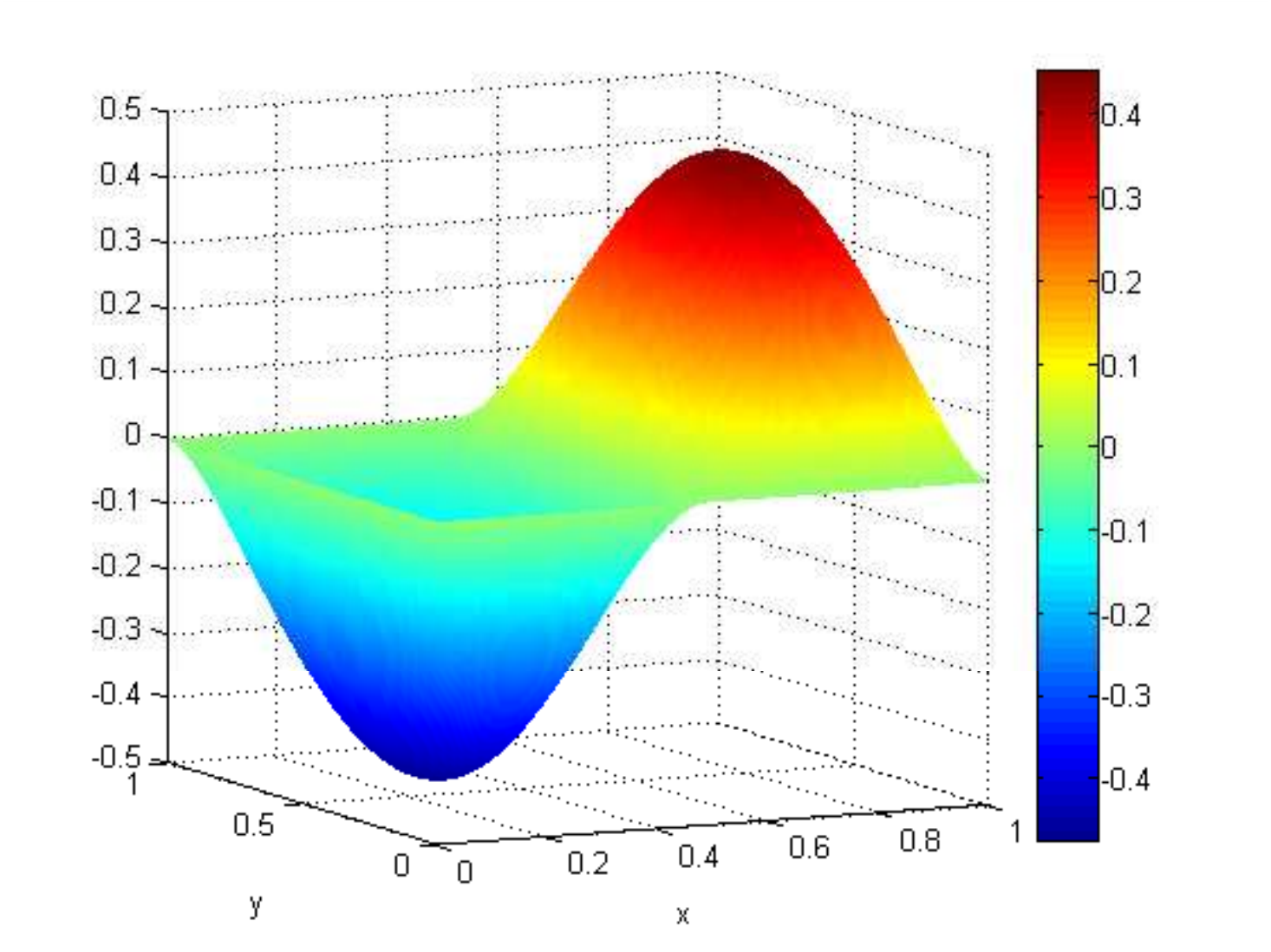}} \\
\subfloat[$\beta=3$]{\includegraphics[width=0.3\textwidth]{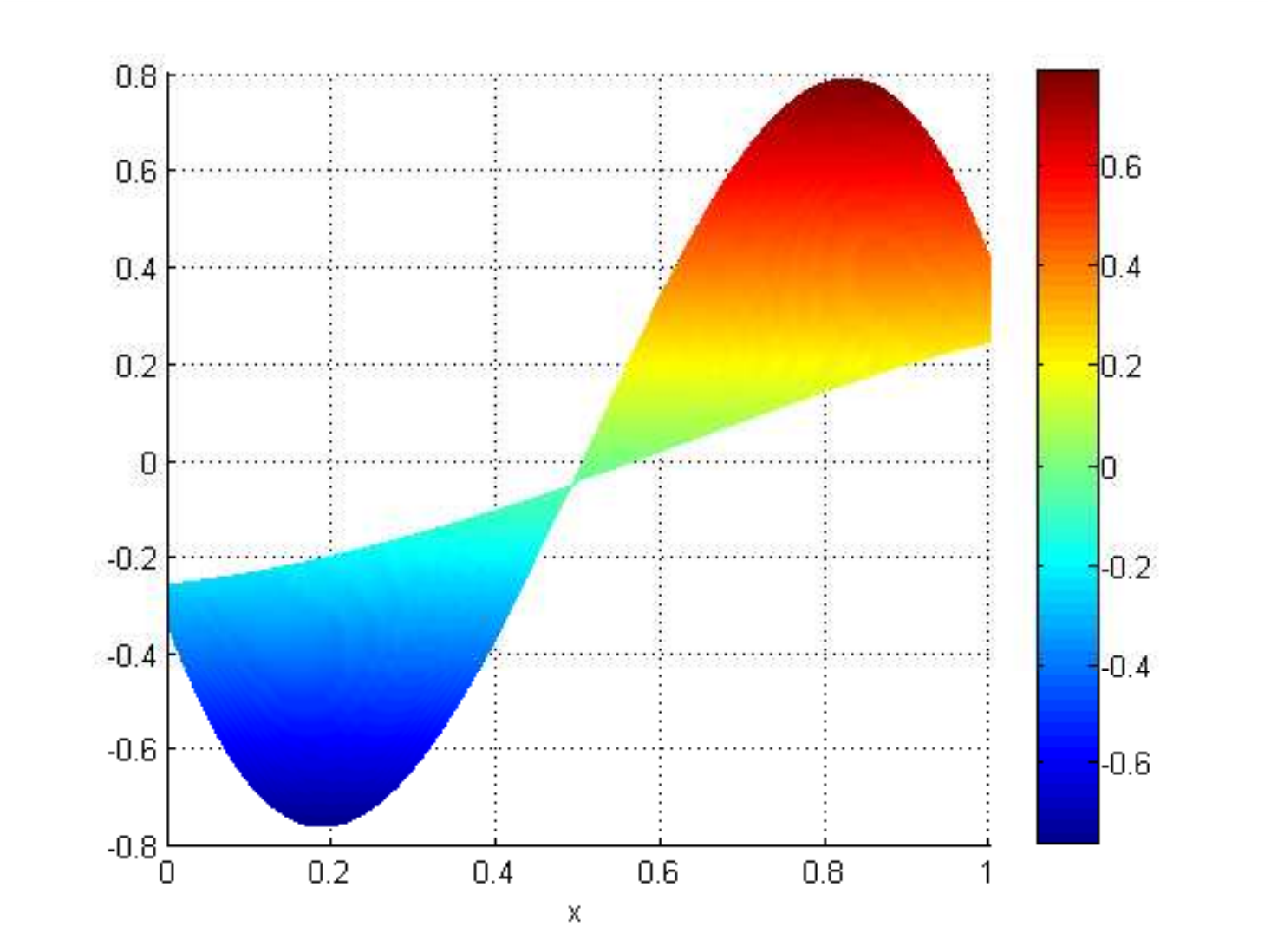}}
\hfill
\subfloat[$\beta=5$]{\includegraphics[width=0.3\textwidth]{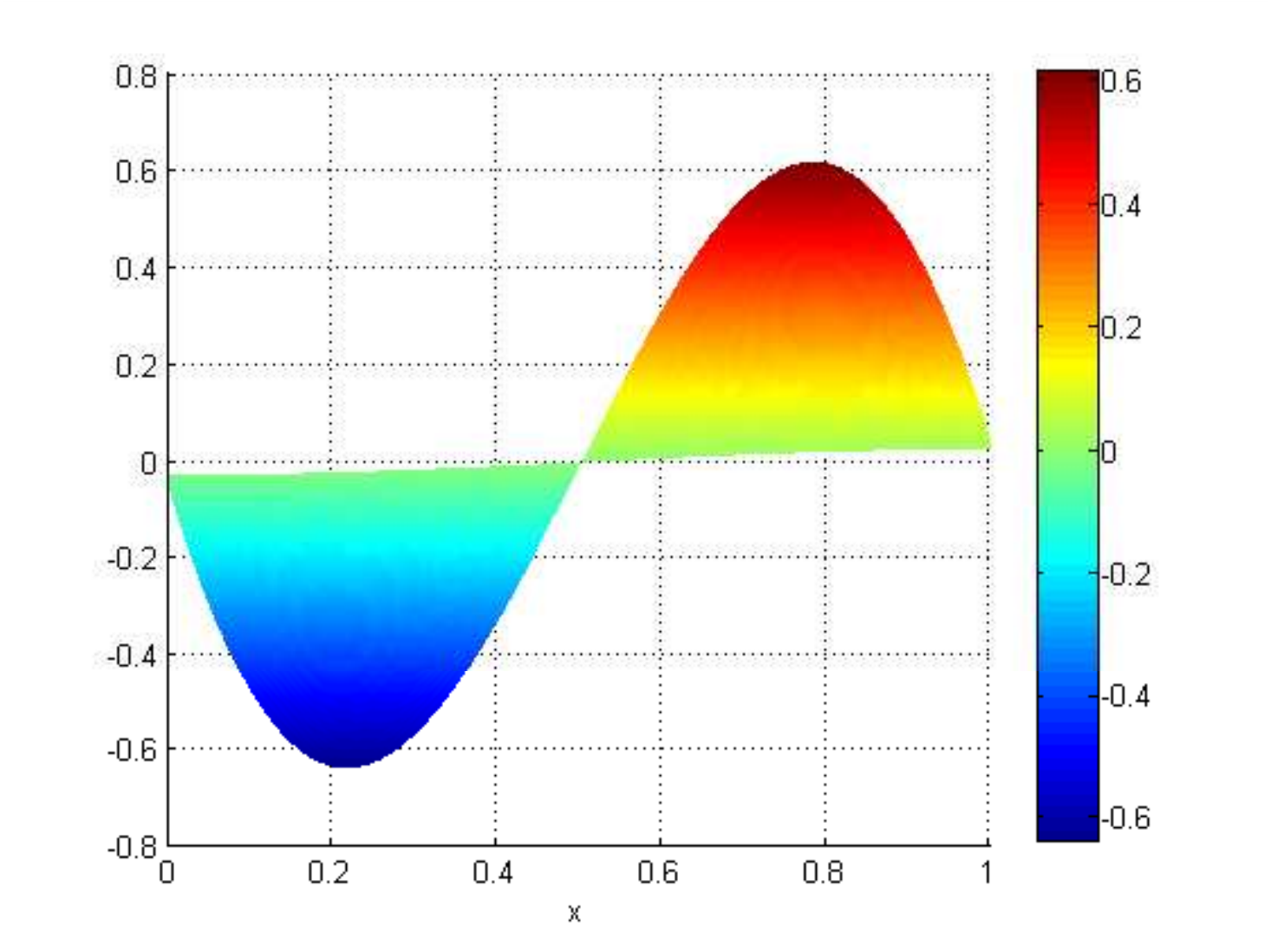}}
\hfill
\subfloat[$\beta=100$]{\includegraphics[width=0.3\textwidth]{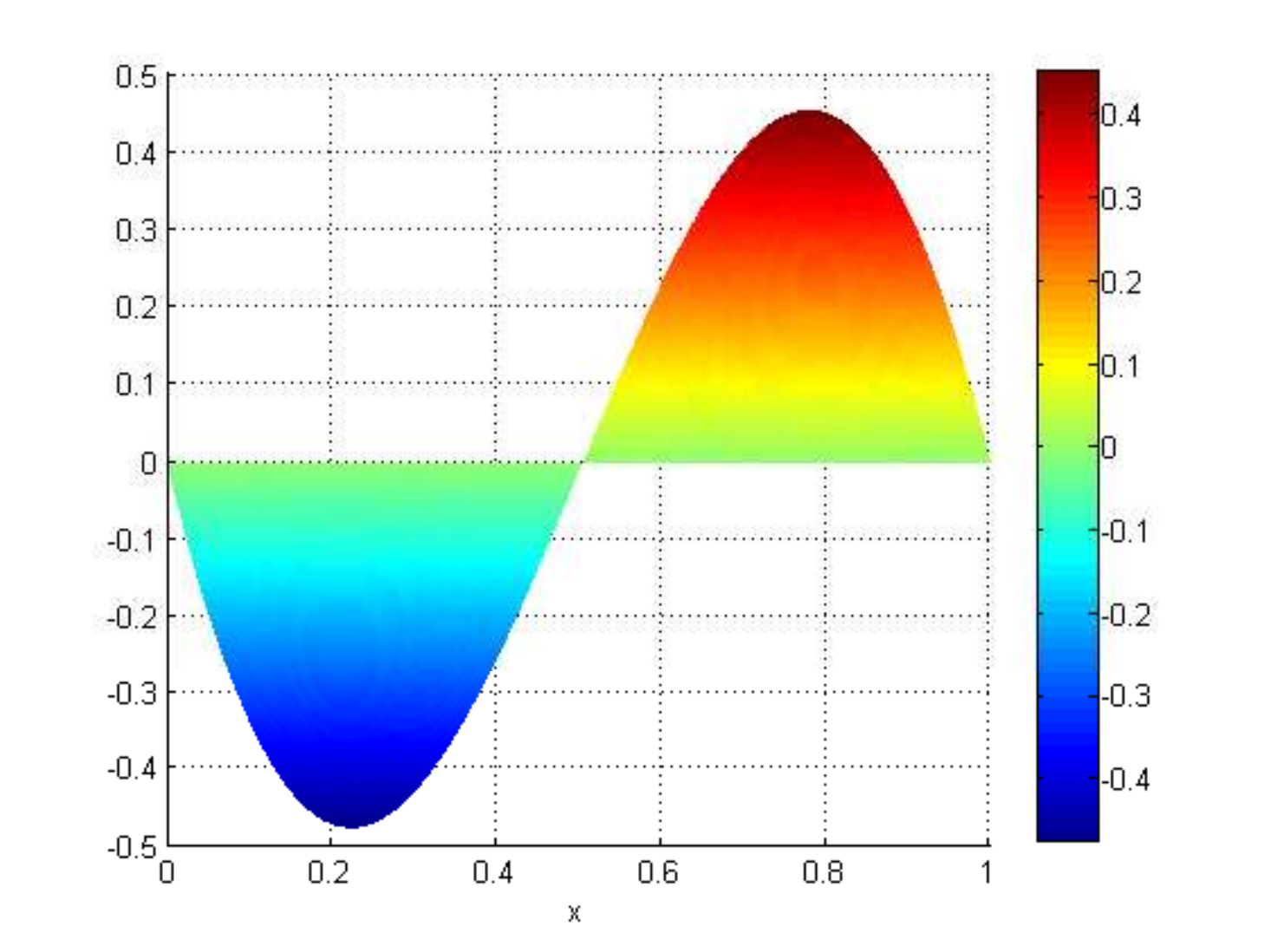}} 
\caption{Profiles of the discrete solution $u_h^k$ to \eqref{CHLW:lim} at iteration $k = 200$ for $\alpha = 1$ and $\beta \in \{0, 1.5, 2 \}$ (top row) and $\beta \in \{3,5,100\}$ (middle row) all viewed in the same camera line of sight.  (Bottom row) Profiles of the discrete solution for $\beta \in \{3,5,100\}$ with a different camera line of sight.  Note that the range of the $z$-axis can be different with each plot.}
\label{fig:1:varybeta}
\end{figure}

\paragraph{Effects of varying $\beta$ with $\alpha = 1$.}  Next, we take the same parameter values for $\tau, \eps, \delta, \kappa$, and the same potentials as in \eqref{Ex1:para}, but now with $\alpha = 1$ and $\beta \in \{0, 1.5, 2, 3, 5, 100\}$.  In Figure \ref{fig:1:varybeta} we plot the discrete solution $u_h^k$ at iteration $k = 200$ for $\alpha = 1$ and $\beta \in  \{0, 1.5, 2, 3, 5, 100\}$.  Similar profiles are observed for the other cases where $\alpha$ and $\beta$ are positive/negative, and so we do not report them here.  We mention that in preliminary testings the profiles for $\beta \in \{-1.5, -2,-2.5\}$ are similar to their positive counterparts but reflected along the line $\{x = 0.5\}$. 

Compared to Figure \ref{fig:1:varyalpha}, it appears that the parameter $\beta$ has a stronger influence on the solution profile than $\alpha$.  Most interestingly, as $\beta$ increases in magnitude, the trace of the discrete solution on the domain boundary $\Gamma_h$ approaches zero, while the range in the bulk domain decreases, as seen in the bottom row of Figure \ref{fig:1:varybeta}.  Whereas for $\beta < 1$, the solution profile exhibits some monotonicity along the $x$-direction.  The transition between the two different profiles occur approximately around $\beta = 1$.

Let us attempt an explanation of the numerical observations.  The similarity of solution profiles for positive/negative values of $\alpha$ and $\beta$ can be attributed to the symmetry of the potential $G$ and the initial condition $u_0$.  In choosing $\delta/\eps = 0.1$ we had hoped that the boundary dynamics might dominate the bulk dynamics, and prioritise driving the trace of $u_h^k$ to the stable minima of $G(\alpha^{-1}(s-\beta))$, which in this case are $\{\beta - \alpha, \beta +\alpha\}$.  For the case $\abs{\beta} \geq 1$ and $\alpha =1$ this is inconsistent with the conservation of boundary mass, and so we believe that the system instead favours to have the trace of $u_h^k$ to be close to zero.  In turn the solution has to accommodate by moving away from the monotonic profile inherited from the initial condition to that observed in Figure \ref{fig:1:varybeta}.  In particular, for our choice of parameter values, potentials and initial condition, it appears that the mass conservation mechanism imposes a stronger influence than the Cahn--Hilliard phase separation mechanism.  

\begin{figure}[h]
\captionsetup[subfigure]{labelformat=empty}
\centering
\subfloat[$\beta=0$]{\includegraphics[width=0.3\textwidth]{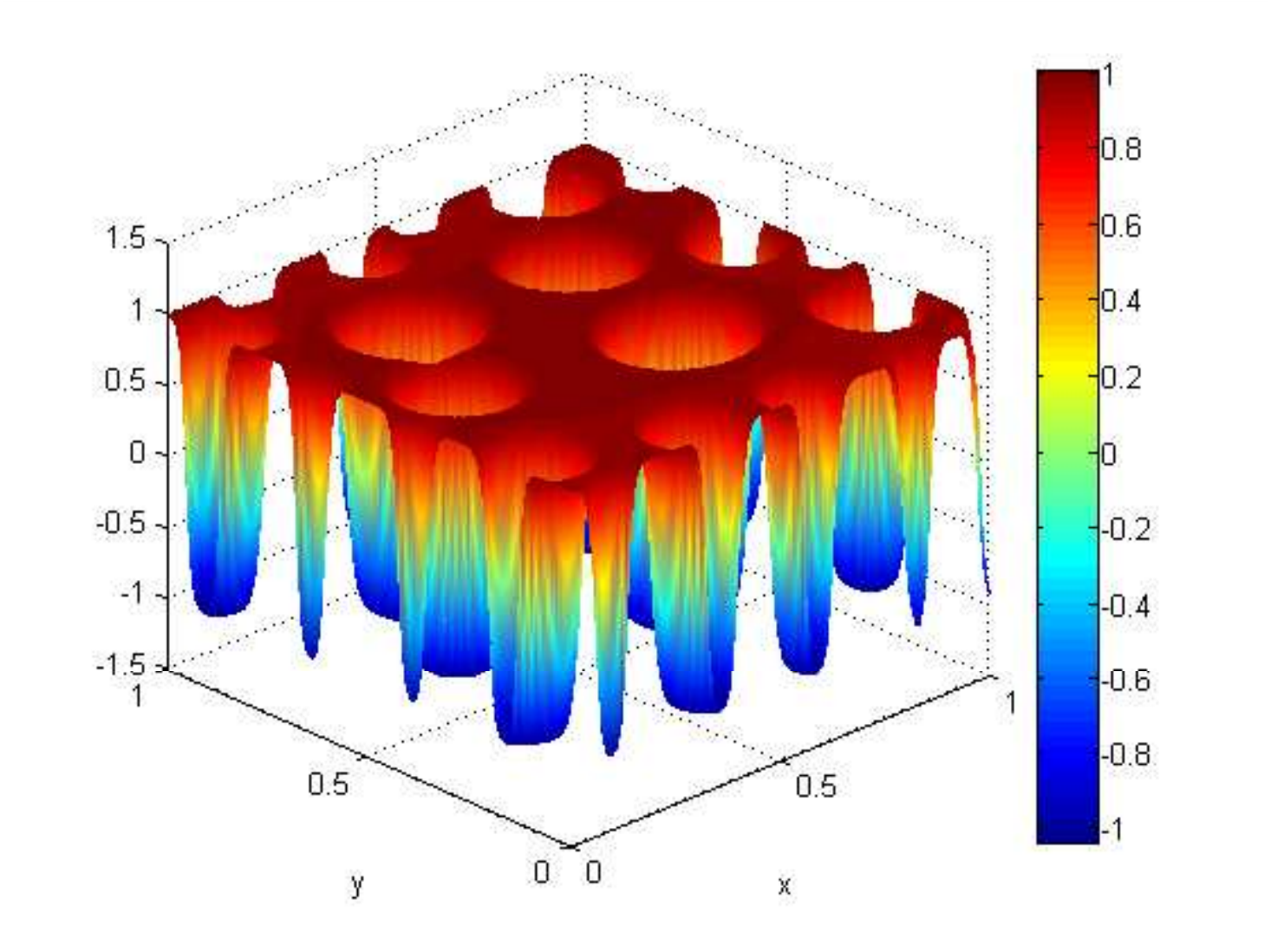}}
\hfill
\subfloat[$\beta=0.3$]{\includegraphics[width=0.3\textwidth]{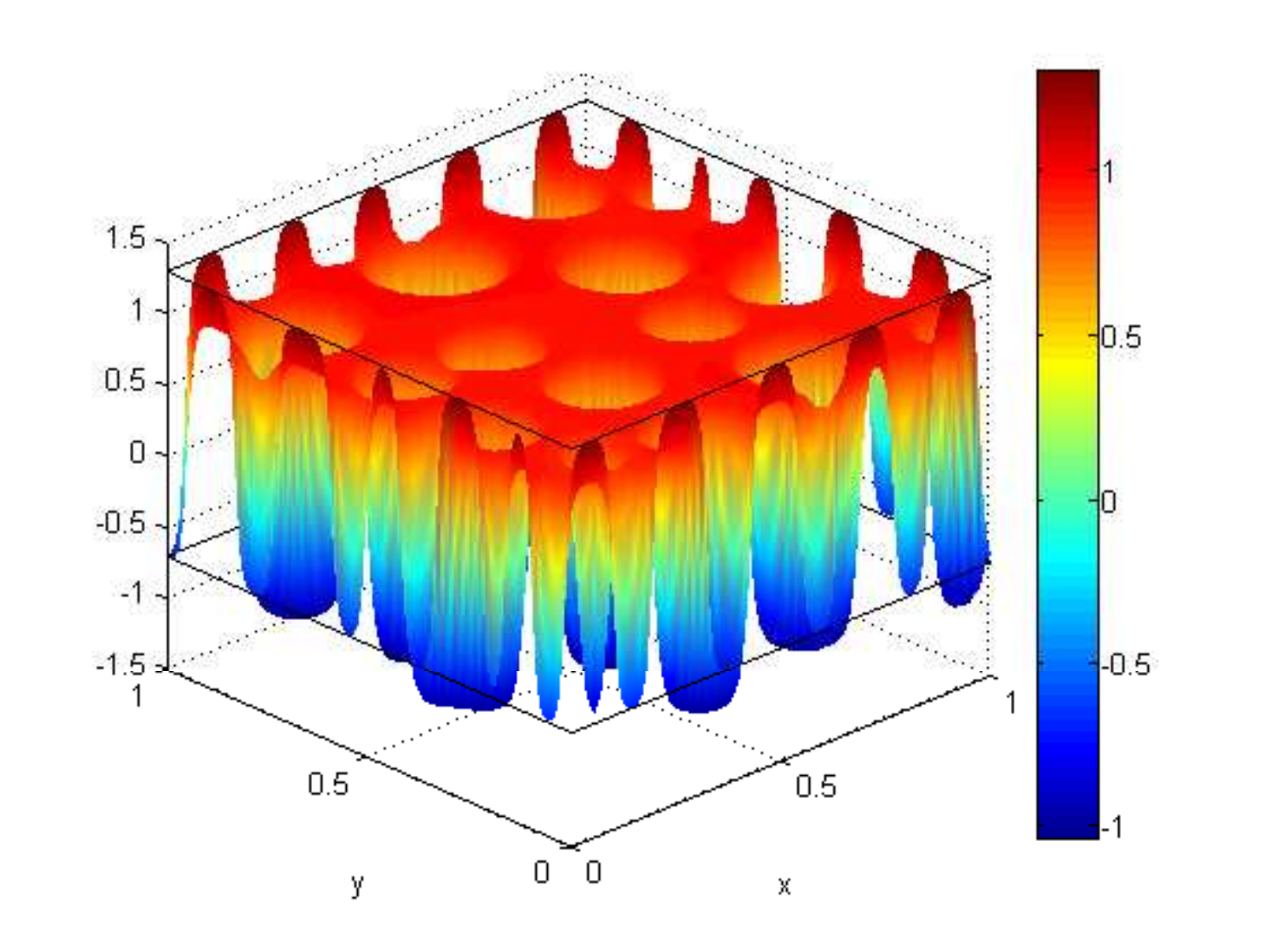}}
\hfill
\subfloat[$\beta=-0.3$]{\includegraphics[width=0.3\textwidth]{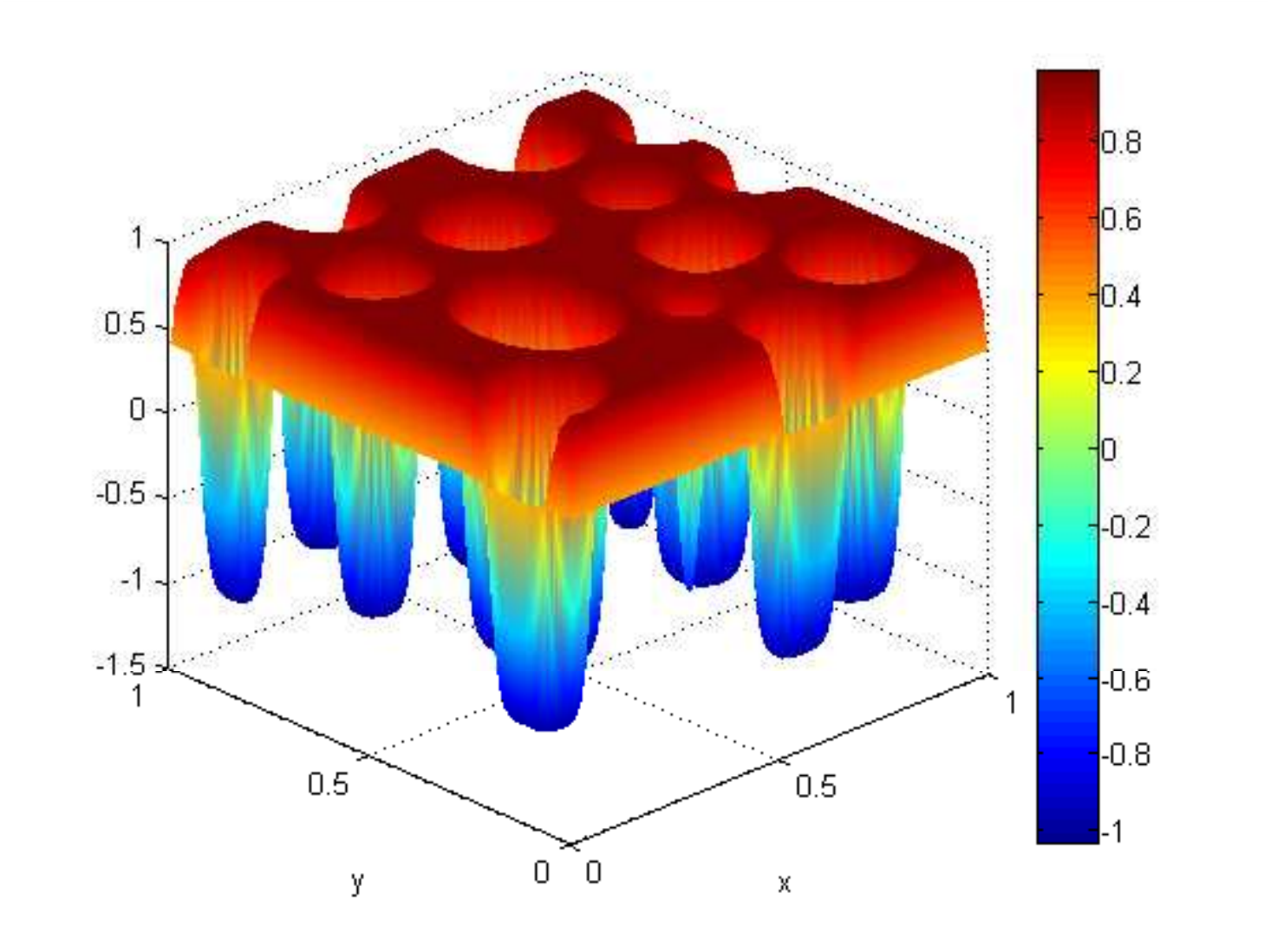}}\\
\subfloat[$\beta=-1$]{\includegraphics[width=0.3\textwidth]{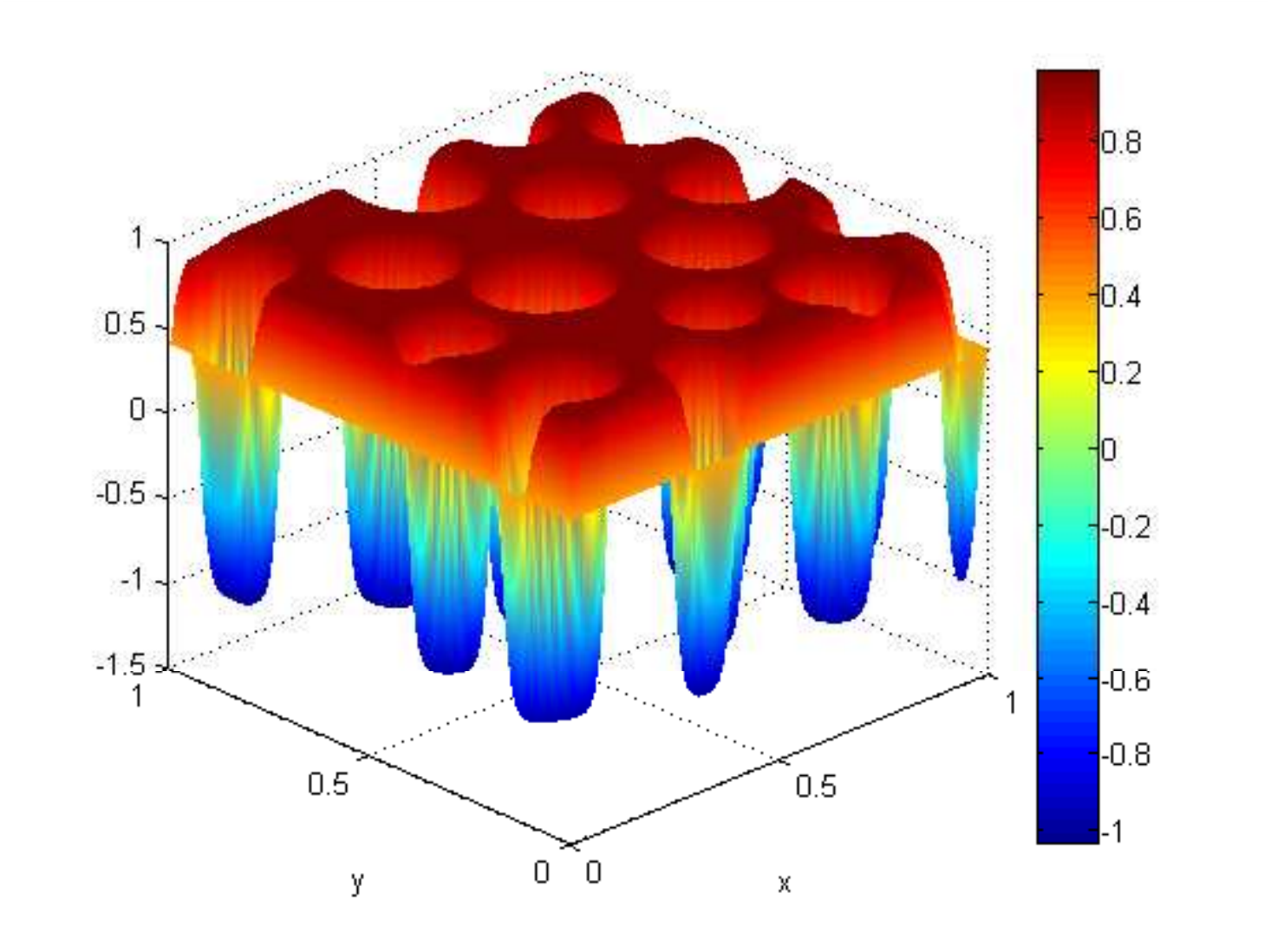}}
\hfill
\subfloat[$\beta=0.9$]{\includegraphics[width=0.3\textwidth]{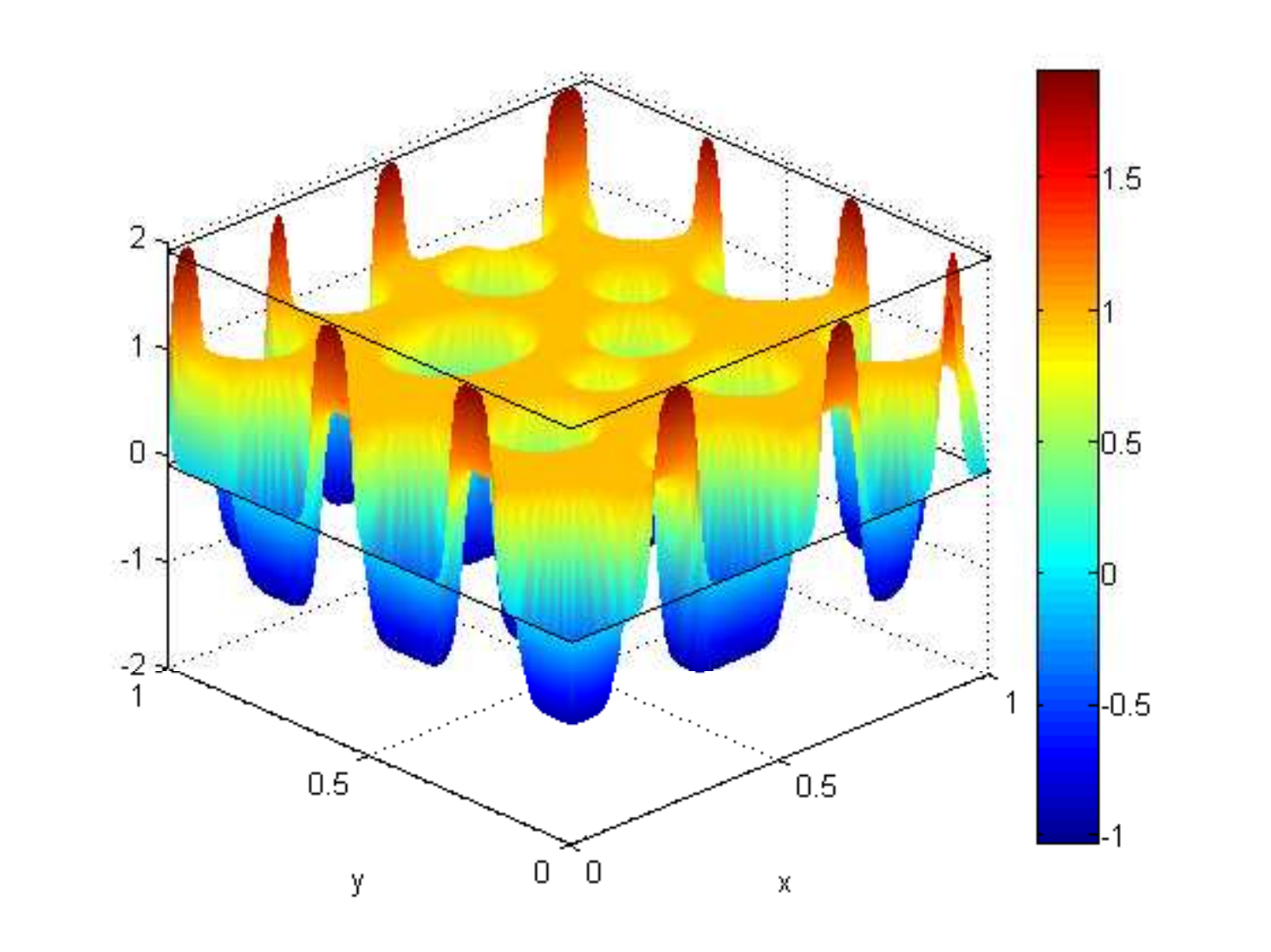}} 
\hfill
\subfloat[$\beta=-0.9$]{\includegraphics[width=0.3\textwidth]{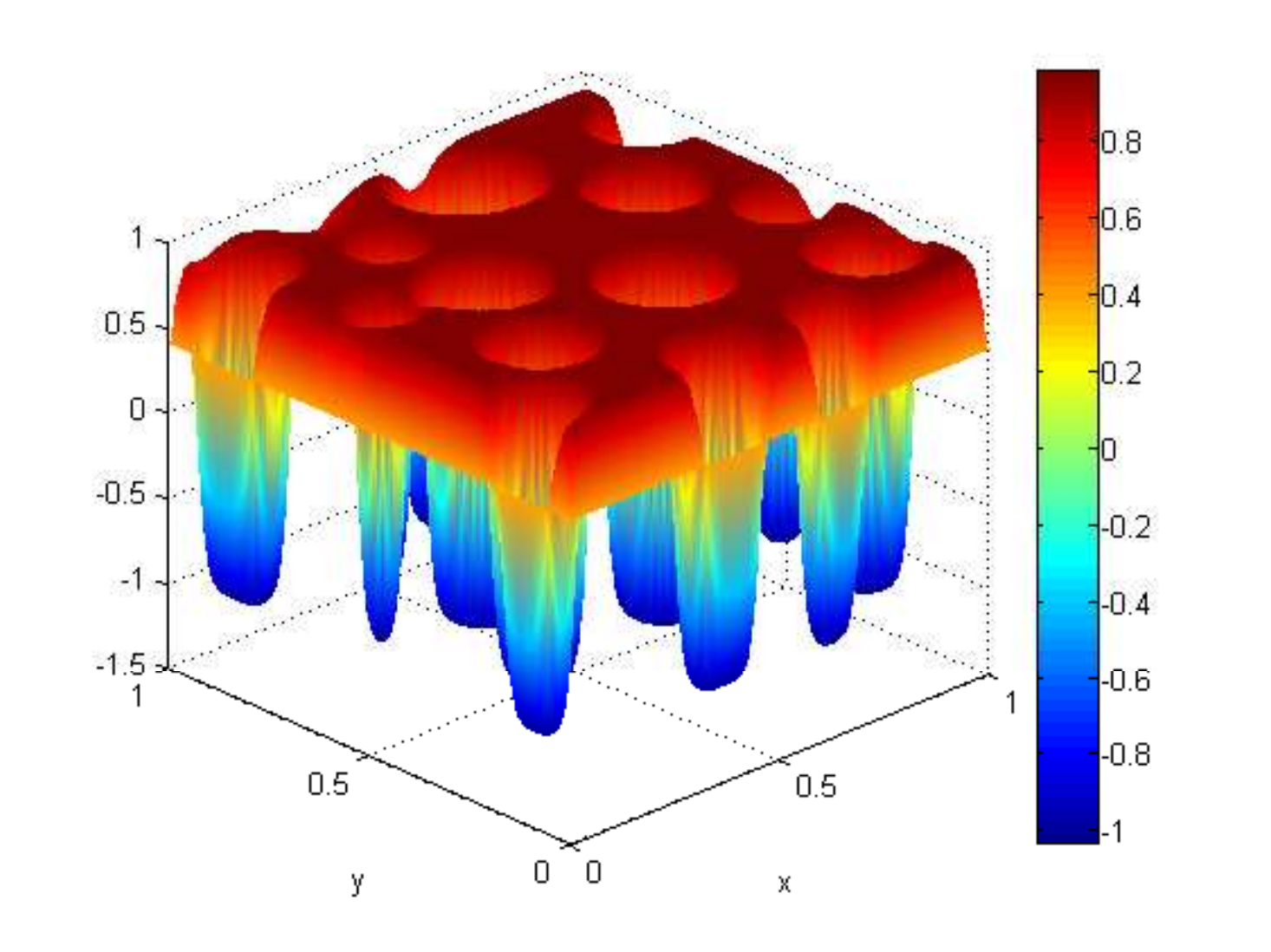}} \\
\subfloat[$\beta=10$]{\includegraphics[width=0.3\textwidth]{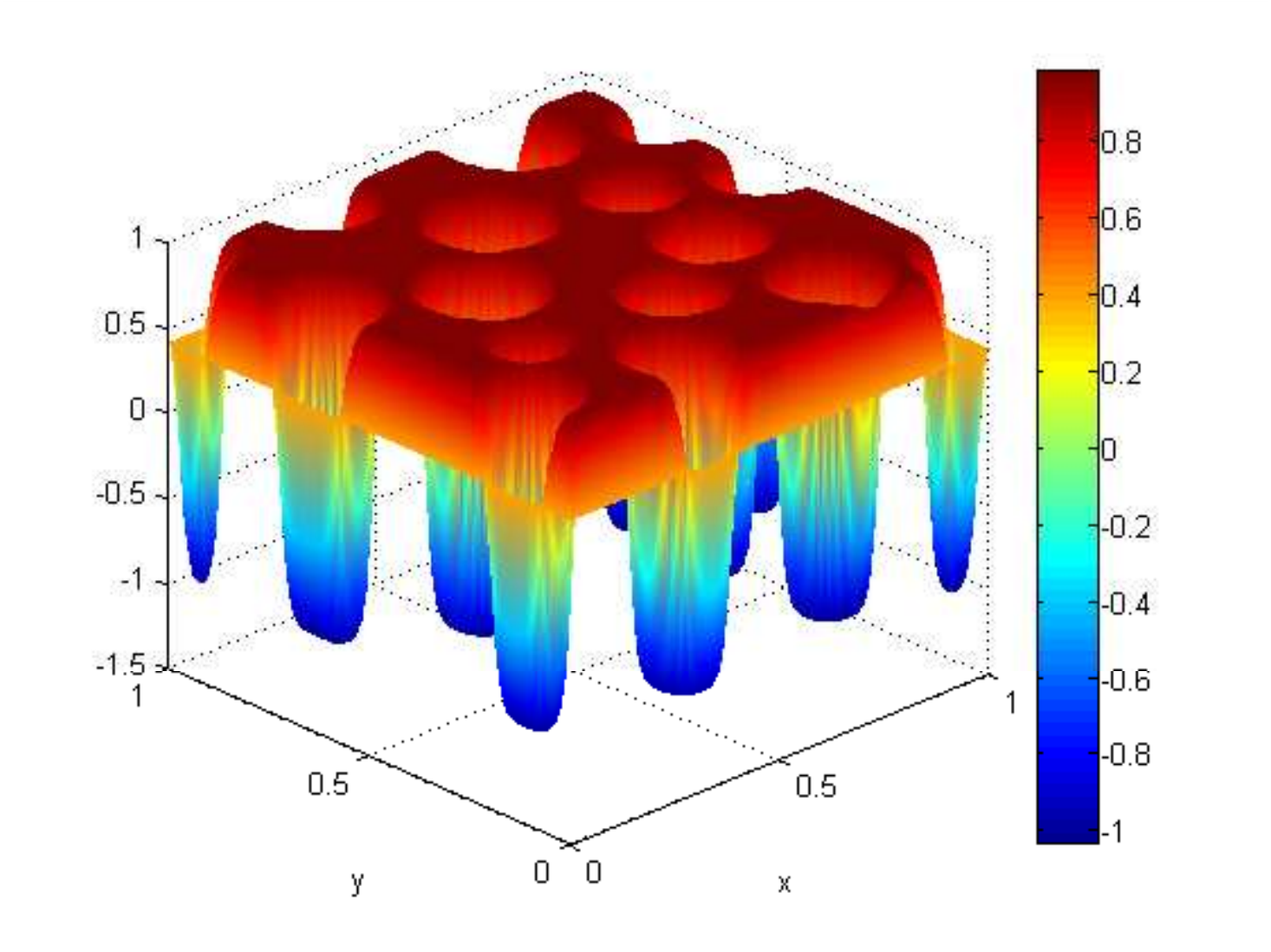}}
\hfill \subfloat[$\beta=2$]{\includegraphics[width=0.3\textwidth]{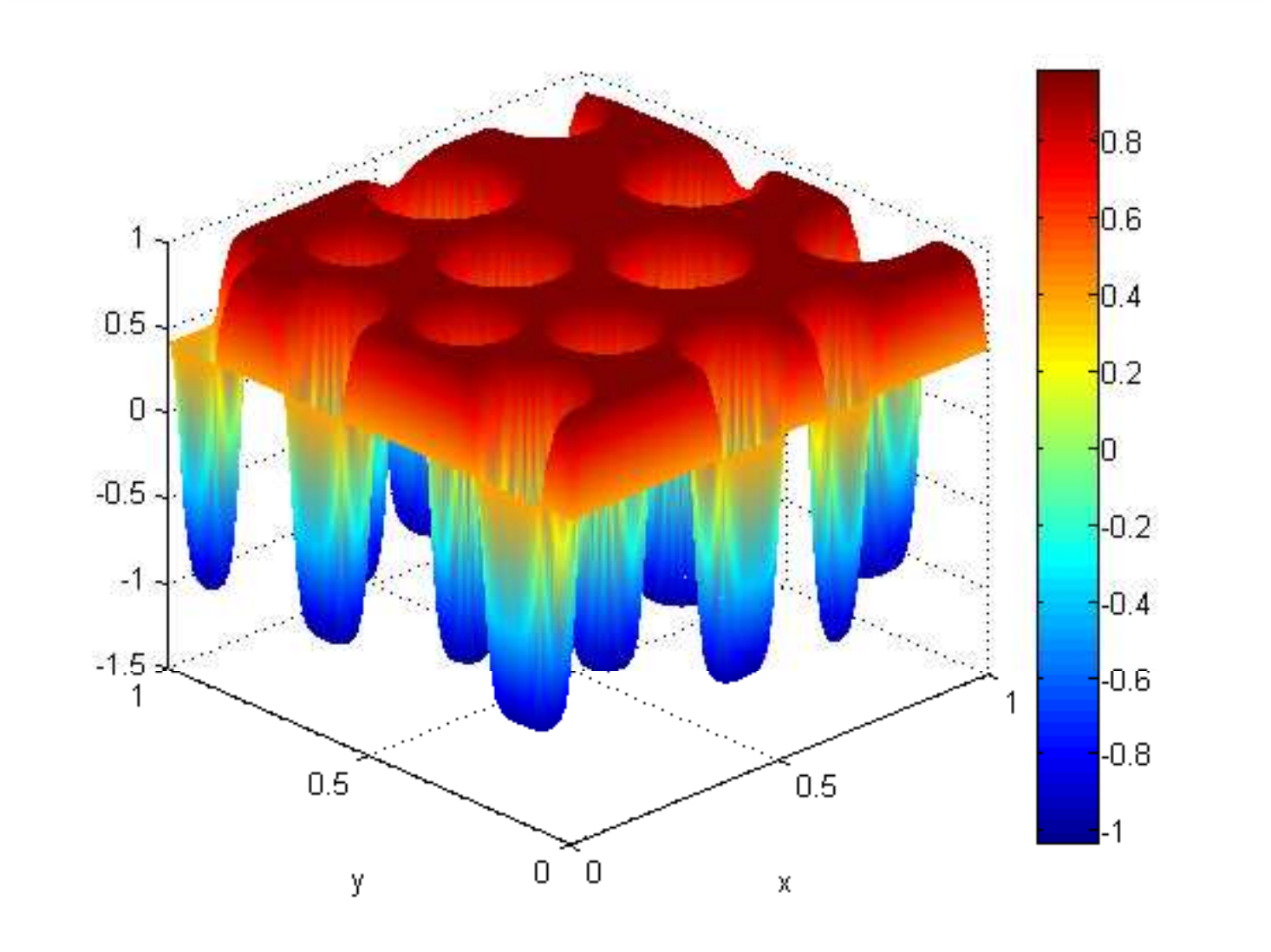}} 
\hfill
\subfloat[$\beta=-2$]{\includegraphics[width=0.3\textwidth]{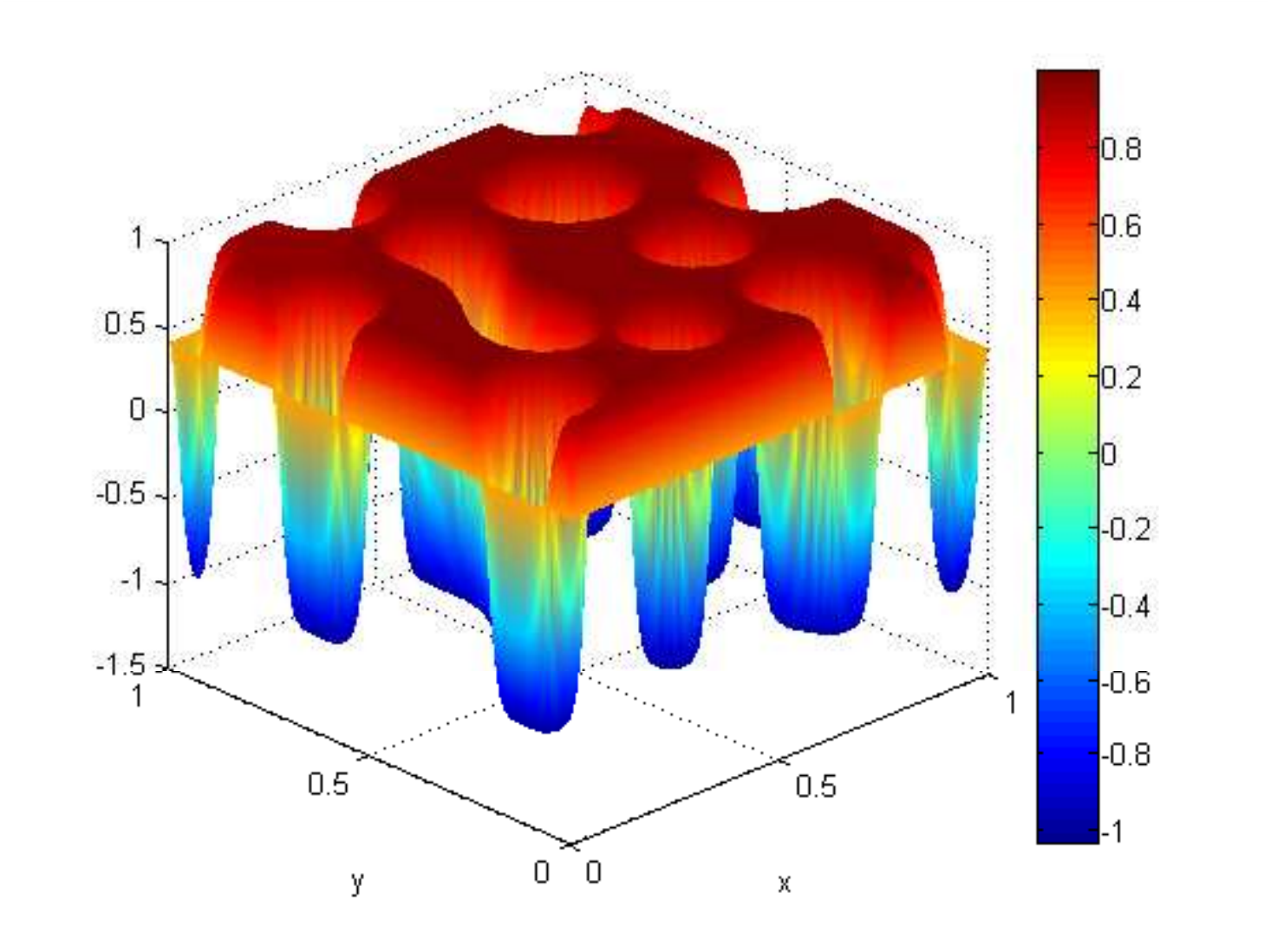}}
\caption{Discrete solution $u_h^k$ to \eqref{CHLW:lim} at iteration $k = 200$ for $\alpha = 1$ and $\beta \in \{0, 0.3, -0.3\}$ (top row) and $\beta \in \{-1, 0.9, -0.9\}$ (middle row) and $\beta \in \{10, 2, -2\}$ (bottom row).  The two black level lines in the plot for $\beta = 0.3$ are $\{z = 1.3\}$ and $\{z = -0.7\}$, and the two black level lines in the plot for $\beta = 0.9$ are $\{z = 1.9\}$ and $\{z = -0.1\}$. }
\label{fig:1:random:bet}
\end{figure}

This competition between the two mechanism is more pronounced by considering an asymmetric initial condition $u_0$, which we take with values uniformly distributed in $[0.3, 0.5]$ such that $\mean{u_0}_\Omega = 0.399$ and $\mean{u_0}_\Gamma = 0.405$ to three decimal points.   We now consider the parameters and potentials to be
\begin{align}\label{num:para}
\tau = 10^{-5}, \quad \eps = \delta = 0.01, \quad \kappa = 1, \quad F(s) = G(s) = \tfrac{1}{4}(s^2-1)^2.
\end{align} 
The values of $\eps$ and $\delta$ are chosen to allow for  the typical Cahn--Hilliard phase separation dynamics, which is observed in Figure \ref{fig:1:random:bet}, where we plot the discrete solution $u_h^k$ at iteration $k = 200$ for $\alpha = 1$ and selected values of $\beta$.

We point out two observations from Figure \ref{fig:1:random:bet}.  The first observation is, unlike Figure \ref{fig:1:varybeta}, there is no symmetry corresponding to the case of positive and negative values of $\beta$, as evident in the second and third columns.  The second observation is that for $\abs{\beta} \geq 1$, the trace of the discrete solution is nearly constant on the domain boundary, with an average value of 0.401 to three decimal points.  In particular, this is in good agreement with the conservation of boundary mass ($\abs{\Gamma}\mean{u_0}_\Gamma = 1.606$ to three decimal points), and is similar to the behaviour in the previous set of experiments.  On the other hand, for $\abs{\beta}$ small there is a remarkable difference in solution behaviour.  In the second column of Figure \ref{fig:1:random:bet} we have the plots for $\beta = 0.3$ and $\beta = 0.9$, in which the Cahn--Hilliard phase separation dynamics can be clearly seen on the domain boundary.  The black level lines in each of these two plots correspond to the stable minima of the surface potential $G(s-\beta)$, which are $\{-0.7, 1.3\}$ for $\beta = 0.3$ and $\{-0.1, 1.9\}$ for $\beta = 0.9$.  However, for $\beta = -0.3$ and $-0.9$, we see that the trace of the discrete solution takes nearly constant values, similarly as in the case for large $\abs{\beta}$ values.  We further comment that the solution behaviour is insensitive to the sign of $\alpha$, since the stable minima are $\{\beta - \alpha, \beta + \alpha \}$ regardless of the sign of $\alpha$.

Based on Figure \ref{fig:1:random:bet}, if $s_1$ and $s_2$ denote the stable minima $\{\beta - \alpha, \beta + \alpha\}$ of $G(\alpha^{-1}(s-\beta))$ with $s_1 < s_2$, we could postulate that Cahn--Hilliard phase separation is more likely to occur on the domain boundary if $\mean{u_0}_\Gamma \in (s_1,s_2)$.  Otherwise, it is likely that the trace of the solution is spatially constant with values close to $\mean{u_0}_\Gamma$.

However, in some simulations not reported here we have observed counterexamples to this conjecture.  In particular, for $(\alpha, \beta) \in \{(0.8,0.9), (0.9,0.9), (1.5, -0.9), (2,-0.9)\}$, the trace of the discrete solution is nearly constant with values close to $\mean{u_0}_\Gamma$, even though $\mean{u_0}_\Gamma \in (s_1, s_2)$ in each of these cases.  Therefore, at present we do not have a robust method to predict whether Cahn--Hilliard phase separation or near constant trace values will occur based on knowledge of the parameters and the initial data.

\subsubsection{Dynamics of the Robin system \eqref{CHLW}}

\paragraph{Convergence as $K \to 0$.} We now demonstrate numerical evidence in support of Theorem \ref{thm:conv}.  The potentials are as described in \eqref{num:para}, and for the following choice of parameters:
\begin{align*}
\tau = 10^{-5}, \quad \eps = \delta = 0.02.
\end{align*}
In this section we report on the errors between the discrete solutions to the limit model \eqref{CHLW:lim} and the Robin model \eqref{CHLW} as $K$ decreases.  For fixed $K_1 > K_2$, let $e(K_i)$ denote the error under consideration at $K_i$, $i = 1,2$.  Then, the experimental order of convergence is defined as
\begin{align*}
\EOC(K_1,K_2) = \frac{\log(e(K_1)/e(K_2))}{\log(K_1/K_2)}.
\end{align*}  
For $p \in \{2,4\}$ and $X \in \{\Omega, \Gamma\}$, the $L^p(0,T;L^2(X))$-norm of a function $f$ is approximated as
\begin{align*}
\norm{f}_{L^p(0,T;L^2(X))} \approx \tau^{\frac{1}{p}} \left  [\sum_{i=1}^N \big ( \int_X I^h \abs{f_h^i}^2 \big)^{\frac{p}{2}} \right ]^{\frac{1}{p}},
\end{align*}
where $f_h^i$ is the discrete approximation of $f$ at iteration $i$.  

\begin{figure}[h!!]
	\centering
	\captionsetup[subfigure]{labelformat=empty}
	%\textbf{Plots of the initial data:}\\[2ex]
	\includegraphics[width=0.4\textwidth]{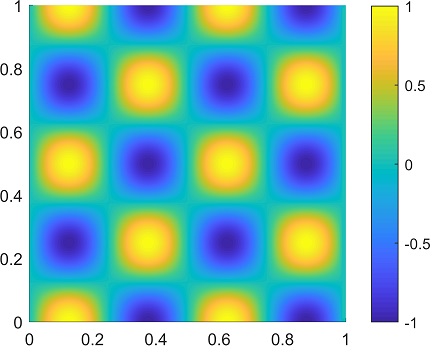}
	\hspace{20pt}
	\includegraphics[width=0.4\textwidth]{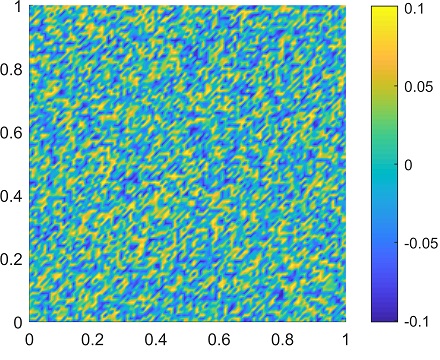}
	\vspace{-5pt}
	\caption{Discrete initial data $u_h^0$ in the cases (a) (left) and (b) (right).}\label{FIG:INI}
\end{figure}

\begin{figure}[h!!]
	\centering
	\captionsetup[subfigure]{position=top,labelformat=empty}
	%\textbf{Simulations for initial datum (a):}\\[-1ex]
	%%%%%%%%%%%%%%%%%%%%%%%%%%%%%%%%%%%%
	%%% Simulation Wave 1: a=1, b=0, k=1
	%%%%%%%%%%%%%%%%%%%%%%%%%%%%%%%%%%%%
	\subfloat[$K=10$]{\includegraphics[width=0.24\textwidth]{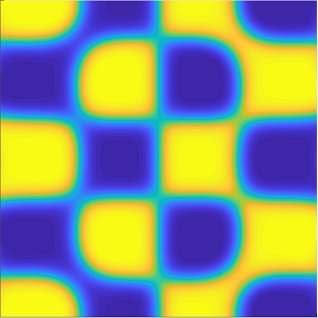}}
	\hfill
	\subfloat[$K=1$]{\includegraphics[width=0.24\textwidth]{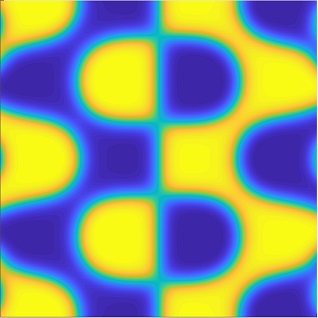}}
	\hfill
	\subfloat[$K=0.1$]{\includegraphics[width=0.24\textwidth]{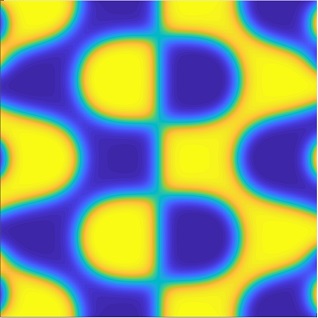}}
	\hfill
	\subfloat[System \eqref{lim:alt}]{\includegraphics[width=0.24\textwidth]{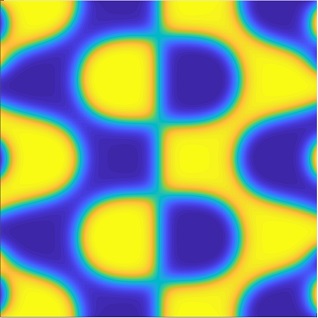}}\\[0.5ex]
%\end{figure}
%\begin{figure}[h!]
%	\centering
%	\vspace{-15pt}
	%%%%%%%%%%%%%%%%%%%%%%%%%%%%%%%%%%%%
	%%% Simulation Wave 2: a=1, b=0, k=1
	%%%%%%%%%%%%%%%%%%%%%%%%%%%%%%%%%%%%
	\includegraphics[width=0.24\textwidth]{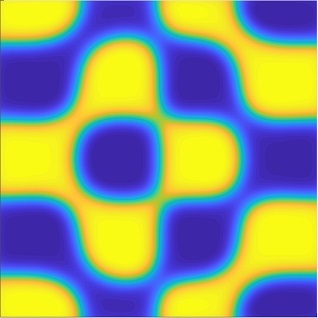}
	\hfill
	\includegraphics[width=0.24\textwidth]{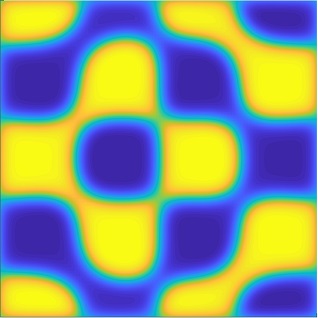}
	\hfill
	\includegraphics[width=0.24\textwidth]{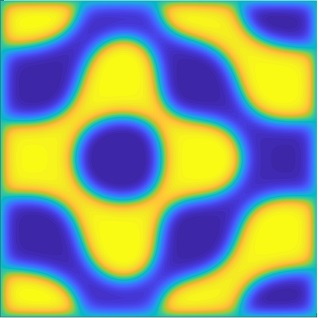}
	\hfill
	\includegraphics[width=0.24\textwidth]{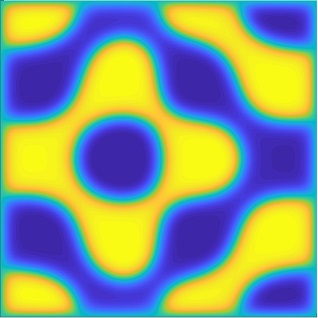}\\[0.5ex]
%\end{figure}
%\begin{figure}[h!]
%	\centering
%	\vspace{-15pt}
	%%%%%%%%%%%%%%%%%%%%%%%%%%%%%%%%%%%%
	%%% Simulation Wave 3: a=1, b=0, k=1
	%%%%%%%%%%%%%%%%%%%%%%%%%%%%%%%%%%%%
	\includegraphics[width=0.24\textwidth]{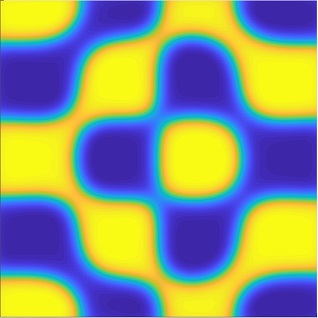}
	\hfill
	\includegraphics[width=0.24\textwidth]{SW3_wave_10.jpg}
	\hfill
	\includegraphics[width=0.24\textwidth]{SW3_wave_10.jpg}
	\hfill
	\includegraphics[width=0.24\textwidth]{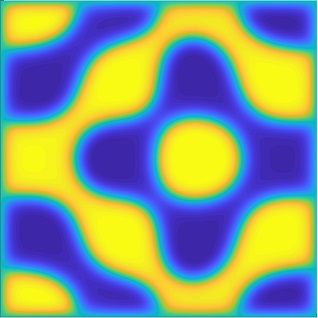}\\
	\vspace{8pt}
	\includegraphics[width=0.42\textwidth]{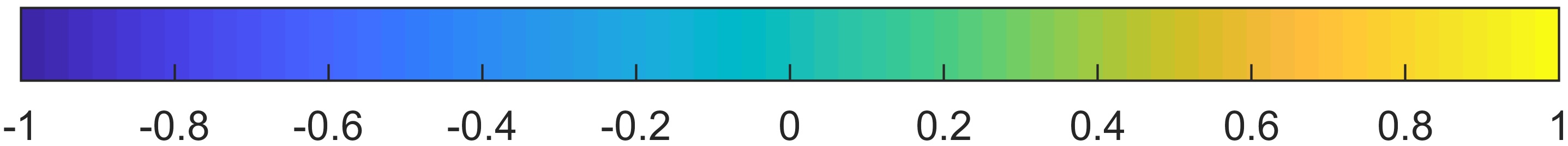}
	\caption{Numerical solution $u_h^k$ for initial data (a) at iteration $k=100$ for $\kappa=1$ as well as $(\alpha,\beta)=(1,0)$ (first row), $(\alpha,\beta)=(2,-4)$ (second row) and $(\alpha,\beta)=(-2,4)$ (third row).}\label{FIG:SIM:A}
\end{figure}

\begin{figure}[h!!]
	\centering
	\captionsetup[subfigure]{position=top,labelformat=empty}
	%\textbf{Simulations for initial datum (b):}\\[-1ex]
	%%%%%%%%%%%%%%%%%%%%%%%%%%%%%%%%%%%%
	%%% Simulation Random 1: a=1, b=0, k=1
	%%%%%%%%%%%%%%%%%%%%%%%%%%%%%%%%%%%%
	\subfloat[$K=10$]{\includegraphics[width=0.24\textwidth]{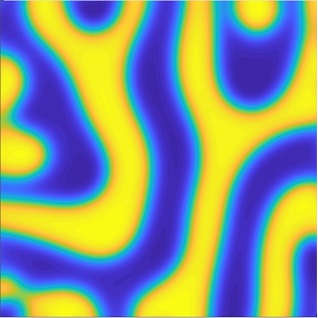}}
	\hfill
	\subfloat[$K=1$]{\includegraphics[width=0.24\textwidth]{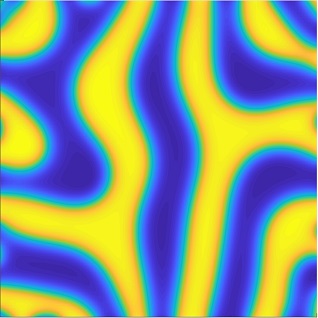}}
	\hfill
	\subfloat[$K=0.1$]{\includegraphics[width=0.24\textwidth]{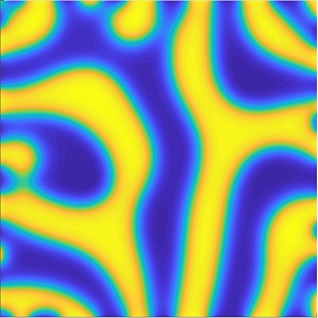}}
	\hfill
	\subfloat[System \eqref{lim:alt}]{\includegraphics[width=0.24\textwidth]{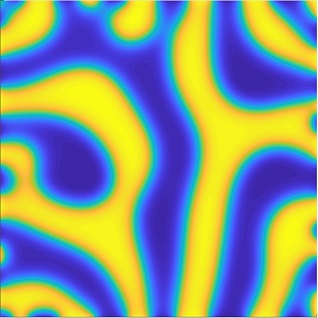}}\\[0.5ex]
%\end{figure}
%\begin{figure}[h!]
%	\centering
%	\vspace{-15pt}
	%%%%%%%%%%%%%%%%%%%%%%%%%%%%%%%%%%%%
	%%% Simulation Random 2: a=1, b=0, k=1
	%%%%%%%%%%%%%%%%%%%%%%%%%%%%%%%%%%%%
	\includegraphics[width=0.24\textwidth]{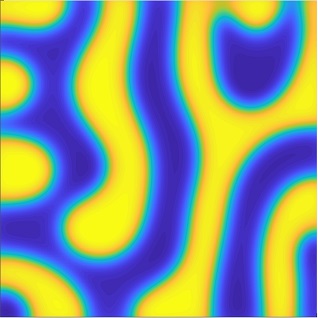}
	\hfill
	\includegraphics[width=0.24\textwidth]{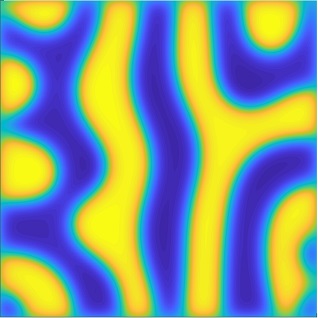}
	\hfill
	\includegraphics[width=0.24\textwidth]{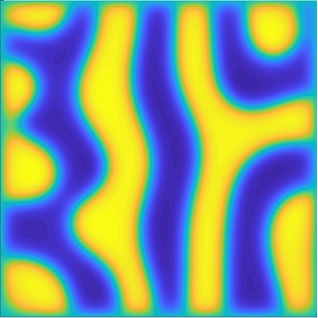}
	\hfill
	\includegraphics[width=0.24\textwidth]{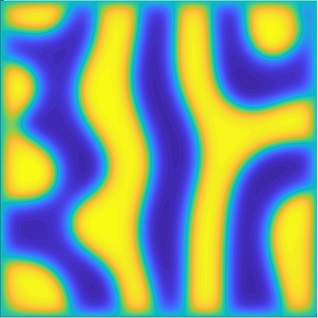}\\[0.5ex]
%\end{figure}
%\begin{figure}[h!]
%	\centering
%	\vspace{-15pt}
	%%%%%%%%%%%%%%%%%%%%%%%%%%%%%%%%%%%%
	%%% Simulation Random 3: a=1, b=0, k=1
	%%%%%%%%%%%%%%%%%%%%%%%%%%%%%%%%%%%%
	\includegraphics[width=0.24\textwidth]{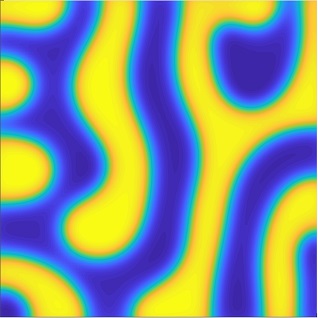}
	\hfill
	\includegraphics[width=0.24\textwidth]{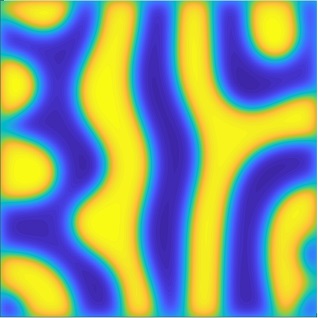}
	\hfill
	\includegraphics[width=0.24\textwidth]{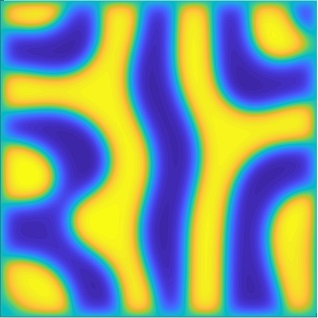}
	\hfill
	\includegraphics[width=0.24\textwidth]{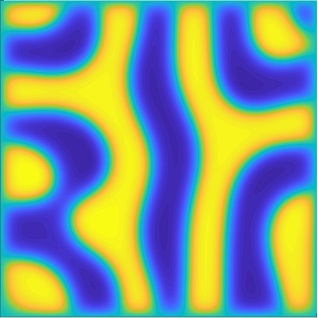}\\
	\vspace{8pt}
	\includegraphics[width=0.42\textwidth]{cb_hor}
	\caption{Numerical solution $u_h^k$ with initial data (b) at iteration $k=100$ for $\kappa=0.4$ as well as $(\alpha,\beta)=(1,0)$ (first row), $(\alpha,\beta)=(2,-4)$ (second row) and $(\alpha,\beta)=(-5,30)$ (third row).}\label{FIG:SIM:B}
\end{figure}

To demonstrate the convergence behaviour for $K\to 0$, we consider two different discrete initial data as illustrated in Figure \ref{FIG:INI}:
\begin{enumerate}[label = (\alph*), ref = (\alph*)]
	\item For any node $x\in\Omega_h$, we set $u_h^0(x) = \sin(4\pi x_1)\cos(4\pi x_2)$.
	\item For any node $x\in\Omega_h$, $u_h^0(x)$ attains a random value in $[-0.1,0.1]$. These random values are uniformly distributed.
\end{enumerate}
In view of Theorem \ref{thm:conv}, we set $v_h^0 = H^{-1}(u_h^0)=\alpha^{-1}(u_0-\beta)$. In particular, the discrete initial data do not depend on $K$.  In Figures \ref{FIG:SIM:A} and \ref{FIG:SIM:B}, we display the discrete solution $u_h^k$ to initial data (a) and (b), respectively, at iteration $k = 100$ for $K \in \{0.1, 1, 10\}$ and three choices of $(\alpha, \beta)$ along with the corresponding discrete solution to the limit system \eqref{lim:alt}.  Visually, we observe good agreement of the discrete solution $u_h^k$ at $K = 0.1$ and the discrete solution to the limit system.  

%\pagebreak[4]

%\paragraph{Simulations for different choices of $\alpha$ and $\beta$.}
%For the following simulations we set $\eps=\delta=0.02$, $\tau=10^{-5}$ and $T=0.001$. 

%\newpage

%\noindent Tables of convergence:\\[1ex]
%\noindent Simulation for initial datum (a) with $(\alpha,\beta)=(1,0)$:\\[1ex]

In Tables \ref{TBL:A:1} to \ref{TBL:B:3}, we collect the error between the discrete solution $(u^K, v^K)$ to \eqref{CHLW} and the discrete solution $(u,v)$  to \eqref{lim:alt} (where $v$ stands for $\alpha^{-1}(u-\beta)$) measured in various norms for several choices of $(\alpha, \beta)$ and for the initial data (a) and (b). In the tables, the labelling of the columns is to be understood as follows: 
\begin{itemize}
	\item The columns $L^2(H^1(\Omega))$ and $L^4(L^2(\Omega))$ display the error $u^K-u$,
	\item the column $L^2(\Sigma_T)$ displays the difference $u^K-(\alpha v^K + \beta)$,
	\item the columns $L^2(H^1(\Gamma))$ and $L^4(L^2(\Gamma))$ display the error $v^K-v$
\end{itemize}
measured in the respective norms on these spaces.
 
In all instances, a linear rate of convergence is observed. In view of Theorem \ref{thm:conv}, this seems to suggest that the error estimate for $u^K - (\alpha v^K + \beta)$ in $L^2(\Sigma_T)$ is sharp, while the rate of convergence for the other terms in \eqref{intro:err} can be improved further.

\begin{table}[h!]
\scriptsize
\begin{tabular}{|C{20pt}||C{40pt}|C{40pt}||C{40pt}|C{40pt}||C{40pt}|C{40pt}|}
	\hline
	\rule{0pt}{9pt}
	$K$ & $L^2(H^1(\Omega))$  & $\EOC$ & $L^4(L^2(\Omega))$ & $\EOC$ & $L^2(\Sigma_T)$ & $\EOC$ \\
	\hline
	\rule{0pt}{9pt}
	$10$       & $2.63e-01$    & -         & $3.98e-02$    & -         & $4.71e-02$    & -       \\ 
	$1$        & $1.20e-01$    & $3.41e-01$    & $1.90e-02$    & $3.22e-01$    & $2.00e-02$    & $3.72e-01$  \\ 
	$10^{-1}$  & $1.47e-02$    & $9.12e-01$    & $2.32e-03$    & $9.12e-01$    & $2.17e-03$    & $9.64e-01$  \\ 
	$10^{-2}$  & $1.53e-03$    & $9.85e-01$    & $2.38e-04$    & $9.90e-01$    & $2.20e-04$    & $9.94e-01$  \\ 
	$10^{-3}$  & $1.53e-04$    & $9.98e-01$    & $2.38e-05$    & $9.99e-01$    & $2.21e-05$    & $9.99e-01$  \\ 
	$10^{-4}$  & $1.53e-05$    & $1.00e+00$    & $2.38e-06$    & $1.00e+00$    & $2.21e-06$    & $1.00e+00$  \\ 
	$10^{-5}$  & $1.53e-06$    & $1.00e+00$    & $2.38e-07$    & $1.00e+00$    & $2.21e-07$    & $1.00e+00$  \\ 
	\hline
\end{tabular}
%\noindent
\begin{tabular}{|C{20pt}||C{40pt}|C{40pt}||C{40pt}|C{40pt}|}
	\hline 
	\rule{0pt}{9pt}
	$K$ & $L^2(H^1(\Gamma))$ & $\EOC$ & $L^4(L^2(\Gamma))$  & $\EOC$ \\
	\hline 
	$10$       & $3.32e-01$    & -         & $7.34e-02$    & -         \\ 
	$1$        & $1.24e-01$    & $4.28e-01$    & $3.16e-02$    & $3.67e-01$    \\ 
	$10^{-1}$  & $1.66e-02$    & $8.73e-01$    & $4.45e-03$    & $8.51e-01$    \\ 
	$10^{-2}$  & $1.77e-03$    & $9.71e-01$    & $4.76e-04$    & $9.71e-01$    \\ 
	$10^{-3}$  & $1.78e-04$    & $9.97e-01$    & $4.79e-05$    & $9.97e-01$    \\ 
	$10^{-4}$  & $1.78e-05$    & $1.00e+00$    & $4.79e-06$    & $1.00e+00$    \\ 
	$10^{-5}$  & $1.78e-06$    & $1.00e+00$    & $4.79e-07$    & $1.00e+00$    \\ 
	\hline
\end{tabular}
\caption{Discrete error measured in various norms between solutions to \eqref{lim:alt} and \eqref{CHLW} for $(\alpha,\beta) = (1,0)$ with initial data (a).}
\label{TBL:A:1}
\end{table}

%\normalsize\bigskip

%\pagebreak[4]

%\noindent Simulation for initial datum (a) with $(\alpha,\beta)=(2,-4)$ (or $(\alpha,\beta)=(-2,4)$, respectively):\\[1ex]
\begin{table}[h!]
\scriptsize
\begin{tabular}{|C{20pt}||C{40pt}|C{40pt}||C{40pt}|C{40pt}||C{40pt}|C{40pt}|}
	\hline
	\rule{0pt}{9pt}
	$K$ & $L^2(H^1(\Omega))$  & $\EOC$ & $L^4(L^2(\Omega))$ & $\EOC$ & $L^2(\Sigma_T)$ & $\EOC$ \\
	\hline
	\rule{0pt}{9pt}
	$10$       & $2.66e-01$    & -         & $3.73e-02$    & -         & $4.99e-02$    & -       \\ 
	$1$        & $1.32e-01$    & $3.03e-01$    & $1.98e-02$    & $2.74e-01$    & $2.62e-02$    & $2.80e-01$  \\ 
	$10^{-1}$  & $1.75e-02$    & $8.78e-01$    & $2.88e-03$    & $8.38e-01$    & $3.42e-03$    & $8.85e-01$  \\ 
	$10^{-2}$  & $1.78e-03$    & $9.94e-01$    & $2.95e-04$    & $9.89e-01$    & $3.44e-04$    & $9.97e-01$  \\ 
	$10^{-3}$  & $1.78e-04$    & $9.99e-01$    & $2.96e-05$    & $9.99e-01$    & $3.44e-05$    & $1.00e+00$  \\ 
	$10^{-4}$  & $1.78e-05$    & $1.00e+00$    & $2.96e-06$    & $1.00e+00$    & $3.44e-06$    & $1.00e+00$  \\ 
	$10^{-5}$  & $1.78e-06$    & $1.00e+00$    & $2.96e-07$    & $1.00e+00$    & $3.44e-07$    & $1.00e+00$  \\ 
	\hline
\end{tabular}
\noindent
\begin{tabular}{|C{20pt}||p{80pt}|C{40pt}||C{40pt}|C{40pt}|}
	\hline 
	\rule{0pt}{9pt}
	$K$ & $L^2(H^1(\Gamma))$ & $\EOC$ & $L^4(L^2(\Gamma))$  & $\EOC$ \\
	\hline 
	\rule{0pt}{9pt}
	$10$       & $1.60e-03$    & -         & $5.92e-04$    & -         \\ 
	$1$        & $6.24e-04$    & $4.08e-01$    & $1.77e-04$    & $5.25e-01$    \\ 
	$10^{-1}$  & $6.35e-05$ ($6.36e-05$)     & $9.92e-01$    & $1.48e-05$    & $1.08e+00$    \\ 
	$10^{-2}$  & $6.36e-06$    & $1.00e+00$    & $1.51e-06$    & $9.91e-01$    \\ 
	$10^{-3}$  & $6.36e-07$    & $1.00e+00$    & $1.52e-07$    & $9.99e-01$    \\ 
	$10^{-4}$  & $6.36e-08$ ($6.37e-08$)   & $1.00e+00$    & $1.52e-08$    & $1.00e+00$    \\ 
	$10^{-5}$  & $6.36e-09$ ($6.37e-09$)   & $1.00e+00$    & $1.52e-09$    & $1.00e+00$    \\ 
	\hline
\end{tabular}%\\[1ex]
\caption{Discrete error measured in various norms between solutions to  \eqref{lim:alt} and \eqref{CHLW} for $(\alpha, \beta) = (2,-4)$ and for $(\alpha, \beta) = (-2,4)$ with initial data (a).  Only the values for $(\alpha,\beta)=(2,-4)$ are listed.  Deviating values for $(\alpha,\beta)=(-2,4)$ are listed in brackets.}
\label{TBL:A:2}
\end{table}

%\normalsize\bigskip

%\noindent Simulation for initial datum (b) with $(\alpha,\beta)=(1,0)$:\\[1ex]
\begin{table}[h!]
\scriptsize
\begin{tabular}{|C{20pt}||C{40pt}|C{40pt}||C{40pt}|C{40pt}||C{40pt}|C{40pt}|}
	\hline
	\rule{0pt}{9pt}
	$K$ & $L^2(H^1(\Omega))$  & $\EOC$ & $L^4(L^2(\Omega))$ & $\EOC$ & $L^2(\Sigma_T)$ & $\EOC$ \\
	\hline
	\rule{0pt}{9pt}
	$10$       & $4.64e-01$    & -         & $8.77e-02$    & -         & $5.08e-02$    & -       \\ 
	$1$        & $2.56e-01$    & $2.58e-01$    & $5.23e-02$    & $2.24e-01$    & $1.85e-02$    & $4.39e-01$  \\ 
	$10^{-1}$  & $8.13e-02$    & $4.98e-01$    & $1.59e-02$    & $5.18e-01$    & $2.51e-03$    & $8.68e-01$  \\ 
	$10^{-2}$  & $1.07e-02$    & $8.83e-01$    & $2.13e-03$    & $8.72e-01$    & $2.67e-04$    & $9.73e-01$  \\ 
	$10^{-3}$  & $1.01e-03$    & $1.02e+00$    & $2.02e-04$    & $1.02e+00$    & $2.69e-05$    & $9.97e-01$  \\ 
	$10^{-4}$  & $1.01e-04$    & $1.00e+00$    & $2.01e-05$    & $1.00e+00$    & $2.69e-06$    & $1.00e+00$  \\ 
	$10^{-5}$  & $1.01e-05$    & $1.00e+00$    & $2.01e-06$    & $1.00e+00$    & $2.69e-07$    & $1.00e+00$  \\ 
	\hline
\end{tabular}
%
%\noindent
\begin{tabular}{|C{20pt}||C{40pt}|C{40pt}||C{40pt}|C{40pt}|}
	\hline 
	\rule{0pt}{9pt}
	$K$ & $L^2(H^1(\Gamma))$ & $\EOC$ & $L^4(L^2(\Gamma))$  & $\EOC$ \\
	\hline 
	\rule{0pt}{9pt}
	$10$       & $5.36e-01$    & -         & $7.82e-02$    & -         \\ 
	$1$        & $1.98e-01$    & $4.33e-01$    & $2.60e-02$    & $4.79e-01$    \\ 
	$10^{-1}$  & $5.84e-02$    & $5.29e-01$    & $8.31e-03$    & $4.95e-01$    \\ 
	$10^{-2}$  & $6.30e-03$    & $9.67e-01$    & $8.82e-04$    & $9.74e-01$    \\ 
	$10^{-3}$  & $6.18e-04$    & $1.01e+00$    & $8.63e-05$    & $1.01e+00$    \\ 
	$10^{-4}$  & $6.17e-05$    & $1.00e+00$    & $8.61e-06$    & $1.00e+00$    \\ 
	$10^{-5}$  & $6.17e-06$    & $1.00e+00$    & $8.61e-07$    & $1.00e+00$    \\ 
	\hline
\end{tabular}
\caption{Discrete error measured in various norms between solutions to \eqref{lim:alt} and \eqref{CHLW} for $(\alpha, \beta) = (1,0)$ with initial data (b).}
\label{TBL:B:1}
\end{table}
%\normalsize\bigskip

%\noindent Simulation for initial datum (b) with $(\alpha,\beta)=(2,-4)$:\\[1ex]
\begin{table}[h!]
\scriptsize
\begin{tabular}{|C{20pt}||C{40pt}|C{40pt}||C{40pt}|C{40pt}||C{40pt}|C{40pt}|}
	\hline
	\rule{0pt}{9pt}
	$K$ & $L^2(H^1(\Omega))$  & $\EOC$ & $L^4(L^2(\Omega))$ & $\EOC$ & $L^2(\Sigma_T)$ & $\EOC$ \\
	\hline
	\rule{0pt}{9pt}
	$10$       & $3.53e-01$    & -         & $6.74e-02$    & -         & $4.04e-02$    & -       \\ 
	$1$        & $1.48e-01$    & $3.77e-01$    & $2.48e-02$    & $4.34e-01$    & $1.90e-02$    & $3.26e-01$  \\ 
	$10^{-1}$  & $1.84e-02$    & $9.05e-01$    & $3.00e-03$    & $9.17e-01$    & $2.36e-03$    & $9.06e-01$  \\ 
	$10^{-2}$  & $1.88e-03$    & $9.91e-01$    & $3.06e-04$    & $9.92e-01$    & $2.37e-04$    & $9.98e-01$  \\ 
	$10^{-3}$  & $1.88e-04$    & $9.99e-01$    & $3.06e-05$    & $9.99e-01$    & $2.37e-05$    & $1.00e+00$  \\ 
	$10^{-4}$  & $1.88e-05$    & $1.00e+00$    & $3.06e-06$    & $1.00e+00$    & $2.37e-06$    & $1.00e+00$  \\ 
	$10^{-5}$  & $1.88e-06$    & $1.00e+00$    & $3.07e-07$    & $1.00e+00$    & $2.37e-07$    & $1.00e+00$  \\ 
	\hline
\end{tabular}
%\noindent
\begin{tabular}{|C{20pt}||C{40pt}|C{40pt}||C{40pt}|C{40pt}|}
	\hline 
	\rule{0pt}{9pt}
	$K$ & $L^2(H^1(\Gamma))$ & $\EOC$ & $L^4(L^2(\Gamma))$  & $\EOC$ \\
	\hline 
	\rule{0pt}{9pt}
	$10$       & $1.74e-03$    & -         & $4.45e-04$    & -         \\ 
	$1$        & $7.68e-04$    & $3.55e-01$    & $1.89e-04$    & $3.72e-01$    \\ 
	$10^{-1}$  & $9.92e-05$    & $8.89e-01$    & $2.64e-05$    & $8.55e-01$    \\ 
	$10^{-2}$  & $1.03e-05$    & $9.85e-01$    & $2.72e-06$    & $9.86e-01$    \\ 
	$10^{-3}$  & $1.03e-06$    & $9.98e-01$    & $2.73e-07$    & $9.99e-01$    \\ 
	$10^{-4}$  & $1.03e-07$    & $1.00e+00$    & $2.73e-08$    & $1.00e+00$    \\ 
	$10^{-5}$  & $1.03e-08$    & $1.00e+00$    & $2.73e-09$    & $1.00e+00$    \\ 
	\hline
\end{tabular}
\caption{Discrete error measured in various norms between solutions to \eqref{lim:alt} and \eqref{CHLW} for $(\alpha, \beta) = (2,-4)$ with initial data (b).}
\label{TBL:B:2}
\end{table}
%\normalsize\bigskip

%\newpage

%\noindent Simulation for initial datum (b) with $(\alpha,\beta)=(-5,30)$:\\[1ex]
\begin{table}[h!]
\scriptsize
\begin{tabular}{|C{20pt}||C{40pt}|C{40pt}||C{40pt}|C{40pt}||C{40pt}|C{40pt}|}
	\hline
	\rule{0pt}{9pt}
	$K$ & $L^2(H^1(\Omega))$  & $\EOC$ & $L^4(L^2(\Omega))$ & $\EOC$ & $L^2(\Sigma_T)$ & $\EOC$ \\
	\hline
	\rule{0pt}{9pt}
	$10$       & $3.75e-01$    & -         & $7.33e-02$    & -         & $4.05e-02$    & -       \\ 
	$1$        & $1.76e-01$    & $3.30e-01$    & $3.46e-02$    & $3.26e-01$    & $1.90e-02$    & $3.29e-01$  \\ 
	$10^{-1}$  & $2.07e-02$    & $9.27e-01$    & $3.45e-03$    & $1.00e+00$    & $2.35e-03$    & $9.08e-01$  \\ 
	$10^{-2}$  & $2.02e-03$    & $1.01e+00$    & $3.32e-04$    & $1.02e+00$    & $2.36e-04$    & $9.98e-01$  \\ 
	$10^{-3}$  & $2.02e-04$    & $1.00e+00$    & $3.31e-05$    & $1.00e+00$    & $2.36e-05$    & $1.00e+00$  \\ 
	$10^{-4}$  & $2.02e-05$    & $1.00e+00$    & $3.31e-06$    & $1.00e+00$    & $2.36e-06$    & $1.00e+00$  \\ 
	$10^{-5}$  & $2.02e-06$    & $1.00e+00$    & $3.31e-07$    & $1.00e+00$    & $2.36e-07$    & $1.00e+00$  \\ 
	\hline
\end{tabular}

%\noindent
\begin{tabular}{|C{20pt}||C{40pt}|C{40pt}||C{40pt}|C{40pt}|}
	\hline 
	\rule{0pt}{9pt}
	$K$ & $L^2(H^1(\Gamma))$ & $\EOC$ & $L^4(L^2(\Gamma))$  & $\EOC$ \\
	\hline 
	\rule{0pt}{9pt}
	$10$       & $4.71e-04$    & -         & $1.19e-04$    & -         \\ 
	$1$        & $2.18e-04$    & $3.34e-01$    & $5.56e-05$    & $3.30e-01$    \\ 
	$10^{-1}$  & $2.76e-05$    & $8.97e-01$    & $7.30e-06$    & $8.82e-01$    \\ 
	$10^{-2}$  & $2.82e-06$    & $9.92e-01$    & $7.36e-07$    & $9.97e-01$    \\ 
	$10^{-3}$  & $2.83e-07$    & $9.99e-01$    & $7.37e-08$    & $1.00e+00$    \\ 
	$10^{-4}$  & $2.83e-08$    & $1.00e+00$    & $7.37e-09$    & $1.00e+00$    \\ 
	$10^{-5}$  & $2.90e-09$    & $9.89e-01$    & $7.37e-10$    & $1.00e+00$    \\ 
	\hline
\end{tabular}
\caption{Discrete error measured in various norms between solutions to \eqref{lim:alt} and \eqref{CHLW} for $(\alpha, \beta) = (-5,30)$ with initial data (b).}
\label{TBL:B:3}
\end{table}
\normalsize\bigskip

\FloatBarrier

\begin{figure}[h!]
	\centering
	\captionsetup[subfigure]{position=top,labelformat=empty}
	%\textbf{Simulations for initial data (a) and (b) with a nonlinear relation $H$:}\\[-1ex]
	%%%%%%%%%%%%%%%%%%%%%%%%%%%%%%%%%%%%
	%%% Simulation Wave H = sin:
	%%%%%%%%%%%%%%%%%%%%%%%%%%%%%%%%%%%%
	\subfloat[$K=10$]{\includegraphics[width=0.24\textwidth]{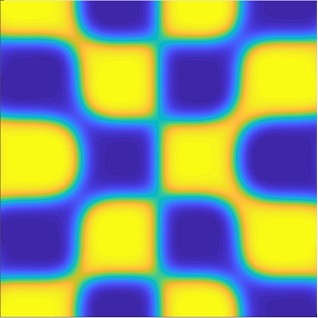}}
	\hfill
	\subfloat[$K=1$]{\includegraphics[width=0.24\textwidth]{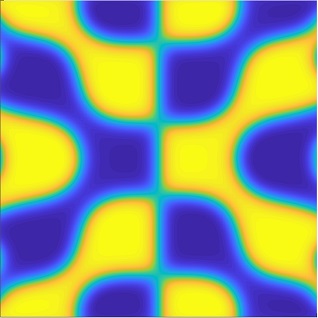}}
	\hfill
	\subfloat[$K=0.1$]{\includegraphics[width=0.24\textwidth]{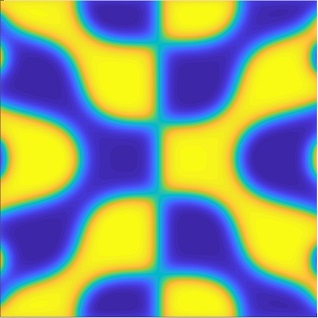}}
	\hfill
	\subfloat[$K=10^{-5}$]{\includegraphics[width=0.24\textwidth]{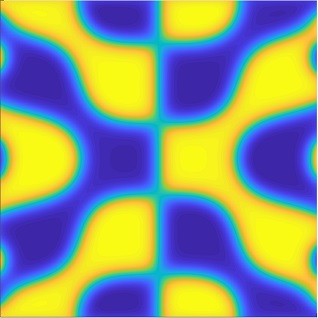}}\\[0.5ex]
%\end{figure}
%\medskip
%\begin{figure}[h!]
%	\centering
%	\vspace{-15pt}
	%%%%%%%%%%%%%%%%%%%%%%%%%%%%%%%%%%%%
	%%% Simulation Random H = ?
	%%%%%%%%%%%%%%%%%%%%%%%%%%%%%%%%%%%%
	\includegraphics[width=0.24\textwidth]{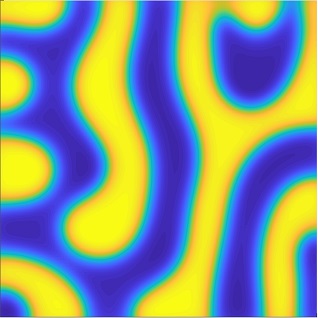}
	\hfill
	\includegraphics[width=0.24\textwidth]{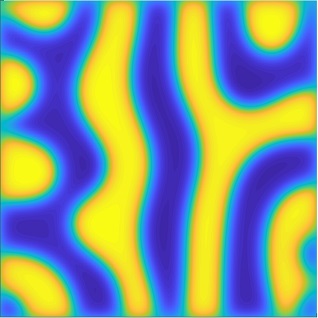}
	\hfill
	\includegraphics[width=0.24\textwidth]{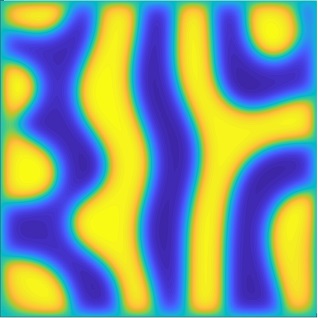}
	\hfill
	\includegraphics[width=0.24\textwidth]{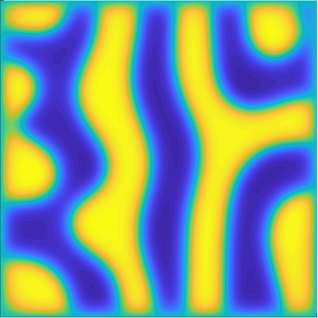}\\
	\vspace{8pt}
	\includegraphics[width=0.42\textwidth]{cb_hor}
	\caption{Numerical solution $u_h^k$ at iteration $k=100$. 
		First row: Simulation for initial datum (a) with $\kappa=1$ and $H(s)=\sin(s)$.
		Second row: Simulation for initial datum (b) with $\kappa=0.4$ and $H(s)=3\cos(s)+2$.}\label{FIG:SIM:NONLIN}
\end{figure}

\paragraph{Simulations for nonlinear transmission relations.} In the last set of experiments we consider nonlinear relations $H$ for the system \eqref{CHLW}. We additionally assume that there exists an interval $I\subseteq\R$ such that $H\vert_I:I\to [-1,1]$ is a bijection. For the simulations (see Figure \ref{FIG:SIM:NONLIN}) we use once more the initial data (a) and (b) (as defined above) where we set $v_h^0 = H\vert_I^{-1}(u_h^0)$. As we did not derive a limit system in the general nonlinear case, we use the solution to system \eqref{CHLW} with $K=10^{-5}$ as the reference solution instead.  The discrete error measured in various norms (as described above) between $(u^K , v^K)$ and $(u^{10^{-5}},  v^{10^{-5}})$ are displayed in Tables \ref{TBL:C:1} and \ref{TBL:C:2}, where we again observe a linear rate of convergence.

\normalsize\bigskip

%\noindent Simulation for initial datum (a) with $H=\sin(\cdot)$:\\[1ex]
\begin{table}[h!]
\scriptsize
\begin{tabular}{|C{20pt}||C{40pt}|C{40pt}||C{40pt}|C{40pt}||C{40pt}|C{40pt}|}
	\hline
	\rule{0pt}{9pt}
	$K$ & $L^2(H^1(\Omega))$  & $\EOC$ & $L^4(L^2(\Omega))$ & $\EOC$ & $L^2(\Sigma_T)$ & $\EOC$ \\
	\hline
	\rule{0pt}{9pt}
	$10$       & $2.52e-01$    & -         & $3.86e-02$    & -         & $7.83e-03$    & -       \\ 
	$1$        & $1.14e-01$    & $3.46e-01$    & $1.79e-02$    & $3.33e-01$    & $2.68e-03$    & $4.66e-01$  \\ 
	$10^{-1}$  & $1.40e-02$    & $9.09e-01$    & $2.26e-03$    & $8.99e-01$    & $3.42e-04$    & $8.93e-01$  \\ 
	$10^{-2}$  & $1.43e-03$    & $9.89e-01$    & $2.31e-04$    & $9.91e-01$    & $3.60e-05$    & $9.78e-01$  \\ 
	$10^{-3}$  & $1.43e-04$    & $1.00e+00$    & $2.29e-05$    & $1.00e+00$    & $3.59e-06$    & $1.00e+00$  \\ 
	$10^{-4}$  & $1.30e-05$    & $1.04e+00$    & $2.08e-06$    & $1.04e+00$    & $3.27e-07$    & $1.04e+00$  \\ 
	$10^{-5}$  & $0.00e+00$    & -    & $0.00e+00$    & -    & $0.00e+00$    & -  \\ 
	\hline
\end{tabular}

%\noindent
\begin{tabular}{|C{20pt}||C{40pt}|C{40pt}||C{40pt}|C{40pt}|}
	\hline 
	\rule{0pt}{9pt}
	$K$ & $L^2(H^1(\Gamma))$ & $\EOC$ & $L^4(L^2(\Gamma))$  & $\EOC$ \\
	\hline 
	\rule{0pt}{9pt}
	$10$       & $2.72e-01$    & -         & $6.22e-02$    & -         \\ 
	$1$        & $9.59e-02$    & $4.53e-01$    & $2.49e-02$    & $3.98e-01$    \\ 
	$10^{-1}$  & $1.23e-02$    & $8.91e-01$    & $3.33e-03$    & $8.74e-01$    \\ 
	$10^{-2}$  & $1.30e-03$    & $9.78e-01$    & $3.51e-04$    & $9.77e-01$    \\ 
	$10^{-3}$  & $1.29e-04$    & $1.00e+00$    & $3.50e-05$    & $1.00e+00$    \\ 
	$10^{-4}$  & $1.18e-05$    & $1.04e+00$    & $3.18e-06$    & $1.04e+00$    \\ 
	$10^{-5}$  & $0.00e+00$    & -    & $0.00e+00$    & -   \\ 
	\hline
\end{tabular}
\caption{Discrete error measured in various norms between $(u^K, v^K)$ and $(u^{10^{-5}}, v^{10^{-5}})$ for $H(s) = \sin(s)$.}
\label{TBL:C:1}
\end{table}
%\normalsize\bigskip

%\noindent Simulation for initial datum (b) with $H=3\cos(\cdot)+2$:\\[1ex]
\begin{table}[h!]
\scriptsize
\begin{tabular}{|C{20pt}||C{40pt}|C{40pt}||C{40pt}|C{40pt}||C{40pt}|C{40pt}|}
	\hline
	\rule{0pt}{9pt}
	$K$ & $L^2(H^1(\Omega))$  & $\EOC$ & $L^4(L^2(\Omega))$ & $\EOC$ & $L^2(\Sigma_T)$ & $\EOC$ \\
	\hline
	\rule{0pt}{9pt}
	$10$       & $3.54e-01$    & -         & $6.75e-02$    & -         & $6.18e-05$    & -       \\ 
	$1$        & $1.48e-01$    & $3.81e-01$    & $2.47e-02$    & $4.36e-01$    & $2.57e-05$    & $3.81e-01$  \\ 
	$10^{-1}$  & $1.86e-02$    & $9.00e-01$    & $3.03e-03$    & $9.12e-01$    & $3.27e-06$    & $8.95e-01$  \\ 
	$10^{-2}$  & $1.90e-03$    & $9.90e-01$    & $3.10e-04$    & $9.91e-01$    & $3.37e-07$    & $9.87e-01$  \\ 
	$10^{-3}$  & $1.89e-04$    & $1.00e+00$    & $3.08e-05$    & $1.00e+00$    & $3.34e-08$    & $1.00e+00$  \\ 
	$10^{-4}$  & $1.72e-05$    & $1.04e+00$    & $2.80e-06$    & $1.04e+00$    & $3.04e-09$    & $1.04e+00$  \\ 
	$10^{-5}$  & $0.00e+00$    & -   & $0.00e+00$    & -    & $0.00e+00$    & -  \\ 
	\hline
\end{tabular}
%\noindent
\begin{tabular}{|C{20pt}||C{40pt}|C{40pt}||C{40pt}|C{40pt}|}
	\hline 
	\rule{0pt}{9pt}
	$K$ & $L^2(H^1(\Gamma))$ & $\EOC$ & $L^4(L^2(\Gamma))$  & $\EOC$ \\
	\hline 
	\rule{0pt}{9pt}
	$10$       & $1.44e-03$    & -         & $3.69e-04$    & -         \\ 
	$1$        & $6.40e-04$    & $3.53e-01$    & $1.58e-04$    & $3.69e-01$    \\ 
	$10^{-1}$  & $8.29e-05$    & $8.88e-01$    & $2.20e-05$    & $8.55e-01$    \\ 
	$10^{-2}$  & $8.58e-06$    & $9.85e-01$    & $2.28e-06$    & $9.86e-01$    \\ 
	$10^{-3}$  & $8.53e-07$    & $1.00e+00$    & $2.26e-07$    & $1.00e+00$    \\ 
	$10^{-4}$  & $7.76e-08$    & $1.04e+00$    & $2.06e-08$    & $1.04e+00$    \\ 
	$10^{-5}$  & $0.00e+00$    & -    & $0.00e+00$    & -   \\ 
	\hline
\end{tabular}
\caption{Discrete errors measured in various norms between $(u^K, v^K)$ and $(u^{10^{-5}}, v^{10^{-5}})$ for $H(s) = 3 \cos(s) + 2$.}
\label{TBL:C:2}
\end{table}

\FloatBarrier

\section{Conclusion}
In this work, we derived and analysed an extension of a recent Cahn--Hilliard model proposed by C.~Liu and H.~Wu with a new type of dynamic boundary conditions.  The extension replaces the classical transmission condition with an affine linear relation between the bulk and surface order parameters.  A natural time discretisation inspired by the associated gradient flow structure allow us to establish the existence of unique weak solutions, and numerical simulations demonstrate that under certain choices of potentials and initial conditions, there seems to be a competition between mass conservation and Cahn--Hilliard phase separation dynamics occurring on the domain boundary.  This behaviour is not observed with the original model of Liu and Wu, and at present we have yet to find a criterion based on the initial data and parameter values that would allow us to predict which of the competing effects is dominant.

As a first step towards nonlinear relations in the dynamic boundary condition, we study a Robin boundary type regularisation.  Similarly, the  regularised system exhibits a gradient flow structure allows us to exploit a natural time discretisation to establish the existence of unique weak solutions.  In comparison, for the case where there is no diffusion for the surface phase field variable, we have to place more restrictions on the surface potential to obtain analogous results.  Meanwhile, for the case of affine linear relations, we can show the weak convergence of solutions to the Robin problem as the regularisation parameter tends to zero, and recover the unique solution to the extended Liu--Wu model.  An error estimate is also derived and is supported by our numerical findings.  Furthermore, our simulations also suggest that similar rates of convergence hold for nonlinear relations in the dynamic boundary condition, and we leave the analytical investigation to future research.

\paragraph{Acknowledgements.}
The work of the second author was supported by a Direct Grant of CUHK, Hong Kong (project 4053288). The kind support of the Fudan Scholar program from Fudan University is also acknowledged.  The authors would like to thank the anonymous referees for their careful reading and useful suggestions which have improved the quality of the manuscript.
\appendix

\section{Energetical variational approach for derivation of \eqref{lim:alt}}
The derivation of \eqref{lim:alt} using the energetically variational approach proceeds in three steps.

\paragraph{Mass conservation.} In the bulk domain $\Omega$, we assume there is a microscopic effective velocity $\bm{v} : \Omega \to \R^3$ satisfying the no-flux boundary condition $\bm{v} \cdot \bm{n} = 0$ on $\Gamma = \pd \Omega$, and $u$ is a locally conserved (sufficiently regular) quantity satisfying the continuity equation
\begin{align}\label{EVA:1}
u_t + \div (u \bm{v}) = 0 \text{ in } Q_T.
\end{align}
Furthermore, we assume that $\bm{w} : \Gamma \to \R^3$ is a microscopic effective tangential velocity field on $\Gamma$ and the boundary dynamics are characterised by a local mass conservation law for the trace of $u$, which we denote again as $u$:
\begin{align}\label{EVA:2}
u_t + \div_\Gamma(u \bm{w}) = 0 \text{ on } \Sigma_T.
\end{align}
Then it is clear that $\frac{d}{dt} \mean{u}_\Omega = 0$ and $\frac{d}{dt} \mean{u}_\Gamma = 0$.

\paragraph{Energy dissipation.} The basic energy law considered is
\begin{align*}
\frac{d}{dt} \mathcal{E}(t) = - \mathcal{D}(t),
\end{align*}
where $\mathcal{E}(t) = E_*(u(t))$ and $\mathcal{D}(t)$ is the rate of energy dissipation defined as
\begin{align*}
\mathcal{D}(t) = \int_\Omega u^2 \abs{\bm{v}}^2 \dx + \int_\Gamma \frac{u^2}{\alpha^2} \abs{\bm{w}}^2 \dG.
\end{align*}

\paragraph{Force balance.} It remains to determine the microscopic velocities $\bm{v}$ and $\bm{w}$ in the continuity equations \eqref{EVA:1} and \eqref{EVA:2} via the energetic variational approach.  We denote by $\Omega_0^X$ the reference configuration with Lagrangian coordinate $X$, and $\Omega_x^t$ the deformed configuration with Eulerian coordinate $x$, which are related via the bulk flow map $x(X,t): \Omega_0^X \to \Omega_x^t$ defined as a solution to the ordinary differential system
\begin{align*}
\frac{d}{dt} x(X,t) = \bm{v}(x(X,t),t), \quad x(X,0) = X
\end{align*} 
induced by the bulk microscopic velocity $\bm{v}$.  Similarly, we introduce the reference boundary $\Gamma_0^X$, the deformed boundary $\Gamma_x^t$ and the surface flow map $x_s = x_s(X_s,t)$ induced by $\bm{w}$ with Lagrangian and Eulerian surface coordinates $X_s$ and $x_s$, respectively.  We mention that the surface flow map is the solution to an analogous ordinary differential system with the microscopic velocity $\bm{w}$, and is not related to the bulk flow map. Now, let the total action functional $\mathcal{A}$ be defined as the negative integral of $\mathcal{E}(t)$ over $[0,T]$. By computations similar to \cite[Sec.~7.1]{LW}, the variation of $\mathcal{A}$ with respect to the spatial variables $x$ in $\Omega$ and $x_s$ on $\Gamma$ yields
\begin{align*}
\delta_{(x,x_s)} \mathcal{A} = - \int_0^T \int_{\Omega_t^x} u \nabla \mu \cdot \tilde{\bm{y}} \dx \dt - \int_0^T \int_{\Gamma_t^x} u \gradg (\mu_s + \eps \pdnu u) \cdot \tilde{\bm{y}}_s \dG \dt
\end{align*}
for arbitrary smooth vector fields $\tilde{\bm{y}} = \tilde{\bm{y}}(x,t)$ and $\tilde{\bm{y}}_s = \tilde{\bm{y}}_s(x_s, t)$ satisfying $\tilde{\bm{y}} \cdot \bn = \tilde{\bm{y}}_s \cdot \bn = 0$ on $\Gamma_t^x$. Here, $\mu$ and $\mu_s$ are defined as
\begin{align*}
\mu = - \eps \Lx u + \eps^{-1} F'(u), \quad \mu_s = - \frac{\delta \kappa}{\alpha^2} \LB u + \frac{1}{ \delta \alpha} G'(\alpha^{-1}(u-\beta)).
\end{align*}
Invoking the principle of least action, the conservative force balance equations are identified by solving $\delta \mathcal{A} = 0$, which yields the bulk and surface conservative forces:
\begin{align*}
\bm{f}_{con}^b = - u \nabla \mu, \quad \bm{f}_{con}^{s} = - u \gradg (\mu_s + \pdnu u).
\end{align*}
Furthermore, via Onsager's maximum dissipation principle for irreversible dissipative processes in the regime of linear response, taking the variation $\delta_{(\bm{v}, \bm{w})} \mathcal{R}$ of the Rayleigh dissipation function $\mathcal{R} = \frac{1}{2} \mathcal{D}$ with respect to the rate functions (i.e., the microscopic velocities $\bm{v}$ and $\bm{w}$) gives the dissipative force balance equations. In our case, they read as
\begin{align*}
\delta_{(\bm{v}, \bm{w})} \tfrac{1}{2} \mathcal{D} = \bm{0} \quad \Rightarrow \quad \bm{f}_{diss}^b = - u^2 \bm{v}, \quad \bm{f}_{diss}^s = - \alpha^{-2} u^2 \bm{w}.
\end{align*}
The classical Newton's force balance law $\bm{f}_{con} + \bm{f}_{diss} = \bm{0}$ leads to the force balance relations in the bulk $\Omega$ and on the boundary $\Gamma$, which additionally allows us to infer an expression for the microscopic velocities:
\begin{align*}
u \nabla \mu + u^2 \bm{v} = \bm{0} \text{ in } \Omega_t^x, \quad u \gradg (\mu_s + \eps \pdnu u) + \alpha^{-2} u^2 \bm{w} = \bm{0} \text{ on } \Gamma_t^x.
\end{align*}
Eventually, in combination with \eqref{EVA:1} and \eqref{EVA:2} we conclude that
\begin{align*}
&\begin{aligned}
&u_t = \Lx \mu, \quad &&\mu = - \eps \kappa \Lx u + \eps^{-1} F'(u) & \text{ in } \Omega, \\
&u_t = \LB \phi, \quad &&\phi = \alpha^{2}(\mu_s + \eps \pdnu u) = - \delta \kappa \LB u + \alpha \delta^{-1} G'(\alpha^{-1}(u-\beta)) + \eps \alpha^2 \pdnu u & \text{ on } \Gamma, 
\end{aligned}\\
&\pdnu \mu = u \bm{v} \cdot \bn = 0 \quad \text{ on } \Gamma,
\end{align*}
which coincides with the system \eqref{lim:alt}.

\pagebreak[2]

\footnotesize
%\nocite{*} % Show all bib entries - both cited and uncited; comment this line to view only cited bib entries;
\bibliographystyle{plain}
%\bibliography{CH_LiuWu_Robin}%

\begin{thebibliography}{999}

\bibitem{Alt} H.W.~Alt.  {\it Linear Functional Analysis - An Application-Oriented Introduction}.  Springer-Verlag London, 2016.

\bibitem{Buffa} A.~Buffa, M.~Costabel and D.~Sheen.  On traces for H(curl,$\Omega$) in Lipschitz domains.  {\it J. Math. Anal. Appl.} {\bf 276} (2002) 845--867

\bibitem{CC} L.~Calatroni and P.~Colli.  Global solution to the Allen-Cahn equation with singular potentials and dynamic boundary conditions.  {\it Nonlinear Anal.} {\bf 79} (2013) 12--27

\bibitem{CherfilsGM} L.~Cherfils, S.~Gatti and A.~Miranville.  A variational approach to a Cahn--Hilliard model in a domain with nonpermeable walls.  {\it J. Math. Sci.} {\bf 189} (2013) 604--636

\bibitem{CP} L.~Cherfils and M.~Petcu.  A numerical analysis of the Cahn--Hilliard equation with non-permeable walls.  {\it Numer. Math.} {\bf 128} (2014) 517--549

\bibitem{CPP} L.~Cherfils, M.~Petcu and M.~Pierre.  A numerical analysis of the Cahn--Hilliard equation with dynamic boundary conditions. {\it Discrete Contin. Dyn. Syst.} {\bf 27} (2010) 1511--1533

\bibitem{ColliF} P.~Colli and T.~Fukao.  Equation and dynamic boundary condition of Cahn--Hilliard type with singular potentials.  {\it Nonlinear Anal.} {\bf 127} (2015) 413--433

\bibitem{CFL} P.~Colli, T.~Fukao and K.F.~Lam.  On a coupled bulk-surface Allen--Cahn system with an affine linear transmission condition and its approximation by a Robin boundary condition.  {\it Nonlinear Anal.} {\bf 184} (2019) 116--147

\bibitem{CFW} P.~Colli, T.~Fukao and H.~Wu.  On a transmission problem for equation and dynamic boundary condition of Cahn--Hilliard type with nonsmooth potentials.  arXiv:1907.13279 [math.AP] (2019)

\bibitem{Evans} L.C.~Evans.  Partial Differential Equations.  AMS, Providence, R.I. 2010

\bibitem{FP} X.~Feng and A.~Prohl.  Error analysis of a mixed finite element method for the Cahn--Hilliard equation.  {\it Numer. Math.} {\bf 99} (2004) 47--84

\bibitem{FP2} X.~Feng and A.~Prohl.  Numerical analysis of the Cahn--Hilliard equation and approximation for the Hele-Shaw problem.  {\it Interfaces Free Bound.} {\bf 7} (2005) 1--28

\bibitem{Fis1} H.P.~Fischer, P.~Maass and W.~Dieterich.  Novel surface modes in spinodal decomposition.  {\it Phys. Rev. Lett.} {\bf 79} (1997) 893--896

\bibitem{Fis2} H.P.~Fischer, J.~Reinhard, W.~Dieterich, J.-F.~Gouyet, P.~Maass, A.~Majhofer and D.~Reinel.  Time-dependent density functional theory and the kinetics of lattice gas systems in contact with a wall.  {\it J. Chem. Phys.} {\bf 108} (1998) 3028--3037

\bibitem{FYW} T.~Fukao, S.~Yoshikawa and S.~Wada.  Structure-preserving finite difference schemes for the Cahn--Hilliard equation with dynamic boundary conditions in the one-dimensional case.  {\it Commun. Pure Applied Anal.} {\bf 16} (2017) 1915--1938

\bibitem{GGM} C.G.~Gal, M.~Grasselli and A.~Miranville.  Cahn--Hilliard--Navier--Stokes systems with moving contact lines. {\it Calc. Var. Partial Differ. Equ.} {\bf 55} (2016) article 50

\bibitem{GGW} C.G.~Gal, M.~Grasselli and H.~Wu.  Global weak solutions to a diffuse interface model for incompressible two-phase flows with moving contact lines and different densities.  {\it Arch. Rational Mech. Anal.} {\bf 234} (2019) 1--56

\bibitem{garckeelas} H. Garcke, On Cahn--Hilliard systems with elasticity. {\it Proc. Roy. Soc. Edinburgh} {\bf 133 A} (2003), pp. 307--331.

\bibitem{GK} H.~Garcke and P.~Knopf.  Weak solutions of the Cahn--Hilliard system with dynamic boundary conditions: A gradient flow approach.  {\it SIAM J. Math. Anal.} {\bf 52} (2020) 340--369

\bibitem{GMS} G.R.~Goldstein, A.~Miranville and G.~Schimperna.  A Cahn--Hilliard model in a domain with non-permeable walls.  {\it Phys. D} {\bf 240} (2011) 754--766

\bibitem{GT} F.~Guill\'en-Gonz\'alez and G.~Tierra.  On linear schemes for a Cahn--Hilliard diffuse interface model.  {\it J. Comput. Phys.} {\bf 243} (2013) 140--171

\bibitem{IMP} H.~Israel, A.~Miranville and M.~Petcu.  Numerical analysis of a Cahn--Hilliard type equation with dynamic boundary conditions.  {\it Ricerche Mat.} {\bf 64} (2015) 25--50

\bibitem{Kenzler} R.~Kenzler, F.~Eurich, P.~Maass, B.~Rinn, J.~Schropp, E.~Bohl and W.~Dieterich.  Phase separation in confined geometries: solving the Cahn--Hilliard equation with generic boundary conditions.  {\it Comput. Phys. Commun.} {\bf 133} (2001) 139--157

\bibitem{KL} B.~Kov\'acs and C.~Lubich.  Numerical analysis of parabolic problems with dynamic boundary conditions. {\it IAM J. Numer. Anal.} {\bf 37} (2017) 1--39

\bibitem{LW} C.~Liu and H.~Wu. An Energetic Variational Approach for the Cahn--Hilliard Equation with Dynamic Boundary Conditions: Model Derivation and Mathematical Analysis.  {\it Arch. Rational Mech. Anal.} {\bf 233} (2019) 167--247

\bibitem{Metzger} S.~Metzger.  An efficient and convergent finite element scheme for Cahn--Hilliard equations with dynamic boundary conditions.  arXiv:1908.04910 [math.NA]  (2019)

\bibitem{Miran} A.~Miranville.  The Cahn--Hilliard Equation: Recent Advances and Applications.  CBMS-NSF Regional Conference Series in Applied Mathematics.  SIAM 2019

\bibitem{Qian} T.Z.~Qian, X.P.~Wang and P.~Sheng.  A variational approach to moving contact line hydrodynamics.  {\it J. Fluid Mech.} {\bf 564} (2006) 333--360

\bibitem{SXY} J.~Shen, J.~Xu and J.~Yang.  The scalar auxiliary variable (SAV) approach for gradient flows.  {\it J. Comput. Phys.} {\bf 353} (2018) 407--416

\bibitem{SY} J.~Shen and X.~Yang.  Numerical approximation of Allen-Cahn and Cahn-Hilliard equations.  {\it Discrete Contin. Dyn. Syst.} {\bf 28} (2010) 1669--1691

\bibitem{Trau} D.~Trautwein.  Finite-Elemente Approximation der Cahn-Hilliard-Gleichung mit Neumann- und dynamischen Randbedingungen.  Bachelor thesis, University of Regensburg (2018)

\end{thebibliography}

\end{document}